\documentclass[12pt]{article}
\usepackage{latexsym, amssymb, amsmath, amscd, amsfonts, epsfig, graphicx, colordvi,verbatim,ifpdf,extarrows}
\usepackage{amsfonts, amsmath, amssymb}
\usepackage{amssymb,amsfonts,amsmath,latexsym,epsfig,cite, psfrag,eepic,color}
\usepackage{amscd,graphics}
\usepackage{latexsym, amssymb,  amsmath,amscd, amsfonts, epsfig, graphicx, colordvi,amsthm}

\usepackage{graphicx}
\usepackage{epstopdf}
\usepackage{color}
\usepackage{ifpdf}
\usepackage{fancybox}
\usepackage[font=small,labelfont=bf,labelsep=none]{caption}
\usepackage{float}

\newtheorem{thm}{Theorem}[section]
\newtheorem{prop}[thm]{Proposition}
\newtheorem{conj}[thm]{Conjecture}

\newtheorem{defi}[thm]{Definition}
\newtheorem{lem}[thm]{Lemma}

\def\pf{\noindent{\it Proof.} }
\def\qed{\nopagebreak\hfill{\rule{4pt}{7pt}}
\medbreak}

\numberwithin{equation}{section}

\def\qed{\nopagebreak\hfill{\rule{4pt}{7pt}}
\medbreak}

\newlength{\boxedparwidth}
  {\begin{center} \begin{tabular}{|@{\hspace{.315in}}c@{\hspace{.15in}}|}
                  \hline \\ \begin{minipage}[t]{\boxedparwidth}
                  \setlength{\parindent}{.25in}}%
  {\end{minipage} \\ \\ \hline \end{tabular} \end{center}}

\begin{document}

\begin{center}

 {\Large \bf Overpartitions and Bressoud's Conjecture, I}
\end{center}

\begin{center}
 {Thomas Y. He}$^{1}$, {Kathy Q. Ji}$^{2}$  and
  {Alice X.H. Zhao}$^{3}$ \vskip 2mm

$^{1}$ School of Mathematical Sciences, Laurent Mathematics Center and V.C. \&  V.R. Key Lab, Sichuan Normal University, Chengdu 610066, P.R. China\\[6pt]

$^{2}$ Center for Applied Mathematics,  Tianjin University, Tianjin 300072, P.R. China\\[6pt]

  $^{3}$ College of Science, Tianjin University of Technology, Tianjin 300384, P.R. China

   \vskip 2mm

    $^1$heyao@sicnu.edu.cn, $^2$kathyji@tju.edu.cn,  $^3$zhaoxh@email.tjut.edu.cn
\end{center}

\vskip 6mm   {\bf Abstract.}  In  1980,  Bressoud conjectured a combinatorial identity $A_j=B_j$ for $j=0$ or $1$,  where  the function $A_j$ counts the number of partitions with certain   congruence conditions and the function $B_j$ counts the number of partitions with certain difference conditions. Bressoud's conjecture  specializes to a wide variety of well-known theorems in the theory of partitions. Special cases of his conjecture have been subsequently proved by Bressoud, Andrews, Kim and Yee. Recently, Kim resolved  Bressoud's conjecture  for the case $j=1$. In this paper,  we introduce  a  new partition function $\overline{B}_j$ which can be viewed as an overpartition analogue of the partition function $B_j$ introduced by Bressoud. By means of
Gordon markings, we build bijections to  obtain a relationship between $\overline{B}_1$ and $B_0$ and a relationship between $\overline{B}_0$ and $B_1$. Based on these former relationships, we further give  overpartition analogues of many classical partition theorems including Euler's partition theorem, the Rogers-Ramanujan-Gordon identities, the Bressoud-Rogers-Ramanujan identities,  the Andrews-G\"ollnitz-Gordon identities and the Bressoud-G\"ollnitz-Gordon identities.

\noindent {\bf Keywords}: Bressoud's conjecture, Overpartitions, Euler's partition theorem, Rogers-Ramanujan identities, G\"ollnitz-Gordon identities,  Bailey pairs, Gordon markings

\noindent {\bf AMS Classifications}: 05A17,11P81,11P84

\section{Introduction}
 Bressoud \cite{Bressoud-1980} proved and conjectured some partition identities involving the partition function $B_j$, which counts the number of partitions with certain difference conditions (see Definition \ref{Bressoud-b-defi-0000000001}). The main objective of this paper is to introduce  a  new partition function $\overline{B}_j$
 which can be regarded as an overpartition analogue of the partition function $B_j$. We establish a relationship between $\overline{B}_1$ and
 $B_0$ and a relationship between $\overline{B}_0$ and $B_1$. Based on these two relationships,
 we obtain overpartition analogues of many classical partition theorems including Euler's partition theorem, the Rogers-Ramanujan-Gordon identities, the Bressoud-Rogers-Ramanujan identities,  the Andrews-G\"ollnitz-Gordon identities and the Bressoud-G\"ollnitz-Gordon identities. It should be noted that the relationship between  $\overline{B}_1$ and
 $B_0$ plays a crucial role in the proof of Bressoud's conjecture for $j=0$ in the subsequent
 paper \cite{He-Ji-Zhao}.

Let us  recall some common notation and terminologies on partitions from \cite[Chapter 1]{Andrews-1976}. A partition $\pi$ of a positive integer $n$ is a finite non-increasing sequence of positive integers $\pi=(\pi_1,\pi_2,\ldots,\pi_\ell)$
such that $\sum_{i=1}^{\ell} \pi_i=n$. An overpartition  of  $n$ is a partition of $n$ such that the first occurrence of a part can be overlined.

For example, there are five partitions of $4$:
\[(4),(3,1),(2,2),(2,1,1),(1,1,1,1),\]
whereas there are  fourteen overpartitions of $4$:
\[(4),(\bar{4}),(3,1),(\bar{3},1),(3,\bar{1}),(\bar{3},\bar{1}),(2,2),(\bar{2},2),
\]
\[(2,1,1),(\bar{2},1,1),(2,\bar{1},1),(\bar{2},\bar{1},1),(1,1,1,1),(\bar{1},1,1,1).\]

We impose the following order on the parts of an overpartition:
\begin{equation}\label{order}
1<\bar{1}<2<\bar{2}<\cdots.
\end{equation}

Let $\pi=(\pi_1,\pi_2,\ldots,\pi_\ell)$ be an ordinary partition (resp. an overpartition) with
$\pi_1\geq\pi_2\geq\cdots\geq \pi_\ell\geq 1$.
 The number of parts of $\pi$ is called the length of $\pi$, denoted $\ell(\pi)$. The weight of $\pi$ is the sum of parts, denoted $|\pi|$.

In 1961, Gordon \cite[p.~394]{Gordon-1961}   found  an infinite family of  combinatorial generalizations of the Rogers-Ramanujan identities, which has been known as the Rogers-Ramanujan-Gordon theorem.

\begin{thm}[Rogers-Ramanujan-Gordon]\label{R-R-B-1}
For $k\geq r\geq 1$, let ${B}_1(-;1,k,r;n)$ denote the number of partitions $\pi=(\pi_1,\pi_2,\ldots,\pi_\ell)$  of $n$, where $\pi_i\geq\pi_{i+k-1}+2$ for $1\leq i\leq \ell-k+1$, and at most $r-1$ of the $\pi_i$ are equal to $1$. For $k\geq r\geq 1$, let ${A}_1(-;1,k,r;n)$ denote the number of partitions of $n$ into parts $\not\equiv0,\pm r\pmod{2k+1}$.
 Then, for   $k\geq r\geq 1$   and $n\geq0$,
\begin{equation*}\label{R-R-B-1-e}
{A}_1(-;1,k,r;n)={B}_1(-;1,k,r;n).
\end{equation*}

\end{thm}

An analytic proof  of Theorem \ref{R-R-B-1}   was  given by Andrews \cite{Andrews-1974}.  He discovered  the following generating function version of Theorem \ref{R-R-B-1}, which has been called the Andrews-Gordon identity: For $k\geq r\geq 1$,
\begin{equation}\label{R-R-A1}
\sum_{N_1\geq  \cdots \geq N_{k-1}\geq0}\frac{q^{N_1^2+\cdots+N^2_{k-1}+N_r+\cdots+N_{k-1}}}
{(q;q)_{N_1-N_2}\cdots(q;q)_{N_{k-2}-N_{k-1}}(q;q)_{N_{k-1}}}=\frac{(q^r,q^{2k-r+1},q^{2k+1};q^{2k+1})_\infty}{(q;q)_\infty}.
\end{equation}

From now on, we assume that $|q|<1$ and adopt the standard notation \cite{Andrews-1976}:
\[(a;q)_\infty=\prod_{i=0}^{\infty}(1-aq^i), \quad (a;q)_n=\frac{(a;q)_\infty}{(aq^n;q)_\infty},\]
and
\[(a_1,a_2,\ldots,a_m;q)_\infty=(a_1;q)_\infty(a_2;q)_\infty\cdots(a_m;q)_\infty.
\]

In 1979, Bressoud \cite{Bressoud-1979} extended the Rogers-Ramanujan-Gordon theorem to even moduli, which  has been called the Bressoud-Rogers-Ramanujan theorem.

\begin{thm}[Bressoud-Rogers-Ramanujan] \label{R-R-B}
For $k> r\geq 1$, let ${B}_0(-;1,k,r;n)$ denote the number of partitions $\pi=(\pi_1,\pi_2,\ldots,\pi_\ell)$
of $n$, where $\pi_i\geq\pi_{i+k-1}+2$ for $1\leq i\leq \ell-k+1$, at most $r-1$ of the $\pi_i$ are
equal to $1$, and for $1\leq i\leq \ell-k+2$, if $\pi_i\leq\pi_{i+k-2}+1$, then
\[\pi_i+\cdots+\pi_{i+k-2}\equiv r-1\pmod{2}.\]
For $k> r\geq 1$, let ${A}_0(-;1,k,r;n)$ denote the number of partitions of $n$
into parts $\not\equiv0,\pm r\pmod{2k}$. Then, for   $k> r\geq 1$   and $n\geq0$,
\begin{equation*}
{A}_0(-;1,k,r;n)={B}_0(-;1,k,r;n).
\end{equation*}
\end{thm}

Furthermore, Bressoud \cite{Bressoud-1980} obtained the following generating function version  of Theorem \ref{R-R-B}: For $k> r\geq 1$,
\begin{equation}\label{R-R-BB1}
\sum_{N_1\geq  \cdots \geq N_{k-1}\geq0}\frac{q^{N_1^2+\cdots+N^2_{k-1}+N_r+\cdots+N_{k-1}}}
{(q;q)_{N_1-N_2}\cdots(q;q)_{N_{k-2}-N_{k-1}}(q^{2};q^{2})_{N_{k-1}}}=\frac{(q^r,q^{2k-r},q^{2k};q^{2k})_\infty}{(q;q)_\infty}.
\end{equation}

Motivated by the Rogers-Ramanujan-Gordon identities,   Andrews \cite{Andrews-1967} found an infinite family of the combinatorial generalizations of the G\"{o}llnitz-Gordon identities, which has been referred to as the Andrews-G\"ollnitz-Gordon theorem.

\begin{thm}[Andrews-G\"ollnitz-Gordon]\label{Gordon-Gollnitz-odd}
For $k\geq r\geq 1$, let $B_1(1;2,k,r;n)$ denote the number of partitions  $\pi=(\pi_1,\pi_2,\ldots,\pi_\ell)$ of $n$  such that no odd part is repeated,
 where $\pi_i\geq\pi_{i+k-1}+2$ with strict inequality if $\pi_i$  is even for $1\leq i\leq \ell-k+1$, and at most $r-1$ of the $\pi_i$ are less than or equal to $2$. For $k\geq r\geq 1$, let $A_1(1;2,k,r;n)$ denote the number of partitions of $n$ into parts $\not\equiv2\pmod4$ and $\not\equiv0,\pm(2r-1)\pmod{4k}$. Then, for   $k\geq r\geq1$ and $n\geq0$,
\[A_1(1;2,k,r;n)=B_1(1;2,k,r;n).\]
\end{thm}

In 1980, Bressoud \cite{Bressoud-1980} extended the Andrews-G\"ollnitz-Gordon theorem to even moduli, which has been called the Bressoud-G\"ollnitz-Gordon theorem.
\begin{thm}[Bressoud-G\"ollnitz-Gordon]\label{Gordon-Gollnitz-even} For $k> r\geq 1$, let $B_0(1;2,k,r;n)$ denote the number of partitions  $\pi=(\pi_1,\pi_2,\ldots,\pi_\ell)$ of $n$ such that no odd part is repeated, where $\pi_i\geq\pi_{i+k-1}+2 $  with strict inequality if $\pi_i$ is even for $1\leq i\leq \ell-k+1$,  at most $r-1$ of the $\pi_i$ are less than or equal to $2$, and for $1\leq i\leq \ell-k+2$, if $\pi_i\leq\pi_{i+k-2}+2$ with strict inequality if $\pi_i$ is odd, then
   \[\pi_i+\cdots+\pi_{i+k-2}\equiv r-1+V_\pi(\pi_i)\pmod{2},\]
where $V_\pi(t)$   denotes the number of odd parts not exceeding $t$ in $\pi$. For $k> r\geq 1$, let $A_0(1;2,k,r;n)$ denote the number of partitions of $n$ into parts not congruent to $ 2k-1\pmod{4k-2}$ may be repeated, no part is congruent to $2\pmod 4$, no part is multiple of $8k-4$, and no part is congruent to $\pm(2r-1)\pmod {4k-2}$. Then, for  $k> r\geq1$ and $n\geq0$,
\[A_0(1;2,k,r;n)=B_0(1;2,k,r;n).\]
\end{thm}

Bressoud \cite{Bressoud-1980} derived the following generating function versions of Theorem \ref{Gordon-Gollnitz-odd} and Theorem \ref{Gordon-Gollnitz-even}:  For $j=0$ or $1$ and $(2k+j)/2> r\geq 1$,
\begin{equation}\label{Gollnitz-A1}
\begin{split}
&\sum_{N_1\geq\cdots\geq N_{k-1}\geq0}\frac{(-q^{1-2N_1};q^2)_{N_1}q^{2(N^2_1+\cdots+N^2_{k-1}+N_r+\cdots+N_{k-1})}}{(q^2;q^2)_{N_1-N_2} \cdots(q^2;q^2)_{N_{k-2}-N_{k-1}}(q^{4-2j};q^{4-2j})_{N_{k-1}}}\\[5pt]
&=\frac{(q^2;q^4)_\infty(q^{2r-1},q^{4k-2r-1+2j},q^{4k-2+2j};q^{4k-2+2j})_\infty}{(q;q)_\infty}.
\end{split}
\end{equation}

For $j=0$ or $1$, it is evident that the generating function of $A_j(1;2,k,r;n)$ defined in Theorem \ref{Gordon-Gollnitz-odd} and Theorem \ref{Gordon-Gollnitz-even}  equals the right-hand side of \eqref{Gollnitz-A1}. Hence, the sum on the left-hand side of \eqref{Gollnitz-A1} can be considered as the generating function of
$B_j(1;2,k,r;n)$ defined in Theorem \ref{Gordon-Gollnitz-odd} and Theorem \ref{Gordon-Gollnitz-even}. More precisely,  when $j=1$, the identity \eqref{Gollnitz-A1} can be viewed as the generating function version of Theorem \ref{Gordon-Gollnitz-odd}, and when $j=0$, the identity \eqref{Gollnitz-A1} can be seen as the generating function version of Theorem \ref{Gordon-Gollnitz-even}.

 Bressoud  obtained   a far-reaching partition
 theorem utilizing an extension of  Watson's $q$-analogue of Whipple's theorem (see \cite[Theorem 1]{Bressoud-1980}).  Throughout this paper, we assume that $\alpha_1,\alpha_2,\ldots, \alpha_\lambda$ and $\eta$ are integers such that
\begin{equation}\label{cond-alpha}
0<\alpha_1<\cdots<\alpha_\lambda<\eta, \quad \text{and} \quad \alpha_i=\eta-\alpha_{\lambda+1-i}\quad \text{for} \quad 1\leq i\leq \lambda.
\end{equation}
When $\lambda$ is odd, observing that
 $\eta=\alpha_{(\lambda+1)/2}+\alpha_{\lambda+1-(\lambda+1)/2}=2\alpha_{(\lambda+1)/2}$, we see that $\eta$ must be even in such case.

\begin{thm}[Bressoud]\label{BRESSOUD-EQN}
 For $j=0$ or $1$ and $(2k+j)/2> r\geq\lambda\geq0$,
\begin{equation}\label{Bressoud-conj-e}
 \begin{split}
 &\sum_{N_1\geq\cdots \geq N_{k-1}\geq0}\frac{q^{\eta(N_1^2+\cdots+N_{k-1}^2+N_r
 +\cdots+N_{k-1})}}{(q^\eta;q^\eta)_{N_1-N_2}
 \cdots(q^\eta;q^\eta)_{N_{k-2}-N_{k-1}}(q^{(2-j)\eta};q^{(2-j)\eta})_{N_{k-1}}}\\[5pt]
 &\hskip 1.5cm\times\prod_{s=1}^{\lambda}(-q^{\eta-\alpha_s-\eta N_s};q^\eta)_{N_s}\prod_{s=2}^{\lambda}(-q^{\eta -\alpha_s+\eta N_{s-1}};q^\eta)_\infty\\[10pt]
  &=\frac{(-q^{\alpha_1},\ldots,-q^{\alpha_\lambda};q^{\eta})_\infty(q^{\eta(r-\frac{\lambda}{2})},q^{\eta
 (2k-r-\frac{\lambda}{2}+j)},q^{\eta(2k-\lambda+j)}
 ;q^{\eta(2k-\lambda+j)})_\infty}{(q^\eta;q^\eta)_\infty}.
 \end{split}
 \end{equation}
 \end{thm}
This theorem reduces to many infinite families of identities. For example, setting $\lambda=0,\ \eta=1$ and $j=1$ or $0$, we recover \eqref{R-R-A1} and \eqref{R-R-BB1} respectively. Setting $\lambda=1,\ \eta=2$ and $\alpha_1=1$, we come to \eqref{Gollnitz-A1}. To give a  combinatorial interpretation of \eqref{Bressoud-conj-e}, Bressoud introduced two partition functions.

\begin{defi}[Bressoud]\label{Bressoud-b-defi-0000000001}  For $j=0$ or $1$ and $k\geq r\geq\lambda\geq0$, define the partition function
$B_j(\alpha_1,\ldots,\alpha_\lambda;\eta,k,r;n)$  to be the number of partitions $\pi=(\pi_1,\pi_2,\ldots,\pi_\ell)$ of $n$   satisfying the following conditions{\rm{:}}
\begin{itemize}
\item[{\rm (1)}] For $1\leq i\leq \ell$, $\pi_i\equiv0,\alpha_1,\ldots,\alpha_\lambda\pmod{\eta}${\rm{;}}

\item[{\rm (2)}] Only multiples of $\eta$ may be repeated{\rm{;}}

\item[{\rm (3)}]  For $1\leq i\leq\ell-k+1$, $ \pi_i\geq\pi_{i+k-1}+\eta$  with strict inequality if $\eta\mid\pi_i${\rm{;}}

\item[{\rm (4)}] At most $r-1$ of the $\pi_i$ are less than or equal to $\eta${\rm{;}}

\item[{\rm (5)}] For $1\leq i\leq\ell-k+2$, if $\pi_i\leq\pi_{i+k-2}+\eta$ with strict inequality if $\eta\nmid\pi_i$,  then \[[\pi_i/\eta]+\cdots+[\pi_{i+k-2}/\eta]\equiv r-1+V_\pi(\pi_i)\pmod{2-j},\]
  where $V_\pi(t)$ denotes the number of parts not exceeding $t$ which are not divisible by $\eta$ in $\pi$ and $[\ ]$ denotes the greatest integer function.
  \end{itemize}
 \end{defi}

\begin{defi}[Bressoud]\label{defi-A_j} For $j=0$ or $1$ and $(2k+j)/2> r\geq\lambda\geq0$, define the partition function $A_j(\alpha_1,\ldots,\alpha_\lambda;\eta,k,r;n)$ to be the number of partitions of $n$ into parts congruent to $0,\alpha_1,\ldots,\alpha_\lambda\pmod\eta$ such that
\begin{itemize}
\item[{\rm (1)}] If $\lambda$ is even, then only multiples of $\eta$ may be repeated and no part is congruent to $0,\pm\eta(r-\lambda/2)  \pmod{\eta(2k-\lambda+j)}${\rm{;}}
\item[{\rm (2)}] If $\lambda$ is odd and $j=1$, then only multiples of ${\eta}/{2}$ may be repeated, no part is congruent to $\eta\pmod{2\eta}$, and no part is congruent to $0,\pm{\eta}(2r-\lambda)/{2} \pmod {\eta(2k-\lambda+1)}${\rm{;}}
 \item[{\rm (3)}]If $\lambda$ is odd and $j=0$, then only multiples of ${\eta}/{2}$ which are not congruent to ${\eta}(2k-\lambda)/{2}\pmod{\eta(2k-\lambda)}$ may be repeated, no part is congruent to $\eta\pmod{2\eta}$,  no part is congruent to $0\pmod {2\eta(2k-\lambda)}$, and no part is congruent to $\pm{\eta}(2r-\lambda)/{2} \pmod {\eta(2k-\lambda)}$.
  \end{itemize}
  \end{defi}

  Bressoud \cite{Bressoud-1980}  posed the following conjecture.

 \begin{conj}[Bressoud] \label{Bressoud-conjecture-j} For $j=0$ or $1$,
 $(2k+j)/2> r\geq\lambda\geq0$  and $n\geq 0$,
 \begin{equation*}\label{Bressoud-conj-1}
A_j(\alpha_1,\ldots,\alpha_\lambda;\eta,k,r;n)=B_j(\alpha_1,
\ldots,\alpha_\lambda;\eta,k,r;n).
 \end{equation*}
\end{conj}

This conjecture specializes to many infinite families of combinatorial identities. For example, setting $\lambda=0,\ \eta=1$  and $j=1$ or $0$, we find that it reduces to Theorem \ref{R-R-B-1} and Theorem  \ref{R-R-B} respectively. For $\lambda=1,\ \eta=2$, $\alpha_1=1$ and $j=1$ or $0$, we see that it  boils down to Theorem \ref{Gordon-Gollnitz-odd} and Theorem \ref{Gordon-Gollnitz-even} respectively.

As remarked by Bressoud \cite{Bressoud-1980}, it is not difficult to see that the generating function of $A_j(\alpha_1,\ldots,\alpha_\lambda;\eta,k,r;n)$ is equal to the right-hand side of \eqref{Bressoud-conj-e}.
 \begin{thm}[Bressoud]\label{Bressoud-conj-epr-thm}
For  $j=0$ or $1$ and  $(2k+j)/2> r\geq\lambda\geq0$,
 \begin{equation}\label{Bressoud-conj-epr}
 \begin{split}
 & \sum_{n\geq0}{A}_j(\alpha_1,\ldots,\alpha_\lambda;\eta,k,r;n)q^n\\[5pt]
  &=\frac{(-q^{\alpha_1},\ldots,-q^{\alpha_\lambda};q^{\eta})_\infty(q^{\eta(r-\frac{\lambda}{2})},q^{\eta
 (2k-r-\frac{\lambda}{2}+j)},q^{\eta(2k-\lambda+j)}
 ;q^{\eta(2k-\lambda+j)})_\infty}{(q^\eta;q^\eta)_\infty}.
 \end{split}
 \end{equation}
 \end{thm}

Nevertheless,  it does not seem easy to prove that the left-hand side of \eqref{Bressoud-conj-e} is indeed the generating function of $B_j(\alpha_1,\ldots,\alpha_\lambda;\eta,k,r;n)$. In this regard, Bressoud \cite{Bressoud-1980}  posed the following conjecture.

\begin{conj}[Bressoud] \label{Bressoud-gen-b-e} For $j=0$ or $1$  and $(2k+j)/2> r\geq\lambda\geq0$,
\begin{equation*}
 \begin{split}
 &\sum_{n\geq0}B_j(\alpha_1,\ldots,\alpha_\lambda;\eta,k,r;n)q^n\\
 &=\sum_{N_1\geq\cdots\geq N_{k-1}\geq0}\frac{q^{\eta(N_1^2+\cdots+N_{k-1}^2+N_r+\cdots+N_{k-1})}}{(q^\eta;q^\eta)_{N_1-N_2}\cdots(q^\eta;q^\eta)_{N_{k-2}-N_{k-1}}(q^{(2-j)\eta};q^{(2-j)\eta})_{N_{k-1}}}\\
 &\hskip1cm\times\prod_{s=1}^{\lambda}(-q^{\eta-\alpha_s-\eta N_s};q^\eta)_{N_s}\prod_{s=2}^{\lambda}(-q^{\eta -\alpha_s+\eta N_{s-1}};q^\eta)_\infty.
 \end{split}
 \end{equation*}
\end{conj}

Andrews \cite{Andrews-1974m}  proved Conjecture \ref{Bressoud-conjecture-j}   for $\eta=\lambda+1$ and $j=1$.  Kim and Yee \cite{Kim-Yee-2014} showed that the conjecture holds for $j=1$ and $\lambda=2$.
In fact, they proved that Conjecture \ref{Bressoud-gen-b-e}
is true for $j=1$ and $\lambda=2$ with the aid of Gordon markings introduced by  Kur\c{s}ung\"oz   \cite{Kursungoz-2010a, Kursungoz-2010}.
Recently, Kim  \cite{Kim-2018} resolved  Conjecture \ref{Bressoud-conjecture-j} for the case $j=1$. To this end, she established the following theorem.

\begin{thm}[Bressoud-Kim]\label{berssoud-kim-1}
For $k\geq r\geq\lambda\geq0$,
\begin{equation}\label{Bressoud-conj-defi-e-11}
 \begin{split}
 &\sum_{n\geq0}B_1(\alpha_1,\ldots,\alpha_\lambda;\eta,k,r;n)q^n\\[5pt]
 &=\frac{(-q^{\alpha_1},\ldots,-q^{\alpha_\lambda};q^{\eta})_\infty
 (q^{\eta(r-\frac{\lambda}{2})},q^{\eta(2k-r-\frac{\lambda}{2}+1)},
 q^{\eta(2k-\lambda+1)};q^{\eta(2k-\lambda+1)})_\infty}{(q^\eta;q^\eta)_\infty}.
 \end{split}
 \end{equation}
 \end{thm}

 It is clear that Conjecture \ref{Bressoud-conjecture-j} for $j=1$
 is an immediate consequence of \eqref{Bressoud-conj-epr} and \eqref{Bressoud-conj-defi-e-11}.

The main objective of this paper is to give  overpartition analogues of
the partition function $B_j$ and the partition function $A_j$ introduced by
Bressoud and to establish overpartition analogues of some classical
partition theorems. The overpartition analogues of classical partition theorems
have caught much attention, see, for example, Chen, Sang and Shi \cite{Chen-Sang-Shi-2011,Chen-Sang-Shi-2013,Chen-Sang-Shi-2013a}, Choi, Kim and Lovejoy \cite{Choi-Kim-Lovejoy-2017}, Corteel and Lovejoy \cite{Corteel-Lovejoy-2009}, Corteel, Lovejoy and Mallet \cite{Corteel-2008}, Corteel and Mallet \cite{Corteel-2007}, Dousse \cite{Dousse-2015,Dousse-2017}, Goyal \cite{Goyal-2018}, He, Ji, Wang and Zhao \cite{He-Ji-Wang-Zhao-2019}, He, Wang and Zhao \cite{He-Wang-Zhao-2017}, Kur\c{s}ung\"{o}z \cite{Kursungoz-2016}, Lovejoy \cite{Lovejoy-2003,Lovejoy-2004,Lovejoy-2005,Lovejoy-2007,Lovejoy-2010},  Lovejoy and Mallet \cite{Lovejoy-Mallet-2008}, Raghavendra and Padmavathamma \cite{P-R-2009}, and Sang and Shi \cite{Sang-Shi-2015}.

Lovejoy \cite{Lovejoy-2003}  established overpartition analogues of the Rogers-Ramanujan-Gordon theorem for the cases $i = 1$ and $i = k$, and the general case was obtained by
Chen, Sang and Shi \cite{Chen-Sang-Shi-2013}.  In Theorem \ref{R-R-G-o} and for the rest of this paper,
we adopt the following convention: For positive integers $t$ and $b$, we define $t\pm b$ (resp. $\overline{t}\pm b$) as a non-overlined part (resp. an overlined part) of size  $t\pm b$. The parts in an overpartition are ordered as in \eqref{order}.

 \begin{thm}[Chen-Sang-Shi]\label{R-R-G-o}
For $k\geq r\geq1$, let $\overline{B}_1(-;1,k,r;n)$ denote the number of overpartitions $\pi=(\pi_1,\pi_2,\ldots,\pi_\ell)$ of $n$, where  $\pi_i\geq\pi_{i+k-1}+1$ with strict inequality if $\pi_i$ is non-overlined for $1\leq i\leq \ell-k+1$, and at most $r-1$ of the $\pi_i$ are equal to $1$. For $k>r\geq1$, let $\overline{A}_1(-;1,k,r;n)$ denote the number of overpartitions of $n$ such that non-overlined parts $\not\equiv0,\pm r\pmod{2k}$, and for $k=r$, let $\overline{A}_1(-;1,k,k;n)$ denote the number of overpartitions of $n$ into parts not divisible by $k$. Then, for $k\geq r\geq1$ and    $n\geq0$,
\begin{equation*}
\overline{A}_1(-;1,k,r;n)=\overline{B}_1(-;1,k,r;n).
\end{equation*}
\end{thm}

Chen, Sang and Shi \cite{Chen-Sang-Shi-2013}  gave the following generating function version of Theorem \ref{R-R-G-o}: For  $k\geq r \geq 1$,
\[
  \begin{split}
  &\sum_{N_1\geq\cdots\geq N_{k-1}\geq0}
 \frac{q^{N_{1}^{2}+\cdots+N_{k-1}^{2}+ N_r+\cdots+N_{k-1}}(1+q^{- N_{r}})(-q^{1- N_{1}};q)_{N_{1}-1}}
 {(q;q)_{N_{1}-N_{2}}\cdots(q;q)_{N_{k-2}-N_{k-1}}(q
 ;q)_{N_{k-1}}}\\[5pt]
  &=\frac{(-q;q)_{\infty}(q^{r},
q^{2k-r}
  ,q^{2k};q^{2k})_{\infty}}
{(q;q)_{\infty}}.
\end{split}
\]

Corteel, Lovejoy and Mallet \cite{Corteel-2008} established an overpartition
analogue  of the Bressoud-Rogers-Ramanujan theorem for the case $i = 1$,
and the general case was obtained by  Chen, Sang and Shi \cite{Chen-Sang-Shi-2013a}.

 \begin{thm}[Chen-Sang-Shi]\label{R-R-G-o-0}
For $k\geq r\geq1$, let $\overline{B}_0(-;1,k,r;n)$ denote the number of overpartitions $\pi=(\pi_1,\pi_2,\ldots,\pi_\ell)$ of $n$, where  $\pi_i\geq\pi_{i+k-1}+1$ with strict inequality if $\pi_i$ is non-overlined for $1\leq i\leq \ell-k+1$,  at most $r-1$ of the $\pi_i$ are equal to $1$,  and for $1\leq i\leq \ell-k+2$, if $\pi_i\leq\pi_{i+k-2}+1$ with strict inequality if   $\pi_i$ is overlined, then
 \[\pi_i+\cdots+\pi_{i+k-2}\equiv r-1+\overline{V}_\pi(\pi_i)\pmod{2}.\]
For $k\geq r\geq1$, let $\overline{A}_0(-;1,k,r;n)$ denote the number of overpartitions of $n$ such that non-overlined parts $\not\equiv0,\pm r\pmod{2k-1}$. Then, for  $k\geq r\geq1$ and  $n\geq0$,
\begin{equation*}
\overline{A}_0(-;1,k,r;n)=\overline{B}_0(-;1,k,r;n).
\end{equation*}

\end{thm}

In Theorem \ref{R-R-G-o-0} and for the rest of this article,
  $\overline{V}_\pi(t)$ (resp.   $\overline{V}_\pi(\overline{t})$), as used by Corteel, Lovejoy and Mallet \cite{Corteel-2008},  stands for the number of overlined parts not exceeding $t$ (resp.  $\overline{t}$) in $\pi$. 
For example,
for an overpartition
$\pi=(\overline{7},7,6,\overline{5},5,\overline{2})$,
we have $\overline{V}_\pi({5})=1$ and  $\overline{V}_\pi(\overline{5})=2$.

The following generating function version of Theorem \ref{R-R-G-o-0} was given by Sang and Shi \cite{Sang-Shi-2015}: For  $k> r \geq 1$,
\begin{equation*}\label{R-R-G-o-0-eqn}
  \begin{split}
  &\sum_{N_1\geq\cdots\geq N_{k-1}\geq0}
 \frac{q^{N_{1}^{2}+\cdots+N_{k-1}^{2}+ N_r+\cdots+N_{k-1}}(1+q^{- N_{r}})(-q^{1- N_{1}};q)_{N_{1}-1}}
 {(q;q)_{N_{1}-N_{2}}\cdots(q;q)_{N_{k-2}-N_{k-1}}(q^{2}
 ;q^{2})_{N_{k-1}}}\\[5pt]
  &=\frac{(-q;q)_{\infty}(q^{r},
q^{2k-r-1}
  ,q^{2k-1};q^{2k-1})_{\infty}}
{(q;q)_{\infty}}.
\end{split}
\end{equation*}

 In this paper, we introduce two new partition functions $\overline{B}_j(\alpha_1,\ldots,\alpha_\lambda;\eta, k,r;n)$ and $\overline{A}_j(\alpha_1,\ldots,\alpha_\lambda;\eta, k,r;n)$ and build connections
 between $\overline{B}_j(\alpha_1,\ldots,\alpha_\lambda;\eta, k,r;n)$ and ${B}_j(\alpha_1\break$$,\ldots,\alpha_\lambda;\eta, k,r;n)$.

 \begin{defi}\label{defi-O-B} For $j=0$ or $1$ and  $k\geq r\geq \lambda\geq0$, define $\overline{B}_j(\alpha_1,\ldots,\alpha_\lambda;\eta,k,r;n)$ to be the number of overpartitions $\pi=(\pi_1,\pi_2,\ldots,\pi_\ell)$ of $n$   subject to the following conditions{\rm{:}}
 \begin{itemize}
  \item[{\rm (1)}] For $1\leq i\leq \ell$, $\pi_i\equiv0,\alpha_1,\ldots,\alpha_\lambda\pmod{\eta}${\rm{;}}

 \item[{\rm (2)}] Only  multiples of $\eta$ may be non-overlined{\rm{;}}

 \item[{\rm (3)}]  For $1\leq i\leq \ell-k+1$, $\pi_i\geq\pi_{i+k-1}+\eta$ with strict inequality if  $\pi_i$  is non-overlined{\rm{;}}

 \item[{\rm (4)}] At most $r-1$ of the $\pi_i$ are less than or equal to $\eta${\rm{;}}

 \item[{\rm (5)}] For $1\leq i\leq \ell-k+2$, if $\pi_i\leq\pi_{i+k-2}+\eta$ with strict inequality if   $\pi_i$ is overlined, then
 \[\left[\pi_i/\eta\right]+\cdots+\left[\pi_{i+k-2}/\eta\right]\equiv r-1+\overline{V}_{\pi}(\pi_i)\pmod{2-j}.\]
\end{itemize}
\end{defi}

\begin{defi}\label{defi-O-A}
 For $j=0$ or $1$ and  $(2k-j)/2\geq r\geq \lambda\geq0$, define   the partition function $\overline{A}_j(\alpha_1,\ldots,\alpha_\lambda;\eta,k,r;n)$ to be the number of overpartitions of $n$ satisfying
$\pi_i\equiv0,\alpha_1,\ldots,\alpha_\lambda\pmod{\eta}$ such that
 \begin{itemize}
 \item[{\rm (1)}] If $\lambda$ is even, then only multiplies of $\eta$ may be non-overlined and  there is no non-overlined part congruent to $0,\pm\eta(r-\lambda/2) \pmod {\eta(2k-\lambda+j-1)}${\rm;}

 \item[{\rm (2)}]  If $\lambda$ is odd and $j=1$, then only multiples of ${\eta}/{2}$ may be non-overlined, no non-overlined part is congruent to ${\eta}(2k-\lambda)/{2}\pmod{\eta(2k-\lambda)}$, no non-overlined part is congruent to $\eta \pmod{2\eta}$, no non-overlined part is congruent to $0 \pmod{2\eta(2k-\lambda)}$, no non-overlined part is congruent to $\pm{\eta}(2r-\lambda)/{2} \pmod {\eta(2k-\lambda)}$, and no overlined part is congruent to ${\eta}/{2}\pmod \eta$ and not congruent to ${\eta}(2k-\lambda)/{2}\pmod{\eta(2k-\lambda)}${\rm;}
 \item[{\rm (3)}]  If $\lambda$ is odd and $j=0$, then only multiples of  ${\eta}/{2}$  may be non-overlined, no non-overlined part is congruent to $\eta \pmod{2\eta}$,  no non-overlined part is congruent to $0,\pm{\eta}(2r-\lambda)/{2} \pmod {\eta(2k-\lambda-1)}$, and no overlined part is congruent to ${\eta}/{2}\pmod \eta$.
  \end{itemize}
    \end{defi}

Observe that for an overpartition $\pi$ counted by $\overline{B}_j(\alpha_1,\ldots,\alpha_\lambda;\eta,k,r;n)$ (resp. $\overline{A}_{j}(\alpha_1,\break$ $\ldots,\alpha_\lambda;\eta,k,r;n)$) without overlined parts divisible by $\eta$, if we change the overlined parts in $\pi$ to non-overlined parts, then we get an ordinary partition counted by ${B}_j(\alpha_1,\ldots,\alpha_\lambda;\eta,k,\break $$ r;n)$ (resp. ${A}_{1-j}(\alpha_1,\ldots,\alpha_\lambda;\eta,k-1+j,r;n)$). Hence we  say that $\overline{B}_j(\alpha_1,\ldots,\alpha_\lambda;\eta,k,r;n)$ (resp. $\overline{A}_{j}(\alpha_1,\ldots,\alpha_\lambda;\eta,k,r;n)$) can be considered as an overpartition analogue of  ${B}_j(\alpha_1,\break$$ \ldots,\alpha_\lambda;\eta,
k,r;n)$ (resp.  ${A}_j(\alpha_1,\ldots,\alpha_\lambda;\eta,k,r;n)$). In this case,  $\overline{V}_\pi(t)$ reduces to the notation   ${V}_\pi(t)$ introduced by Bressoud \cite{Bressoud-1980}.

By means of Gordon markings, we build bijections to obtain   the following   relationships between $\overline{B}_j(\alpha_1,\ldots,\alpha_\lambda;\eta,k,r;n)$ and ${B}_j(\alpha_1,\ldots,\alpha_\lambda;\eta,k,r;n)$.
\begin{thm}\label{rel-over1}For   $k\geq r\geq \lambda\geq0$ and $k>\lambda$,
 \begin{equation*}\label{b-0-over}
 \begin{split}
 &\sum_{n\geq0}\overline{B}_1(\alpha_1,\ldots,\alpha_\lambda;
 \eta,k,r;n)q^n=(-q^\eta;q^\eta)_\infty
 \sum_{n\geq0}B_0(\alpha_1,\ldots,\alpha_\lambda;
 \eta,k,r;n)q^n.
 \end{split}
\end{equation*}
\end{thm}

  \begin{thm} \label{rel-over2}For   $k> r\geq \lambda\geq0$ and $k-1>\lambda$,
 \begin{equation*}\label{b-0-over-000}
 \begin{split}
 &\sum_{n\geq0}\overline{B}_0(\alpha_1,\ldots,\alpha_\lambda;
 \eta,k,r;n)q^n=(-q^\eta;q^\eta)_\infty
 \sum_{n\geq0}B_1(\alpha_1,\ldots,\alpha_\lambda;
 \eta,k-1,r;n)q^n.
 \end{split}
\end{equation*}
For  $k-1> \lambda$,
 \begin{equation*}\label{b-0-overkk}
 \begin{split}
 &\sum_{n\geq0}\overline{B}_0(\alpha_1,\ldots,\alpha_\lambda;
 \eta,k,k;n)q^n=(-q^\eta;q^\eta)_\infty
 \sum_{n\geq0}B_1(\alpha_1,\ldots,\alpha_\lambda;
 \eta,k-1,k-1;n)q^n.
 \end{split}
\end{equation*}
\end{thm}

We also derive the generating function of $\overline{A}_j(\alpha_1,\ldots,\alpha_\lambda;\eta,k,r;n)$.

\begin{thm}\label{gf-overlineA_j} For  $j=0$ or $1$ and  $(2k-j)/2\geq r\geq\lambda\geq0$,
 \begin{equation}\label{gf-overlineA_j-eq}
 \begin{split}
 & \sum_{n\geq0}\overline{A}_{j}(\alpha_1,\ldots,\alpha_\lambda;\eta,k,r;n)q^n\\[5pt]
  &=\frac{(-q^{\alpha_1},\ldots,-q^{\alpha_\lambda},-q^{\eta};q^{\eta})_\infty(q^{\eta(r-\frac{\lambda}{2})},q^{\eta
 (2k-r-\frac{\lambda}{2}-1+j)},q^{\eta(2k-\lambda-1+j)}
 ;q^{\eta(2k-\lambda-1+j)})_\infty}{(q^\eta;q^\eta)_\infty}.
 \end{split}
 \end{equation}
 \end{thm}

By Theorem \ref{R-R-B} and  Theorem \ref{rel-over1} with $\lambda=0$ and $\eta=1$, we find that for $k>r\geq 1$,
\[
\sum_{n\geq 0}\overline{B}_1(-;1,k,r;n)q^n=\frac{(-q;q)_{\infty}(q^{r},
q^{2k-r}
  ,q^{2k};q^{2k})_{\infty}}
{(q;q)_{\infty}}.
\]
Combining with  Theorem \ref{gf-overlineA_j} with $\lambda=0$, $\eta=1$ and $j=1$,  we can  recover Theorem \ref{R-R-G-o} for $k>r\geq 1$.  By Theorem \ref{Gordon-Gollnitz-even} and  Theorem \ref{rel-over1} with $\lambda=1$ and $\eta=2$, we find that for $k>r\geq 1$,
\[
\sum_{n\geq 0}\overline{B}_1(1;2,k,r;n)q^n=\frac{(-q^2;q^2)_{\infty}(-q;q^2)_{\infty}(q^{2r-1},
q^{4k-2r-1}
  ,q^{4k-2};q^{4k-2})_{\infty}}
{(q^2;q^2)_{\infty}}.
\]
Applying Theorem  \ref{gf-overlineA_j}  with $\lambda=1$, $\eta=2$ and $j=1$, we   obtain a new overpartition analogue of the Andrews-G\"ollnitz-Gordon theorem.

\begin{thm}\label{Over-bre-121} For $k> r\geq 1$,
let $\overline{B}_1(1;2,k,r;n)$ denote the number of overpartitions $\pi=(\pi_1,\pi_2,\ldots,\pi_\ell)$ of $n$  such that {\rm (1)} only even parts may be non-overlined{\rm{;}} {\rm (2)} $\pi_i\geq\pi_{i+k-1}+2$ with strict inequality if $\pi_i$ is  non-overlined for $1\leq i\leq \ell-k+1${\rm{;}} {\rm (3)} at most $r-1$ of the $\pi_i$ are less than or equal to $2$.

For $k> r\geq 1$,
let   $\overline{A}_1(1;2,k,r;n)$ denote the number of  overpartitions of $n$ such that {\rm (1)} no non-overlined part is congruent to $2k-1\pmod{4k-2}${\rm{;}} {\rm (2)} no non-overlined part is congruent to $2 \pmod{4}${\rm{;}} {\rm (3)} no non-overlined part is congruent to $0 \pmod{8k-4}${\rm{;}} {\rm (4)} no non-overlined part is congruent to $\pm(2r-1) \pmod {4k-2}${\rm{;}} {\rm (5)} no overlined part is congruent to $1\pmod 2$ and not congruent to $2k-1\pmod{4k-2)}$. Then, for $k> r\geq1$ and  $n\geq0$,
  \begin{equation*}\label{main000000}
 \overline{A}_1(1;2,k,r;n)=\overline{B}_1(1;2,k,r;n).
 \end{equation*}
\end{thm}
The generating function version of Theorem \ref{Over-bre-121} will be given in our subsequent paper \cite{He-Ji-Zhao}. It should be  mentioned  that  Lovejoy \cite{Lovejoy-2004}   obtained an  overpartition analogue of the Andrews-G\"ollnitz-Gordon theorem for $r=k$ and He, Ji, Wang and Zhao \cite{He-Ji-Wang-Zhao-2019} found an overpartition analogue for  the general case.

In view of Theorems \ref{berssoud-kim-1}, \ref{rel-over2} and Theorem \ref{gf-overlineA_j} for $j=0$,  we   obtain the following  overpartition analogue of Bressoud's Conjecture \ref{Bressoud-conjecture-j} for $j=0$.

\begin{thm}\label{Over-bre-1}
For  $k\geq r\geq \lambda\geq0$,  $k-1>\lambda$ and $n\geq 0$, we have
  \begin{equation*}\label{main}
 \overline{A}_0(\alpha_1,\ldots,\alpha_\lambda;
 \eta,k,r;n)=\overline{B}_0(\alpha_1,\ldots,\alpha_\lambda;\eta,k,r;n).
 \end{equation*}
\end{thm}

The generating function version   of  Theorem \ref{Over-bre-1} can be derived with the aid of Bailey pairs.

\begin{thm}\label{G-B-O-1} For  $k\geq r>\lambda \geq 0$,
\begin{align}\label{G-B-O-1-eq}
  &\sum_{N_1\geq\cdots\geq N_{k-1}\geq0}
 \frac{q^{\eta(N_{1}^{2}+\cdots+N_{k-1}^{2}+ N_r+\cdots+N_{k-1})}(1+q^{-\eta N_{r}})(-q^{\eta-\eta N_{\lambda+1}};q^{\eta})_{N_{\lambda+1}-1}
 }
 {(q^{\eta};q^{\eta})_{N_{1}-N_{2}}\cdots(q^\eta;q^\eta)_{N_{k-2}-N_{k-1}}(q^{2\eta}
 ;q^{2\eta})_{N_{k-1}}}\nonumber\\
 &\hskip 1.5cm\times
(-q^{\eta+\eta N_{\lambda}};q^{\eta})_{\infty} \prod_{s=1}^\lambda(-q^{\eta-\alpha_{s}-\eta N_{s}};q^{\eta})_{N_{s}}
\prod_{s=2}^\lambda(-q^{\eta-\alpha_{s}+\eta N_{s-1}};q^{\eta})_{\infty}\nonumber\\[10pt]
&=\frac{(-q^{\alpha_{1}},\ldots
 -q^{\alpha_{\lambda}},-q^{\eta};q^{\eta})_{\infty}(q^{(r-\frac{\lambda}{2})\eta},
q^{(2k-r-\frac{\lambda}{2}-1)\eta}
  ,q^{(2k-\lambda-1)\eta};q^{(2k-\lambda-1)\eta})_{\infty}}
{(q^{\eta};q^{\eta})_{\infty}}.
\end{align}
\end{thm}

Combining Theorem \ref{gf-overlineA_j} for $j=0$, Theorem \ref{Over-bre-1} and Theorem \ref{G-B-O-1}, we obtain the following generating function of $\overline{B}_0(\alpha_1,\ldots,\alpha_\lambda;\eta,k,r;n)$, which can be regarded as the overpartition analogue of Bressoud's Conjecture \ref{Bressoud-gen-b-e}  for $j=0$.

\begin{thm}\label{sum-side-conj} For  $k\geq r>\lambda \geq 0$ and $k-1>\lambda$,
\begin{equation*}
 \begin{split}
 &\sum_{n\geq0}\overline{B}_0(\alpha_1,\ldots,
 \alpha_\lambda;
 \eta,k,r;n)q^n\\[5pt]
  &=\sum_{N_1\geq\cdots\geq N_{k-1}\geq0}
 \frac{q^{\eta(N_{1}^{2}+\cdots+N_{k-1}^{2}+ N_r+\cdots+N_{k-1})}(1+q^{-\eta N_{r}})(-q^{\eta-\eta N_{\lambda+1}};q^{\eta})_{N_{\lambda+1}-1}
 }
 {(q^{\eta};q^{\eta})_{N_{1}-N_{2}}\cdots
 (q^{\eta};q^{\eta})_{N_{k-2}-N_{k-1}}(q^{2\eta}
 ;q^{2\eta})_{N_{k-1}}}\\[5pt]
 &\hskip1cm\times
 (-q^{\eta+\eta N_{\lambda}};q^{\eta})_{\infty}\prod_{s=1}^\lambda(-q^{\eta-\alpha_{s}-\eta N_{s}};q^{\eta})_{N_{s}}
\prod_{s=2}^\lambda(-q^{\eta-\alpha_{s}+\eta N_{s-1}};q^{\eta})_{\infty}.
 \end{split}
 \end{equation*}
 \end{thm}

Theorem \ref{Over-bre-1} and Theorem \ref{G-B-O-1}  specialize to  overpartition analogues of a number of classical partition theorems. Setting $\lambda=0,\ \eta=1,\ k=3$ and $r=2$, we obtain an overpartition analogue of   Euler's partition theorem \cite{Euler-1748}. Recall that Euler's partition theorem states that for $n\geq 1$,  the number of partitions of $n$ into  odd parts  equals the number of partitions of $n$ into distinct parts.
\begin{thm}\label{euler-o}
Let $\overline{B}_0(-;1,3,2;n)$ denote the number of overpartitions $\pi=(\pi_1,\pi_2,\ldots,\break$$\pi_\ell)$ of  $n$, where $\pi_i\geq\pi_{i+2}+1$ with strict inequality if $\pi_i$  is non-overlined for $1\leq i\leq \ell-2$, and for $1\leq i\leq \ell-1$,
if $\pi_i\leq\pi_{i+1}+1$ with strict inequality if  $\pi_i$ is overlined, then $\pi_i+\pi_{i+1}\equiv 1+\overline{V}_\pi(\pi_i)\pmod{2}$.
Let $\overline{A}_0(-;1,3,2;n)$ denote the number of overpartitions of $n$  such that no non-overlined part is congruent to $\equiv 0,\pm2\pmod 5$. Then, for $n\geq0$,
\[\overline{A}_0(-;1,3,2;n)=\overline{B}_0(-;1,3,2;n).\]
\end{thm}

The generating function version takes the form:
\[\sum_{N_1\geq N_2\geq0}\frac{q^{N_1^2+N_2^2+N_2}(1+q^{-N_2})(-q^{1-N_1};q)_{N_1-1}}{(q;q)_{N_1-N_2}(q^2;q^2)_{N_2}}=\frac{(-q;q)_\infty(q^2,q^3,q^5;q^5)_\infty}{(q;q)_\infty}.\]

For an overpartition $\pi=(\pi_1, \pi_2,\ldots, \pi_\ell)$  counted by $\overline{B}_0(-;1,3,2;n)$,    if there are no overlined parts in $\pi$, then   $\overline{V}_\pi(\pi_i)=0$ for $1\leq i\leq \ell$. This implies that $\pi_i+\pi_{i+1}$ is odd if $\pi_i\leq \pi_{i+1}+1$. Hence we deduce that   $\pi_i>\pi_{i+1}$ for $1\leq i\leq \ell-1$. Therefore, $\pi$ is a  partition into distinct parts. For this reason,  Theorem \ref{euler-o} can be perceived  as an overpartition analogue of Euler's partition theorem.

Putting $\lambda=0$ and $\eta=1$ in Theorem \ref{Over-bre-1}, we are led to the overpartition analogue of the Bressoud-Rogers-Ramanujan theorem due to Chen, Sang and Shi \cite{Chen-Sang-Shi-2013a}. In a similar way,  Theorem \ref{G-B-O-1}  yields  the generating function version found by Sang and Shi \cite{Sang-Shi-2015}. Setting $\lambda=1$ and $\eta=2$ in Theorem \ref{Over-bre-1}, we find an overpartition analogue of the Bressoud-G\"ollnitz-Gordon theorem.

\begin{thm}\label{Over-bre-1210}
For $k>2$ and $k\geq r\geq  1$, let $\overline{B}_0(1;2,k,r;n)$ denote the number of overpartitions $\pi=(\pi_1,\pi_2,\ldots,\pi_\ell)$ of $n$  such that  only even parts may be non-overlined, $\pi_i\geq\pi_{i+k-1}+2$ with strict inequality if $\pi_i$  is  non-overlined for $1\leq i\leq \ell-k+1$, at most $r-1$ of the $\pi_i$ are less than or equal to $2$, and for $1\leq i\leq \ell-k+2$, if $\pi_i\leq\pi_{i+k-2}+2$ with strict inequality if   $\pi_i$ is overlined, then
 \[\left[\pi_i/2\right]+\cdots+\left[\pi_{i+k-2}/2\right]\equiv r-1+\overline{V}_\pi(\pi_i)\pmod{2}.\]
For $k>2$ and $k\geq r\geq  1$, let   $\overline{A}_0(1;2,k,r;n)$ denote the number of overpartitions of $n$ such that  no non-overlined part is congruent to $2 \pmod{4}$, no non-overlined part is congruent to $0,\pm(2r-1) \pmod {4k-4}$, and no overlined part is congruent to $1\pmod2$. Then, for $k>2$, $k\geq r\geq 1$ and $n\geq 0$,
  \begin{equation*}\label{main000000-0}
 \overline{A}_0(1;2,k,r;n)=\overline{B}_0(1;2,k,r;n).
 \end{equation*}
\end{thm}

 He, Wang and Zhao\cite{He-Wang-Zhao-2017} established an overpartition analogue of the Bressoud-G\"ollnitz-Gordon theorem. Putting $\lambda=1$ and $\eta=2$ in Theorem \ref{G-B-O-1}, we get the generating function version of Theorem \ref{Over-bre-1210}: For  $k\geq r>1$,
\begin{equation*}
  \begin{split}
    &\sum_{N_1\geq\cdots\geq N_{k-1}\geq0}
 \frac{q^{2(N_{1}^{2}+\cdots+N_{k-1}^{2}+ N_r+\cdots+N_{k-1})}(1+q^{-2 N_{r}})}
 {(q^{2};q^{2})_{N_{1}-N_{2}}\cdots(q^2;q^2)_{N_{k-2}-N_{k-1}}(q^{4}
 ;q^{4})_{N_{k-1}}}\\[5pt]
 &\hskip1cm \times(-q^{1-2 N_{1}};q^{2})_{N_{1}}(-q^{2-2 N_{2}};q^{2})_{N_{2}-1}
 (-q^{2+2 N_{1}};q^{2})_{\infty}\\[5pt]
 &=\frac{(-q;q^{2})_{\infty}(-q^{2};q^{2})_{\infty}(q^{2r-1},
q^{4k-2r-3}
  ,q^{4k-4};q^{4k-4})_{\infty}}
{(q^{2};q^{2})_{\infty}}.
  \end{split}
\end{equation*}

This paper is organized as follows.  In Section 2, we present a
 proof of Theorem \ref{gf-overlineA_j}. In Section 3, we  introduce the   notions of Gordon marking, reverse Gordon marking, and $(k-1)$-bands  of an overpartition counted by $\overline{B}_1(\alpha_1,\ldots,\alpha_\lambda;
 \eta,k,r;n)$. Furthermore, we give a criterion to determine whether an overpartition    counted by $\overline{B}_1(\alpha_1,\ldots,\alpha_\lambda; \eta,k,r;n)$   is   counted by $\overline{B}_0(\alpha_1,\ldots,\alpha_\lambda; \eta,k,r;n)$ as well.  In Section 4, we   define   the forward move and the backward move based on the Gordon marking and the reverse Gordon marking of an overpartition counted by $\overline{B}_1(\alpha_1,\ldots,\alpha_\lambda;
 \eta,k,r;n)$. These operations allow us to
 provide a combinatorial proof of  Theorem
 \ref{rel-over1}. Section 5 is devoted to   a
 combinatorial proof of Theorem \ref{rel-over2}. In
 Section 6, we give a proof of Theorem
 \ref{G-B-O-1} with the aid of Bailey pairs.
 In Section 7, we discuss possible directions for future work.

 \section{Proof of Theorem \ref{gf-overlineA_j}}

 As mentioned in the introduction, the function $\overline{A}_{j}(\alpha_1,\ldots,\alpha_\lambda;\eta,k,r;n)$ can be viewed as the overpartition analogue of ${A}_{j}(\alpha_1,\ldots,\alpha_\lambda;\eta,k,r;n)$ introduced by Bressoud \cite{Bressoud-1980}. Similar to the case for ${A}_{j}(\alpha_1,\ldots,\alpha_\lambda;\eta,k,r;n)$, it is not difficult to establish the generating function of $\overline{A}_{j}(\alpha_1,\ldots,\alpha_\lambda;\eta,k,r;n)$ stated as in Theorem \ref{gf-overlineA_j}. For completeness, we include  a detailed derivation.

 \noindent{\it Proof of Theorem \ref{gf-overlineA_j}.} Clearly, the right-hand side of \eqref{gf-overlineA_j-eq} can be  interpreted as the
 generating function of $\overline{A}_{j}(\alpha_1,\ldots,\alpha_\lambda;\eta,k,r;n)$  when $\lambda$ is even. It remains to show that the right-hand side of \eqref{gf-overlineA_j-eq} is also the
 generating function of $\overline{A}_{j}(\alpha_1,\ldots,\alpha_\lambda;\eta,k,r;n)$   when $\lambda$ is odd. When $\lambda$ is odd, it is clear from   \eqref{cond-alpha} that   $\eta=\alpha_{(\lambda+1)/2}+\alpha_{\lambda+1-(\lambda+1)/2}=2\alpha_{(\lambda+1)/2}$. This implies  that $\eta$ must be even in this event.

 When $j=1$,  by definition,   we have  for  $k>r\geq \lambda\geq 0$,
 \begin{equation} \label{eqn-a-1-1-1-1-1at}
 \begin{split}
 &\sum_{n\geq0}\overline{A}_1(\alpha_1,\ldots,\alpha_\lambda;\eta,k,r;n)q^n\\[5pt]
 &\qquad =(-q^{\alpha_1},\ldots,-q^{\alpha_{(\lambda-1)/2}},-q^{\alpha_{(\lambda+3)/2}},
 \ldots,-q^{\alpha_\lambda},-q^{\eta};q^{\eta})_\infty(-q^{\eta(2k-\lambda)/2};q^{\eta(2k-\lambda)})_\infty\\[8pt]
 &\qquad \qquad \times\frac{(q^{\eta(2r-\lambda)/2},
 q^{\eta(4k-2r-\lambda)/2},q^{\eta(2k-\lambda)/2};
 q^{\eta(2k-\lambda)})_\infty(q^{2\eta(2k-\lambda)};q^{2\eta(2k-\lambda)})_\infty(q^{\eta};q^{2\eta})_\infty}{(q^{\eta/2};q^{\eta/2})_\infty}.
 \end{split}
 \end{equation}
Since $\eta$ is even, we find that
  \begin{equation}\label{eqn-a-1-1-1-1-1a}
 \frac{(q^{\eta};q^{2\eta})_\infty}{(q^{\eta/2};q^{\eta/2})_\infty}=\frac{(q^{\eta/2},-q^{\eta/2};q^{\eta})_\infty}{(q^{\eta/2},q^{\eta};q^{\eta})_\infty}=\frac{(-q^{\eta/2};q^{\eta})_\infty}{(q^{\eta};q^{\eta})_\infty},
 \end{equation}
and
 \begin{equation}\label{eqn-a-1-1-1-1-1b}
(-q^{\eta(2k-\lambda)/2},q^{\eta(2k-\lambda)/2};q^{\eta(2k-\lambda)})_\infty(q^{2\eta(2k-\lambda)};q^{2\eta(2k-\lambda)})_\infty=(q^{\eta(2k-\lambda)};q^{\eta(2k-\lambda)})_\infty.
\end{equation}
Substituting \eqref{eqn-a-1-1-1-1-1a}  and \eqref{eqn-a-1-1-1-1-1b} into \eqref{eqn-a-1-1-1-1-1at}, and noting that $\alpha_{(\lambda+1)/2}=\eta/2$,
we obtain that for  $k>r\geq \lambda\geq 0$,
\begin{equation}\label{pf-thm18}
 \begin{split}
 & \sum_{n\geq0}\overline{A}_{1}(\alpha_1,\ldots,\alpha_\lambda;\eta,k,r;n)q^n\\[5pt]
  &\qquad =\frac{(-q^{\alpha_1},\ldots,-q^{\alpha_\lambda},-q^{\eta};q^{\eta})_\infty(q^{\eta(r-\frac{\lambda}{2})},q^{\eta
 (2k-r-\frac{\lambda}{2})},q^{\eta(2k-\lambda)}
 ;q^{\eta(2k-\lambda)})_\infty}{(q^\eta;q^\eta)_\infty}.
 \end{split}
 \end{equation}

When $j=0$, by definition, we have for  $k\geq r\geq \lambda\geq 0$,
 \begin{equation}\label{eqn-over-a-1-1-1-1-0}
 \begin{split}
 &\sum_{n\geq0}\overline{A}_0(\alpha_1,\ldots,\alpha_\lambda;\eta,k,r;n)q^n\\[5pt]
 &\qquad =(-q^{\alpha_1},\ldots,-q^{\alpha_{(\lambda-1)/2}},-q^{\alpha_{(\lambda+3)/2}},\ldots,-q^{\alpha_\lambda},-q^{\eta};q^{\eta})_\infty\\[5pt]
 &\qquad \qquad \times\frac{(q^{\eta(2r-\lambda)/2},
 q^{\eta(4k-2r-\lambda-2)/2},
 q^{\eta(2k-\lambda-1)};
 q^{\eta(2k-\lambda-1)})_\infty(q^{\eta};q^{2\eta})_\infty}{(q^{\eta/2};q^{\eta/2})_\infty}.
 \end{split}
 \end{equation}
Substituting  \eqref{eqn-a-1-1-1-1-1a} into \eqref{eqn-over-a-1-1-1-1-0}, we obtain that for $k\geq r\geq \lambda\geq 0$,
 \begin{equation}\label{pf-thm18aa}
 \begin{split}
 & \sum_{n\geq0}\overline{A}_{0}(\alpha_1,\ldots,\alpha_\lambda;\eta,k,r;n)q^n\\[5pt]
  &\qquad =\frac{(-q^{\alpha_1},\ldots,-q^{\alpha_\lambda},-q^{\eta};q^{\eta})_\infty(q^{\eta(r-\frac{\lambda}{2})},q^{\eta
 (2k-r-\frac{\lambda}{2}-1)},q^{\eta(2k-\lambda-1)}
 ;q^{\eta(2k-\lambda-1)})_\infty}{(q^\eta;q^\eta)_\infty}.
 \end{split}
 \end{equation}
Combining \eqref{pf-thm18} and \eqref{pf-thm18aa}, we conclude that \eqref{gf-overlineA_j-eq} holds when $\lambda$ is odd. This completes the proof. \qed

\section{The (reverse) Gordon marking and $(k-1)$-bands}

The main objective of this section is to  give a criterion to determine whether an overpartition counted by $\overline{B}_1(\alpha_1,\ldots,\alpha_\lambda;
 \eta,k,r;n)$  is also  counted by $\overline{B}_0(\alpha_1,\ldots,\alpha_\lambda;  \eta,k,r;n)$. Let $j=0$ or $1$ and let $\lambda,\ k$ and $r$    be integers such that  $k\geq r\geq \lambda\geq0$.  Let $\overline{\mathcal{B}}_j(\alpha_1,\ldots,\alpha_\lambda;\eta,k,r)$ denote the set of overpartitions  counted by $\overline{B}_j(\alpha_1,\ldots,\alpha_\lambda;
 \eta,k,r;n)$ for  $n\geq 0$.  Let $ {\mathcal{B}}_j(\alpha_1,\ldots,\alpha_\lambda;\eta,k,r)$ denote the set of partitions  counted by ${B}_j(\alpha_1,\ldots,\break$
 $\alpha_\lambda;
 \eta,k,r;n)$ for $n\geq 0$.  As   mentioned in the introduction,   we could use   $ {\mathcal{B}}_j(\alpha_1,\ldots,\alpha_\lambda;\eta,k,\break $$r)$ to denote the set of overpartitions  in $\overline{\mathcal{B}}_j(\alpha_1,\ldots,\alpha_\lambda;\eta,k,r)$  without overlined parts divisible by $\eta$.   For $\pi \in \overline{\mathcal{B}}_j(\alpha_1,\ldots,\alpha_\lambda;\eta,k,r)$, we call $\pi$ a $\overline{B}_j$-overpartition for short.

The Gordon marking of an ordinary partition  was  introduced by Kur\c{s}ung\"oz in \cite{Kursungoz-2010a,Kursungoz-2010}. Kim \cite{Kim-2018}  introduced   the Gordon marking of an ordinary partition in  ${\mathcal{B}}_1(\alpha_1,\ldots,\alpha_\lambda;\break$ $\eta,k,r)$, which generalizes the Gordon marking of an ordinary partition.  The Gordon marking of an overpartition was   defined by  Chen, Sang and Shi \cite{Chen-Sang-Shi-2013}.   Now  we  define   the Gordon marking  of a $\overline{B}_1$-overpartition. Bear in mind that the parts in an overpartition are ordered as follows:
\begin{equation*}
1<\bar{1}<2<\bar{2}<\cdots.
\end{equation*}
 For positive integers $t$ and $b$,  we define $t\pm b$ (resp. $\overline{t}\pm b$) as a non-overlined part of size  $t\pm b$ (resp. an overlined part of size ${t\pm b}$).

\begin{defi}[Gordon marking]\label{Gordon-marking}
For  $k\geq r\geq \lambda\geq0$, let $\pi=(\pi_1,\pi_2,\ldots,\pi_\ell)$ be an overpartition satisfying {\rm{(1)}} and {\rm{(2)}} in Definition \ref{defi-O-B}. Assign a positive integer to each part of $\pi$ as follows{\rm{:}} First,  assign $1$ to  $\pi_\ell$. Then, for each $\pi_i$, assign $s$  to $\pi_i$, where $s$ is the smallest positive integer that is not used to mark the parts  $\pi_m$ such that $m>i$ and $\pi_m\geq \pi_i-\eta$ with strict inequality if $\pi_i$ is overlined. Denote the Gordon marking of $\pi$ by $G(\pi)$.
\end{defi}

It can be seen that for each $\pi_i$, the mark of $\pi_i$ is the smallest positive integer that is not used to mark the parts after $\pi_i$ in $[\pi_i-\eta, \pi_i]$ (resp. $(\pi_i-\eta, \pi_i)$) if $\pi_i$ is non-overlined (resp. overlined). Assume that $\pi_\ell$ is assigned with 1. Then the part $\pi_m$ of $\pi$ is in $[\pi_i-\eta, \pi_i]$ (resp. $(\pi_i-\eta, \pi_i)$) means that $\pi_i-\eta\leq \pi_m\leq \pi_i$ (resp. $\pi_i-\eta< \pi_m< \pi_i$). For notational convenience, we use $I(\pi_i-\eta, \pi_i)$ to denote   the interval $[\pi_i-\eta, \pi_i]$  if $\pi_i$ is non-overlined,
or the interval  $(\pi_i-\eta, \pi_i)$ if $\pi_i$ is overlined.

For example, let  $\pi$ be an overpartition in $\overline{\mathcal{B}}_1(1,5,9;10,5,4)$  given by
\begin{equation}\label{mark-exa-1}
\begin{split}
&\pi=(\overline{80},80,80,\overline{70},70,\overline{69},\overline{60},
60,\overline{55},\overline{51},{50},\overline{49},\overline{45},\overline{41},\overline{39},\overline{35},\\[3pt]
&\ \ \ \ \ \ \ \overline{29},\overline{20},{20},20,\overline{11},\overline{10},\overline{9},
\overline{5},\overline{1}).
\end{split}
\end{equation}
The Gordon marking of $\pi$ is given by
\begin{equation}\label{exa-g-1}
\begin{split}
&G(\pi)=(\overline{80}_1,{ {80}_4},{80}_2,\overline{70}_1,{70}_3,\overline{69}_2,{ \overline{60}_4},
{60}_1,\overline{55}_2,\overline{51}_3,{ {50}_4},\overline{49}_1,\overline{45}_2,\overline{41}_3,\overline{39}_1,\overline{35}_2,\\[3pt]
&\ \ \ \ \ \ \ \ \ \ \ \  \overline{29}_1,{ \overline{20}_4},{20}_3,{20}_2,\overline{11}_1,{ \overline{10}_4},\overline{9}_3,
\overline{5}_2,\overline{1}_1),
\end{split}
\end{equation}
where the subscript of each part  represents the mark in the Gordon marking.

For  $k\geq r\geq \lambda\geq0$, let $\pi=(\pi_1,\pi_2,\ldots,\pi_\ell)$ be an overpartition satisfying {\rm{(1)}} and {\rm{(2)}} in Definition \ref{defi-O-B}.  If  the condition  {\rm(3)} in Definition \ref{defi-O-B} is also fulfilled, that is, for $1\leq i\leq \ell-k+1$, $\pi_i\geq \pi_{i+k-1}+\eta$ with strict inequality if $\pi_i$ is non-overlined, then  for each $\pi_i$, the number of  parts  after $\pi_i$ belonging to $I(\pi_i-\eta, \pi_i)$   is at most $k-2$, so the marks in the  Gordon marking of $\pi$ do not  exceed $k-1$.   For the overpartition $\pi$ in $\overline{\mathcal{B}}_1(1,5,9;10,5,4)$  defined in \eqref{mark-exa-1}, by \eqref{exa-g-1}, we see that the largest mark in $G(\pi)$ is $4$.

If we assign a mark to each part starting with the largest part instead, then the resulting marking will be called the reverse Gordon marking.

\begin{defi}[Reverse Gordon marking]\label{R-Gordon-marking}
For  $k\geq r\geq \lambda\geq0$, let $\pi=(\pi_1,\pi_2,\ldots,\pi_\ell)$ be an overpartition satisfying {\rm{(1)}} and {\rm{(2)}} in Definition \ref{defi-O-B}. Assign a positive integer to each part of $\pi$ as follows{\rm{:}} First assign $1$ to $\pi_1$. Then, for each $\pi_i$, assign $s$ to $\pi_i$, where $s$ is   the smallest positive integer that is not used to mark the parts  $\pi_m$  such that $m<i$ and $\pi_m\leq \pi_i+\eta$ with strict inequality if $\pi_i$ is overlined. Denote the reverse Gordon marking of $\pi$ by $RG(\pi)$.
\end{defi}

Analogously, for each $\pi_i$, the mark of $\pi_i$ is  the smallest positive integer that is not used to mark the parts  before $\pi_i$ belonging to $I(\pi_i, \pi_i+\eta)$. Furthermore, for  $\pi$ in $\overline{\mathcal{B}}_1(\alpha_1,\ldots,\alpha_\lambda;\eta,k,r)$,  the marks in the reverse Gordon marking of $\pi$ do not  exceed $k-1$.

For the  overpartition $\pi$ in $\overline{\mathcal{B}}_1(1,5,9;10,5,4)$  defined in \eqref{mark-exa-1}, the reverse Gordon marking  of $\pi$ reads
\begin{equation*}\label{exa-r-1}
\begin{split}
&RG(\pi)=(\overline{80}_1,{80}_2,{80}_3,\overline{70}_1,{ {70}_4},\overline{69}_2,\overline{60}_1,
{60}_3,\overline{55}_2,{ \overline{51}_4},{50}_1,\overline{49}_3,\overline{45}_2,{ \overline{41}_4},\overline{39}_1,\overline{35}_2,\\[3pt]
&\ \ \ \ \ \ \ \ \ \ \ \ \ \  \overline{29}_1,\overline{20}_2,{20}_3,{ {20}_4},\overline{11}_1,\overline{10}_2,\overline{9}_3,{ \overline{5}_4},\overline{1}_1),
\end{split}
\end{equation*}
from which we see that the largest mark in $RG(\pi)$ is $4$.

We proceed to give a  criterion   to determine whether a $\overline{B}_1$-overpartition  is also a $\overline{B}_0$-overpartition.
Let $\pi=(\pi_1,\pi_2,\ldots,\pi_\ell)$ be an overpartition in $\overline{\mathcal{B}}_1(\alpha_1,\ldots,\alpha_\lambda;\eta,k,r)$.
If there are no  $k-1$ consecutive parts  $\pi_i,\,\pi_{i+1},\ldots, \pi_{i+k-2}$ in $\pi$ such that
 \begin{equation}\label{sequence}
\pi_i\leq\pi_{i+k-2}+\eta\text{ with strict inequality if }
\\ \pi_i \text{ is overlined},
\end{equation}
then by Definition \ref{defi-O-B}, we see that $\pi$ is also in  $\overline{\mathcal{B}}_0(\alpha_1,\ldots,\alpha_\lambda;\eta,k,r)$.
Assume that there exist $k-1$ consecutive parts  $\pi_i,\,\pi_{i+1},\ldots, \pi_{i+k-2}$ in $\pi$ satisfying \eqref{sequence}.
By definition, $\pi$ is not  in  $\overline{\mathcal{B}}_0(\alpha_1,\ldots,\alpha_\lambda;\eta,k,r)$ if
 \begin{equation*}
 \left[\pi_i/\eta\right]+\cdots+\left[\pi_{i+k-2}/\eta\right]\equiv r+\overline{V}_{\pi}(\pi_i)\pmod{2}.
 \end{equation*}
It follows that $\pi$ is   an overpartition in  $\overline{\mathcal{B}}_0(\alpha_1,\ldots,\alpha_\lambda;\eta,k,r)$ if and only if for any $k-1$  consecutive parts   $\pi_i,\,\pi_{i+1},\ldots, \pi_{i+k-2}$ in $\pi$ satisfying \eqref{sequence}, we have
\begin{equation}\label{congruence}
 \left[\pi_i/\eta\right]+\cdots+\left[\pi_{i+k-2}/\eta\right]\equiv r-1+\overline{V}_{\pi}(\pi_i)\pmod{2}.
 \end{equation}

 The above $k-1$ consecutive parts satisfying   \eqref{sequence} will be called  {\it a $(k-1)$-band of $\pi$}  in this sense that the difference between the largest element and the smallest element in a $(k-1)$-band is at most $\eta$.  For the $(k-1)$-band $\{\pi_{i+l}\}_{0\leq l\leq k-2}$, if  $\{\pi_{i+l}\}_{0\leq l\leq k-2}$ satisfy the congruence condition \eqref{congruence},  then we say that the $(k-1)$-band $\{\pi_{i+l}\}_{0\leq l\leq k-2}$ is even. Otherwise, we say that it is odd.

 For example,  let $\pi$ be the overpartition in $\overline{\mathcal{B}}_1(1,5,9;10,5,4)$ defined in \eqref{mark-exa-1}, where $k=5$. There are twelve $4$-bands in  $\pi$. It can be checked that all of them are even.

  \[\{80,{ 80},\overline{70},{ 70}\},\{{ 70},\overline{69},{ \overline{60}},60\},
  \{{ \overline{60}},{60},\overline{55},{ \overline{51}}\},
  \{{60},\overline{55},{ \overline{51}},{ {50}}\},\]
   \[\{\overline{55},{ \overline{51}},{ {50}},\overline{49}\},
   \{{ \overline{51}},{ {50}},\overline{49},\overline{45}\},\{{{50}},\overline{49},
   \overline{45},{ \overline{41}}\},\{\overline{29},{ \overline{20}},{20},{ 20}\},\]
   \[\{{ \overline{20}},{20},{ 20},\overline{11}\},\{20,{ 20},\overline{11},{ \overline{10}}\},
   \{\overline{11},{ \overline{10}},\overline{9},{ \overline{5}}\},\{{ \overline{10}},\overline{9},{ \overline{5}},\overline{1}\}.\]

 For the overpartitions  $\pi$  in $\overline{\mathcal{B}}_1(\alpha_1,\ldots,\alpha_\lambda;\eta,k,r)$,  we see that     $\pi$ is  an overpartition in $\overline{\mathcal{B}}_0(\alpha_1,\ldots,\alpha_\lambda;\eta,k,r)$ if and only if all $(k-1)$-bands  of $\pi$ are even. For  $k\geq r\geq \lambda\geq0$, let $\pi$ be an overpartition in $\overline{\mathcal{B}}_1(\alpha_1,\ldots,\alpha_\lambda;\eta,k,r)$.
For each $(k-1)$-band $\{\pi_{i+l}\}_{0\leq l\leq k-2}$ of $\pi$,
it is easy to see that the marks of $\pi_{i+l}$ are distinct in
the Gordon marking and the reverse Gordon marking of  $\pi$.
Hence there exists one part in  $\{\pi_{i+l}\}_{0\leq l\leq k-2}$ marked with $k-1$ in the Gordon marking and the reverse Gordon marking of  $\pi$.

Next we show that we may restrict our attention
to certain special $(k-1)$-bands  to determine
whether $\pi$ is  an overpartition in $\overline{\mathcal{B}}_0(\alpha_1,\ldots,\alpha_\lambda;\eta,k,r)$.
Such special $(k-1)$-bands will be classified into two kinds depending on the positions of the $(k-1)$-marked parts of the Gordon marking or the reverse Gordon marking in the $(k-1)$-bands.

{\noindent \it The  $(k-1)$-bands of the first kind} will be concerned with the case in which  the $(k-1)$-marked part in the Gordon marking
 is the largest element in the band.
 Assume that there are $N$ parts marked with $k-1$ in $G(\pi)$,
and denote these $(k-1)$-marked parts by $\tilde{g}_1(\pi)> {\tilde{g}}_2(\pi)>\cdots >\tilde{g}_N(\pi)$. For each $(k-1)$-marked part $\tilde{g}_p(\pi)$
in $G(\pi)$, the number of the  $(k-1)$-bands of $\pi$ including $\tilde{g}_p(\pi)$ is at most $k-1$. We claim  that there is a $(k-1)$-band of $\pi$ such that $\tilde{g}_p(\pi)$ is the largest element of this $(k-1)$-band.
Assume that $\tilde{g}_p(\pi)$ is the $s$-th part $\pi_s$ of $\pi$.
By Definition \ref{Gordon-marking}, we deduce that there exist $k-2$
parts $\pi_m$  such that $m>s$  and $\pi_m\geq \pi_s-\eta$ with strict
inequality if $\pi_s$ is overlined. This implies  that
$\pi_s,\pi_{s+1},\ldots,\pi_{s+k-2}$ in $\pi$ satisfy  \eqref{sequence}, that is,
$\{\pi_{s+l}\}_{0\leq l\leq k-2}$ is a $(k-1)$-band of $\pi$. Furthermore, the $(k-1)$-marked part $\tilde{g}_p(\pi)$  is the largest
element of this $(k-1)$-band. So the claim is proved. {\it Such a $(k-1)$-band is called the $(k-1)$-band
induced by $\tilde{g}_p(\pi)$}, {\it denoted  $\{\tilde{g}_p(\pi)\}_{k-1}$}.
Obviously,  $\{\tilde{g}_1(\pi)\}_{k-1}, \{\tilde{g}_2(\pi)\}_{k-1},\ldots,\{\tilde{g}_N(\pi)\}_{k-1}$ of $\pi$
are disjoint.

For example, for   the  overpartition   $\pi$ given in \eqref{mark-exa-1},
  there are five $4$-marked parts in $G(\pi)$, namely,
 $\tilde{g}_1(\pi)=80$, $\tilde{g}_2(\pi)=\overline{60}$,
 $\tilde{g}_3(\pi)={50}$, $\tilde{g}_4(\pi)=\overline{20}$ and $\tilde{g}_5(\pi)=\overline{10}$.  The   $4$-bands induced by  $\tilde{g}_1(\pi),\tilde{g}_2(\pi),\tilde{g}_3(\pi),
 \tilde{g}_4(\pi)$  and $\tilde{g}_5(\pi)$   are  illustrated
 in $G(\pi)$ below:
\begin{equation*} \label{mark-exa-1-s}
\begin{split}
&G(\pi)=(\overline{80}_1,\overbrace{{ {80}_4},
{ {80}_2,
\overline{70}_1,{70}_3}}^{{\{80\}_4}},\overline{69}_2,\overbrace{{ \overline{60}_4},{ {60}_1,\overline{55}_2,\overline{51}_3}}^{{\{\overline{60}\}_4}},\overbrace{{ {50}_4},{ \overline{49}_1,\overline{45}_2,\overline{41}_3}}^{{\{{50}\}_4}},\overline{39}_1,\overline{35}_2,\\[5pt]
&\ \ \ \ \ \ \ \ \ \ \
\overline{29}_1,\underbrace{{ \overline{20}_4},{ {20}_3,
{20}_2,\overline{11}_1}}_
{\{\overline{20}\}_4},\underbrace{{ {\overline{10}_4}},
{ \overline{9}_3,\overline{5}_2,\overline{1}_1}}_{\{\overline{10}\}_4}).
\end{split}
\end{equation*}

 We now consider {\it the  $(k-1)$-bands of the second kind} under the condition that the $(k-1)$-marked part in the reverse Gordon marking is
the smallest part in the band. Assume that there are $M$ parts marked with $k-1$ in $RG(\pi)$, namely,  $\tilde{r}_1(\pi)> \tilde{{r}}_2(\pi)>\cdots >\tilde{r}_M(\pi)$.
By the same reasoning, we see that  there is a $(k-1)$-band of $\pi$ in which $\tilde{r}_p(\pi)$ is the smallest element.  {\it Such a $(k-1)$-band is called the $(k-1)$-band induced by $\tilde{r}_p(\pi)$}, {\it denoted  $\{\tilde{r}_p(\pi)\}_{k-1}$}. Clearly, $\{\tilde{r}_1(\pi)\}_{k-1}, \{\tilde{r}_2(\pi)\}_{k-1},\ldots, \{\tilde{r}_M(\pi)\}_{k-1}$ of $\pi$ are disjoint.

 For example, for the overpartition  $\pi$    given in \eqref{mark-exa-1},   there are five $4$-marked parts in $RG(\pi)$, which are  $\tilde{r}_1(\pi)=70$, $\tilde{r}_2(\pi)=\overline{51}$, $\tilde{r}_3(\pi)=\overline{41}$, $\tilde{r}_4(\pi)={20}$  and $\tilde{r}_5(\pi)=\overline{5}$. The   $4$-bands induced by  $\tilde{r}_1(\pi)$, $\tilde{r}_2(\pi)$, $\tilde{r}_3(\pi)$, $\tilde{r}_4(\pi)$ and  $\tilde{r}_5(\pi)$ are  displayed  below:
\begin{equation*} \label{rmark-exa-1-s}
\begin{split}
&RG(\pi)=(\overline{80}_1,\overbrace{{ {80}_2,{80}_3,\overline{70}_1},
{ {70}_4}}^{{\{70\}_4}},\overline{69}_2,
\overbrace{{ \overline{60}_1,{60}_3,\overline{55}_2}, { \overline{51}_4}}^{{\{\overline{51}\}_4}},\overbrace{{{50}_1,\overline{49}_3,\overline{45}_2}, { \overline{41}_4}}^{{\{\overline{41}\}_4}}, \overline{39}_1,\overline{35}_2,\\[5pt]
&\ \ \ \ \ \ \ \ \ \ \ \ \ \ \underbrace{{ \overline{29}_1,\overline{20}_2,} { {20}_3},{ {20}_4}}
_{\{20\}_4},
\underbrace{{ \overline{11}_1,\overline{10}_2,\overline{9}_3},
{ \overline{5}_4}}_{\{\overline{5}\}_4},\overline{1}_1).
\end{split}
\end{equation*}

It remains to show that  the  $(k-1)$-bands  of $\pi$ induced by the $(k-1)$-marked parts in  $G(\pi)$ or $RG(\pi)$ are enough to  determine whether $\pi$ is  an overpartition in  $\overline{\mathcal{B}}_0(\alpha_1,\ldots,\alpha_\lambda;\eta,\break $$ k,r)$. For this purpose, we need the following  property relating $G(\pi)$ and $RG(\pi)$.

\begin{prop}\label{sequence-length}
For  $k\geq r\geq \lambda\geq0$, let $\pi$ be an overpartition in $\overline{\mathcal{B}}_1
 (\alpha_1,\ldots,\alpha_\lambda;\eta,k,r)$. Assume that there are $N$ parts marked with $k-1$ in the Gordon marking of $\pi$, say, $\tilde{g}_1(\pi)> \tilde{g}_2(\pi)>\cdots >\tilde{g}_N(\pi)$, and there are  $M$ parts marked with $k-1$ in the reverse Gordon marking of $\pi$, say, $\tilde{r}_1(\pi)> \tilde{r}_2(\pi)>\cdots >\tilde{r}_M(\pi)$. Then $N=M$.   Moveover, for each $1\leq i \leq N$, we have $\tilde{g}_i(\pi)\in \{\tilde{r}_i(\pi)\}_{k-1}$ and $\tilde{r}_i (\pi)\in \{\tilde{g}_i(\pi)\}_{k-1}$, where $\{\tilde{g}_i(\pi)\}_{k-1}$ {\rm(}resp. $\{\tilde{r}_i(\pi)\}_{k-1}${\rm)} is the $(k-1)$-band of $\pi$ induced by $\tilde{g}_i(\pi)$ {\rm(}resp. $\tilde{r}_i(\pi)${\rm)}.
\end{prop}

\pf  For $N=0$, there are no $(k-1)$-marked parts  in $G(\pi)$, and so there are no $(k-1)$-bands  in $\pi$. This implies that there are no  $(k-1)$-marked parts   in  $RG(\pi)$. It follows that  $M=0$. Conversely, if $M=0$, then  $N=0$.

We next consider the case $M,N>0$.   We first prove that $M\geq N$.  For each fixed $(k-1)$-marked part $\tilde{g}_i(\pi)$ in $G(\pi)$, where $1\leq i\leq N$,  assume that $\tilde{g}_i(\pi)$  is the $g_i$-th part of $\pi=(\pi_1,\pi_2,\ldots,\pi_\ell)$, that is, $\tilde{g}_i(\pi)=\pi_{g_i}$. Since $\pi_{g_i}$ is the largest element in the $(k-1)$-band induced by $\pi_{g_i}$, we find that the parts
\[  \pi_{g_i}\geq \pi_{g_i+1} \geq \cdots \geq   \pi_{g_i+k-2} \]
are  in the $(k-1)$-band of $\pi$ induced by $\pi_{g_i}$. Moreover, the marks of these parts in $RG(\pi)$ are distinct. It follows that there exists $t_i$ such that $0\leq t_i\leq k-2$ and $ \pi_{g_i+t_i}$ is marked with $k-1$ in $RG(\pi)$. Since the $(k-1)$-bands  $\{\tilde{g}_1(\pi)\}_{k-1}, \{\tilde{g}_2(\pi)\}_{k-1},\ldots, \{\tilde{g}_N(\pi)\}_{k-1}$ of $\pi$ are disjoint,   the parts $ \pi_{g_1+t_1}$, $ \pi_{g_2+t_2}$,\ldots, $ \pi_{g_N+t_N}$ are distinct, which are marked with $k-1$ in $RG(\pi)$. This means that $M\geq N$.  A similar argument   yields   $N\geq M$.  We conclude that $M=N$. Note that the above proof also indicates that $\tilde{r}_i (\pi)=\pi_{g_i+t_i}$, which implies that $\tilde{r}_i (\pi) \in \{\tilde{g}_i(\pi)\}_{k-1}$ for $1\leq i\leq N$. Similarly,
 $\tilde{g}_i(\pi)\in \{\tilde{r}_i(\pi)\}_{k-1}$ for $1\leq i\leq N$, and thus the proof is complete.   \qed

 For example, for  the overpartition $\pi$  given in  \eqref{mark-exa-1}, there are five $4$-marked parts in $G(\pi)$, and in the meantime there are five  $4$-marked parts in $RG(\pi)$.

 We are now in a position to present the main result of this section.

  \begin{thm}\label{parity k-1 sequence-over-g}
  For $k\geq r\geq \lambda\geq0$, $k-1>\lambda$ and $N\geq0$, let $\pi=(\pi_1,\pi_2,\ldots, \pi_\ell)$ be an overpartition in ${{\overline{\mathcal{B}}}}_1(\alpha_1,\ldots,\alpha_\lambda;\eta,k,r)$  with $N$ parts marked with $k-1$ in $G(\pi)$ {\rm{(}}resp. $RG(\pi)${\rm{)}}, say $ \tilde{g}_1(\pi)>  \tilde{g}_2(\pi)>\cdots > \tilde{g}_N(\pi)$ {\rm{(}}resp. $ \tilde{r}_1(\pi)>  \tilde{r}_2(\pi)>\cdots > \tilde{r}_N(\pi)${\rm{)}}.   Then  $\pi$ is an overpartition in $\overline{\mathcal{B}}_0(\alpha_1,\ldots,\alpha_\lambda;\eta,k,r)$ if and only if for all   $1\leq i\leq N$, $\{\tilde{g}_i(\pi)\}_{k-1}$ {\rm{(}}resp. $\{\tilde{r}_i(\pi)\}_{k-1}${\rm{)}}  are even. In particular, for $\lambda=k-1$, the assertion holds if   there are no  overlined parts divisible by $\eta$ in $\pi$.
  \end{thm}

 For example,  the  overpartition $\pi$  given in  \eqref{mark-exa-1} is also an overpartition in    $\overline{\mathcal{B}}_0(1,5,9;10,5,\break$$4)$.  To prove Theorem \ref{parity k-1 sequence-over-g}, we need  the following lemma.

\begin{lem}\label{parity k-1 sequence-over-old} For $k\geq r\geq \lambda\geq0$ and $k-1>\lambda$, let $\pi=(\pi_1,\pi_2,\ldots, \pi_\ell)$ be an overpartition   in ${\mathcal{\overline{B}}}_1(\alpha_1,\ldots,\alpha_\lambda;\eta,k,r)$, and let  $\{\pi_{c+l}\}_{0\leq l\leq k-2}$ and $\{\pi_{d+l}\}_{0\leq l\leq k-2}$   be two $(k-1)$-bands  of $\pi$. If $\pi_c>\pi_d$ and $\pi_c\leq \pi_{d+k-2}+2\eta$ with strict inequality if $\pi_c$ is overlined, then    $\{\pi_{c+l}\}_{0\leq l\leq k-2}$ and $\{\pi_{d+l}\}_{0\leq l\leq k-2}$  are of the same parity. In particular, for $\lambda=k-1$, the assertion holds if   there are no  overlined parts divisible by $\eta$ in $\pi$.
  \end{lem}

 The above lemma enables us to  establish  the following proposition, which, together with Proposition \ref{sequence-length},  leads to Theorem \ref{parity k-1 sequence-over-g}.

\begin{prop} \label{parity k-1 sequence-over}   For $k\geq r\geq \lambda\geq0$ and $k-1>\lambda$, let $\tilde{g}_p(\pi)$ be a $(k-1)$-marked part  in the Gordon marking of an overpartition $\pi$ in ${{\overline{\mathcal{B}}}}_1(\alpha_1,\ldots,\alpha_\lambda;\eta,k,r)$. Then  the $(k-1)$-bands  of $\pi$ including $\tilde{g}_p(\pi)$ and the $(k-1)$-band induced by $\tilde{g}_p(\pi)$ are of the same parity.  In particular, for $\lambda=k-1$, the assertion holds if   there are no  overlined parts divisible by $\eta$ in $\pi$.
\end{prop}

 For example, for  the  overpartition $\pi$ given in  \eqref{mark-exa-1},  the $4$-bands induced by $\tilde{g}_1(\pi)$, $\tilde{g}_2(\pi)$,
 $\tilde{g}_3(\pi)$, $\tilde{g}_4(\pi)$ and $\tilde{g}_5(\pi)$ are all even. Moreover,   for any $1\leq p\leq 5$, the $(k-1)$-bands  including $\tilde{g}_p(\pi)$ and the $(k-1)$-band  induced by $\tilde{g}_p(\pi)$ are of the same parity.

 It should be mentioned that for $\lambda=k-1$, if there is an overlined part divisible by $\eta$ in $\pi$, then a $(k-1)$-band  of $\pi$ including the $(k-1)$-marked part  $\tilde{g}_p(\pi)$   may have a  different parity from that of the $(k-1)$-band induced by   $\tilde{g}_p(\pi)$. For example, let  $\pi=(\overline{21},\overline{20},\overline{15},\overline{11},\overline{9},\overline{5})$ be an overpartition in ${{\overline{\mathcal{B}}}}_1(1,5,9;10, 4,3)$, where $\lambda=3$, $k=4$,  $\eta=10$ and $r=3$. The Gordon marking of $\pi$ is
\[G(\pi)=(\overline{21}_3,\overline{20}_2,\overline{15}_1,\overline{11}_3,\overline{9}_2,\overline{5}_1),\]
from which we see that there are two $3$-marked parts $\tilde{g}_1(\pi)=\overline{21}$ and $\tilde{g}_2(\pi)=\overline{11}$, and so we get   the $3$-band $\{\overline{11},\overline{9},\overline{5}\}$  induced by $\tilde{g}_2(\pi)=\overline{11}$ along with a $3$-band $\{\overline{20},\overline{15},\overline{11}\}$  including $\tilde{g}_2(\pi)=\overline{11}$. Apparently, the $3$-band $\{\overline{11},\overline{9},\overline{5}\}$ is even, since
\[[\overline{11}/10]+[\overline{9}/10]+[\overline{5}/10] = 1\equiv  r-1+\overline{V}_{\pi}(\overline{11})\pmod{2},\]
and the $3$-band $\{\overline{20},\overline{15},\overline{11}\}$ is odd, since
\[[\overline{20}/10]+[\overline{15}/10]+[\overline{11}/10] =4\equiv  r+\overline{V}_{\pi}(\overline{20})\pmod{2}.\]

In the remainder of this section, we  present the  proofs of Lemma \ref{parity k-1 sequence-over-old} and  Proposition  \ref{parity k-1 sequence-over}.

 \noindent{\it Proof of Lemma \ref{parity k-1 sequence-over-old}.} To show that $\{\pi_{d+l}\}_{0\leq l\leq k-2}$ and  $\{\pi_{c+l}\}_{0\leq l\leq k-2}$ are of the same parity, we write
  \begin{equation}\label{cong-x-lem-over}
  [\pi_{d+k-2}/\eta]+\cdots+[\pi_{d}/\eta]\equiv a+\overline{V}_\pi (\pi_{d})\pmod2,
  \end{equation}
   where $a=r-1$ or $r$. Trivially, $\{\pi_{d+l}\}_{0\leq l\leq k-2}$ is even  when $a=r-1$, and $\{\pi_{d+l}\}_{0\leq l\leq k-2}$ is odd when $a=r$.

   We intend to prove that
    \begin{equation}\label{cong-g-lem-over}
  [\pi_{c+k-2}/\eta]+\cdots+[\pi_{c}/\eta]\equiv a+\overline{V}_\pi(\pi_{c})\pmod2,
   \end{equation}
   where $a$ is given as  in \eqref{cong-x-lem-over}.   Since  $\{\pi_{d+l}\}_{0\leq l\leq k-2}$ and $\{\pi_{c+l}\}_{0\leq l\leq k-2}$  are $(k-1)$-bands  of $\pi$,   we have
 \begin{equation}\label{sequence-lemm1}
\begin{split}
&\ \ \ \ \ \ \ \ \ \ \ \ \ \ \ \ \ \ \ \ \
\pi_d\geq \pi_{d+1}\geq \cdots \geq \pi_{d+k-2},\\[5pt]
& \text{where }\pi_d\leq\pi_{d+k-2}+\eta\text{ with strict inequality if }
 \pi_d\text{ is overlined, and}
\end{split}
\end{equation}
 \begin{equation}\label{sequence-lemm2}
 \begin{split}
&\ \ \ \ \ \ \ \ \ \ \ \ \ \ \ \ \ \ 
 \pi_c\geq \pi_{c+1}\geq \cdots \geq \pi_{c+k-2},\\[5pt]
&\text{where }\pi_c\leq\pi_{c+k-2}+\eta\text{ with strict inequality if } \pi_c \text{ is overlined}.
\end{split}
\end{equation}
Under the condition $\pi_c>\pi_d$,   we have  $c<d$. Assume that $d=c+t$ where $t\geq1$. Given that $\pi=(\pi_1,\pi_2,\ldots, \pi_\ell)$ is an overpartition in ${{\overline{\mathcal{B}}}}_1(\alpha_1,\ldots,\alpha_\lambda;\eta,k,r)$,   for $1\leq i\leq \ell-k+1$,
\begin{equation}\label{sequence-lempa3}
\pi_i\geq\pi_{i+k-1}+\eta\text{ with strict inequality if }
\\ \pi_i \text{ is non-overlined}.
 \end{equation}
It follows that there are at most $2k-2$ parts of $\pi$ belonging to $I(\pi_c-2\eta,\pi_c)$. Therefore, by $\pi_{d+k-2}\geq \pi_c-2\eta$,
we deduce that $1\leq t\leq k-1$, and so for $1\leq t \leq l \leq k-2$,
\begin{equation}\label{relation-00}
   \pi_{c+l}=\pi_{d+l-t}.
\end{equation}
Combining \eqref{sequence-lemm1} and \eqref{sequence-lemm2},  we find that $\pi_{c+t}=\pi_d\leq \pi_{d+k-2}+\eta$ with strict inequality if $\pi_{d}$ is overlined.
Noting that $(d+k-2)-(c+t-1)=k-1$,  using \eqref{sequence-lempa3}, we obtain that $\pi_{c+t-1}\geq\pi_{d+k-2}+\eta$ with strict inequality if $\pi_{c+t-1}$ is non-overlined. It follows that  $\pi_{c+t}<\pi_{c+t-1}$.   The same argument yields  $\pi_{d+k-1-t}<\pi_{d+k-2-t}$.  To summarize, the overlapping structure of  $\{\pi_{c+l}\}_{0\leq l\leq k-2}$ and  $\{\pi_{d+l}\}_{0\leq l\leq k-2}$ can be described as follows, depending on two cases.

For $1\leq t< k-1$, we  have
{\footnotesize\[\begin{array}{ccccccccccccccccccc}
\pi_{d+k-2}&\leq &\cdots
&\leq& \pi_{d+k-1-t}&<&\pi_{d+k-2-t}&\leq&\cdots&\leq& \pi_d\\[5pt]
&&&&&&\|&&&&\|&&&&\\[5pt]
&&&&&&\pi_{c+k-2} &\leq& \cdots &\leq &\pi_{c+t}&<&\pi_{c+t-1}&\leq
&\cdots&\leq&\pi_{c}.
\end{array}
\]}
For $t=k-1$, we have
\[\pi_{d+k-2}\leq\cdots\leq \pi_d<\pi_{c+k-2}\leq\cdots\leq \pi_c.\]

We are now ready to prove \eqref{cong-g-lem-over}. By \eqref{relation-00}, we have
      \begin{eqnarray*}
   &&[\pi_{c+k-2}/\eta]+\cdots+[\pi_{c+t}/\eta]+[\pi_{c+t-1}/\eta]+\cdots
   +[\pi_{c}/\eta]\\[5pt]
   &&\quad =
   [\pi_{d+k-2-t}/\eta]+\cdots+[\pi_d/\eta]+[\pi_{c+t-1}/\eta]+\cdots
   +[\pi_{c}/\eta]\\[5pt]
   && \quad=
   [\pi_{d+k-2}/\eta]+\cdots+[\pi_d/\eta]\\[5pt]
   &&\qquad \quad +
  [\pi_{c+t-1}/\eta]+\cdots
   +[\pi_{c}/\eta]-
   \left([\pi_{d+k-2}/\eta]+\cdots+[\pi_{d+k-1-t}/\eta]\right),
   \end{eqnarray*}
 and by \eqref{cong-x-lem-over}, we find that in order to show \eqref{cong-g-lem-over}, it suffices to show that
   \begin{eqnarray} \label{cong-g-x-f-over}
   &&[\pi_{c+t-1}/\eta]+\cdots
   +[\pi_{c}/\eta]-
   \left([\pi_{d+k-2}/\eta]+\cdots+[\pi_{d+k-1-t}/\eta]\right)\nonumber\\[5pt]
   &&\quad \equiv \overline{V}_\pi(\pi_c)-\overline{V}_\pi(\pi_d) \pmod{2}.
   \end{eqnarray}

   We consider the following two cases.

Case 1: $\pi_c$ is non-overlined. In this case,  $\pi_c$ is divisible by $\eta$, so we may write $\pi_c=(b+2)\eta$.  In view of the condition that $\pi_{d+k-2}\geq \pi_c-2\eta=b\eta$, together with \eqref{sequence-lempa3}, we find that $\pi_{d+k-1-t}<\pi_c-\eta=(b+1)\eta$ and $\pi_{c+t-1}>(b+1)\eta$. Hence
\[b\eta \leq \pi_{d+k-2}\leq   \cdots \leq \pi_{d+k-1-t}<(b+1)\eta,\] and
\[(b+1)\eta<\pi_{c+t-1}\leq  \cdots \leq \pi_{c}=(b+2) \eta.\]
This implies that for $k-1-t\leq l\leq k-2$,
$[\pi_{d+l}/\eta]=b$,  and  for  $0\leq l\leq t-1$, $[\pi_{c+l}/\eta]=b+1$  if $\pi_{c+l}$ is overlined, or $[\pi_{c+l}/\eta]=b+2$ if $\pi_{c+l}$ is non-overlined. Consequently,
 \begin{eqnarray*}
   &&[\pi_{c+t-1}/\eta]+\cdots
   +[\pi_{c}/\eta]-
   \left([\pi_{d+k-2}/\eta]+\cdots+[\pi_{d+k-1-t}/\eta]\right)\nonumber\\[5pt]
    &&\quad \equiv   \overline{V}_\pi(\pi_c)-\overline{V}_\pi(\pi_d) \pmod{2},
   \end{eqnarray*}
  and so \eqref{cong-g-x-f-over} is confirmed.

Case 2: $\pi_c$ is overlined. Set   $\alpha_0=0$. Then we may write  $\pi_c=\overline{(b+1)\eta+\alpha_s}$, where $0\leq s\leq \lambda$.    Using the condition that $\pi_{d+k-2}> \pi_c-2\eta=\overline{(b-1)\eta+\alpha_s}$ and the
relation \eqref{sequence-lempa3}, we deduce that $\pi_{c+t-1}\geq\pi_{d+k-2}+\eta>\overline{b\eta+\alpha_s}$ and $\pi_{d+k-1-t}\leq\pi_c-\eta= \overline{b\eta+\alpha_s}$,  so that
 \begin{equation}\label{lemma2.6-pfa}
     \overline{(b-1)\eta+\alpha_s}<\pi_{d+k-2}\leq   \cdots \leq \pi_{d+k-1-t}\leq \overline{b\eta+\alpha_s},
 \end{equation}
  and
  \begin{equation}\label{lemma2.6-pfb}
      \overline{b\eta+\alpha_s}<\pi_{c+t-1}\leq   \cdots \leq \pi_{c}=\overline{(b+1)\eta+\alpha_s}.
  \end{equation}
 Assume that there are $f_1$ parts $\pi_{d+l}$ in \eqref{lemma2.6-pfa} satisfying $\overline{(b-1)\eta+\alpha_{s}}< \pi_{d+l}\leq\break $$ \overline{(b-1)\eta+\alpha_\lambda}.$
 For such  a part $\pi_{d+l}$, we have
 $[\pi_{d+l}/\eta]=b-1$.
 Assume that there are $f_2$ parts $\pi_{d+l}$  in \eqref{lemma2.6-pfa} satisfying
$b\eta\leq \pi_{d+l}\leq \overline{b\eta+\alpha_s}.$  For  such  a part $\pi_{d+l}$,  we have $[\pi_{d+l}/\eta]=b$.

  Assume that there are $f_3$ parts $\pi_{c+l}$  in \eqref{lemma2.6-pfb} satisfying
$\overline{b\eta+\alpha_{s}}< \pi_{c+l}< (b+1)\eta.$ In this case, we have $[\pi_{c+l}/\eta]=b$.
  Assume that there are $f_4$ parts $\pi_{c+l}$   in \eqref{lemma2.6-pfb} satisfying
  $ \pi_{c+l}= (b+1)\eta,$  which gives  $[\pi_{c+l}/\eta]=b+1$. Assume that there are $f_5$ parts
  $\pi_{c+l}$ in \eqref{lemma2.6-pfb} satisfying
  $(b+1)\eta< \pi_{c+l}\leq \overline{(b+1)\eta+\alpha_s},$ which implies  $[\pi_{c+l}/\eta]=b+1$.
To sum up, we get
  \begin{eqnarray*}
   &&[\pi_{c+t-1}/\eta]+\cdots
   +[\pi_{c}/\eta]-
   \left([\pi_{d+k-2}/\eta]+\cdots+[\pi_{d+k-1-t}/\eta]\right)\nonumber\\[5pt]
   &&\quad = bf_3+(b+1)f_4+(b+1)f_5-(b-1)f_1-bf_2,
   \end{eqnarray*}
  and
    \begin{equation}\label{cong-g-x-tem-over}
    \overline{V}_\pi(\pi_c)-\overline{V}_\pi(\pi_d)=f_3+f_5.
    \end{equation}
We proceed to show that $f_1=f_3+f_4$ and $f_2=f_5$.  By means of  \eqref{sequence-lempa3}, we obtain   $f_2+k-t-1+f_3+f_4\leq k-1$, that is, $f_2+f_3+f_4\leq t$. Since $f_1+f_2=t$,  we have
\begin{equation}\label{cong-g-x-tema-over}
f_1\geq f_3+f_4.
\end{equation}
To prove
\begin{equation}\label{cong-g-x-temb-over}
f_1\leq f_3+f_4,
\end{equation}
we consider  three cases:

(1) If $t=k-1$, then  $f_3+f_4+f_5=t=k-1$. Using the condition that $k-1>\lambda$, we have
\begin{equation*}\label{lem4.2bb}
f_1+f_5\leq (\lambda-s)+(s+1)=\lambda+1\leq k-1.
 \end{equation*}
This yields \eqref{cong-g-x-temb-over}.

(2) If $1\leq t<k-1$ {and $\pi_d<(b+1)\eta$}, then  $(b+1)\eta>\pi_d=\pi_{c+t}\geq \pi_{c+k-2}>\overline{b\eta+\alpha_s}$, and so we may write $\pi_d=\overline{b\eta+\alpha_g}$ with $g>s$. It follows that $\pi_{d+k-2}>\overline{(b-1)\eta+\alpha_g}$. Since
\[\overline{b\eta+\alpha_s}<\pi_{d+k-2-t}\leq   \cdots \leq \pi_{d}=\overline{b\eta+\alpha_g},\]
we find that
\begin{equation}\label{lem4.2aaaa}
   k-t-1\leq g-s.
\end{equation} Given the condition that $k-1>\lambda$, we obtain
 \begin{equation}\label{lem4.2aa}
 f_1+f_5\leq (\lambda-g)+(s+1)=(\lambda+1)-g+s\leq k-1-g+s.
  \end{equation}
Combining  \eqref{lem4.2aaaa} and \eqref{lem4.2aa} gives $f_1+f_5\leq t$. Since $f_3+f_4+f_5=t$,  we arrive at \eqref{cong-g-x-temb-over}.

{ (3) If  $1\leq t<k-1$  and $\pi_d\geq(b+1)\eta$, then  $\pi_{d+k-2}\geq\pi_d-\eta\geq b\eta$, and so
\[{b\eta}\leq\pi_{d+k-2}\leq   \cdots \leq \pi_{d+k-1-t}\leq \overline{b\eta+\alpha_s}.
\]
This implies that $f_1=0$, which leads to
\eqref{cong-g-x-temb-over}.
}

Returning to the special case  there are no   overlined parts divisible by $\eta$ in $\pi$, we have $f_5\leq s$, and so
 \eqref{cong-g-x-temb-over} is also valid for $\lambda=k-1$.
 To sum up,  \eqref{cong-g-x-temb-over}  is justified for all cases.  Combining with \eqref{cong-g-x-tema-over}, we conclude that $f_1=f_3+f_4$.

It is now clear that $f_2=f_5$ since $f_1+f_2=f_3+f_4+f_5=t$ and  $f_1=f_3+f_4$. Thus,
 \begin{eqnarray}\label{lemma2.5aa}
   &&[\pi_{c+t-1}/\eta]+\cdots
   +[\pi_{c}/\eta]-
   \left([\pi_{d+k-2}/\eta]+\cdots+[\pi_{d+k-1-t}/\eta]\right)\nonumber\\[5pt]
   &&\qquad= bf_3+(b+1)f_4+(b+1)f_5-(b-1)f_1-bf_2 \nonumber\\[5pt]
   &&\qquad= bf_3+(b+1)f_4+(b+1)f_5-(b-1)(f_3+f_4)-bf_5 \nonumber \\[5pt]
   &&\qquad=f_3+2f_4+f_5 \nonumber \\[5pt]
   &&\qquad \equiv f_3+f_5 \pmod{2}.
   \end{eqnarray}
Substituting  \eqref{cong-g-x-tem-over} into \eqref{lemma2.5aa}, we reach
\eqref{cong-g-x-f-over}, and this completes the proof.      \qed

We conclude this section with a proof of Proposition   \ref{parity k-1 sequence-over} resorting to
 Lemma \ref{parity k-1 sequence-over-old}.

 \noindent{\it Proof of Proposition  \ref{parity k-1 sequence-over}.}  Given a $(k-1)$-marked part $\tilde{g}_p(\pi)$ in $G(\pi)$, we like to show that a $(k-1)$-band of $\pi$
 including $\tilde{g}_p(\pi)$ has the same parity as that of  the $(k-1)$-band of $\pi$ induced by $\tilde{g}_p(\pi)$. Assume that $\tilde{g}_p(\pi)$ is the $g_p$-th part  of $\pi=(\pi_1,\pi_2,\ldots, \pi_\ell)$ in ${{\overline{\mathcal{B}}}}_1(\alpha_1,\ldots,\alpha_\lambda;\eta,k,r)$, that is, $\pi_{g_p}=\tilde{g}_p(\pi)$,
 then $\{\pi_{g_p+l}\}_{0\leq l\leq k-2}$ is the $(k-1)$-band  induced by $\tilde{g}_p(\pi)=\pi_{g_p}$. Assume that $\{\pi_{c+l}\}_{0\leq l\leq k-2}$ is  a $(k-1)$-band of $\pi$ including $\tilde{g}_p(\pi)=\pi_{g_p}$ and   $g_p=c+t$, where $1\leq t\leq k-2$.
Since $\pi$ is an overpartition in ${{\overline{\mathcal{B}}}}_1(\alpha_1,\ldots,\alpha_\lambda;\eta,k,r)$,  we have $\pi_{g_p-1}\geq \pi_{g_p+k-2}+\eta$ with strict inequality if $\pi_{g_p+k-2}$ is non-overlined. By the definition of $(k-1)$-bands, we have
$\pi_{g_p}\leq\pi_{g_p+k-2}+\eta$   with strict inequality if $ \pi_{g_p+k-2} $   is overlined. Thus,  $\pi_{g_p-1}>\pi_{g_p}$, and so $\pi_c\geq \pi_{c+t-1}=\pi_{g_p-1}>\pi_{g_p}$.
The assumption that $\{\pi_{g_p+l}\}_{0\leq l\leq k-2}$ and $\{\pi_{c+l}\}_{0\leq l\leq k-2}$ are $(k-1)$-bands indicates
\[\pi_{g_p+k-2}\geq \pi_{g_p}-\eta=\pi_{c+t}-\eta\geq \pi_{c+k-2}-\eta\geq \pi_c-2\eta, \]
with strict inequality if $\pi_c$ is overlined. Thus, the conditions of Lemma \ref{parity k-1 sequence-over-old} are satisfied, thereby $\{\tilde{g}_p(\pi)\}_{k-1}$ and
$\{\pi_{c+l}\}_{0\leq l\leq k-2}$ are of the same parity. This completes the proof.  \qed

 \section{Proof of Theorem  \ref{rel-over1}}
 The main objective of this section is to give a combinatorial proof of Theorem \ref{rel-over1}. The relationship between  $\overline{B}_1$ and $B_0$ stated in Theorem \ref{rel-over1} plays a crucial role in the proof of Bressoud's conjecture for the case $j=0$ in a subsequent paper \cite{He-Ji-Zhao}.

  Let $\mathcal{D}_\eta$ denote the set of  partitions with distinct parts divisible by $\eta$.      Theorem \ref{rel-over1} is equivalent to   the following combinatorial statement.

\begin{thm}\label{lem-b-0-over} Let $\lambda,k$ and $r$ be   integers  such that $k\geq r\geq \lambda\geq0$ and $k>\lambda$.
There is a bijection $\Phi$ between $\mathcal{D}_\eta\times{\mathcal{B}}_{0}(\alpha_1,\ldots,\alpha_\lambda;\eta,k,r)$ and $\mathcal{\overline{B}}_1(\alpha_1,\ldots,\alpha_\lambda;\eta,k,r)$, namely, for a pair $(\zeta,\mu)\in \mathcal{D}_\eta\times\mathcal{{B}}_{0}
(\alpha_1,\ldots,\alpha_\lambda;\eta,k,r)$, we have $\pi=\Phi(\zeta,\mu)\in\mathcal{\overline{B}}_1(\alpha_1,\ldots,\alpha_\lambda;\eta,k,r)$ such that $|\pi|=|\zeta|+|\mu|$.
\end{thm}

   The bijection $\Phi$ is constructed via merging  $\zeta$ and $\mu$ to produce an overpartition $\pi$ in $\mathcal{\overline{B}}_1(\alpha_1,\ldots,\alpha_\lambda;\eta,k,r)$.  Recall that  $ {\mathcal{B}}_0(\alpha_1,\ldots,\alpha_\lambda;\eta,k,r)$ is the set of overpartitions  in $\overline{\mathcal{B}}_0(\alpha_1,\ldots,\alpha_\lambda;\eta,k,r)$  without overlined parts divisible by $\eta$.
      Assume there are $N$ parts marked with $k-1$ in $RG({\mu})$  and let $\zeta=(\eta\zeta_1,\ldots,\eta\zeta_{c},
\eta\zeta_{c+1}\ldots,\eta\zeta_{c+m})$ be a partition in $\mathcal{D}_\eta$ with $\zeta_1>\cdots>\zeta_c> N\geq\zeta_{c+1}>\cdots>\zeta_{c+m}>0$.  In fact,  the bijection $\Phi$ consists of two steps. The first step  is to
merge the parts
$\eta\zeta_{c+1},\eta\zeta_{c+2},\ldots,\eta\zeta_{c+m}$ and ${\mu}$.   The second step is   to  merge the remaining parts $\eta\zeta_1, \eta\zeta_2,\ldots,\eta\zeta_c$ of $\zeta$ and ${\mu}$ to generate certain ovelined parts   divisible by $\eta$. As will be seen,   the overpartition $\nu$ obtained in the first step is   in   ${{\mathcal{B}}}_1(\alpha_1,\ldots,\alpha_\lambda;\eta,k,r)$. In the meantime,  there are no  overlined parts divisible by $\eta$ in $\nu$. Eventually,  the resulting overpartition $\pi$  of the second step is in $\mathcal{\overline{B}}_1(\alpha_1,\ldots,\alpha_\lambda;\eta,k,r)$.

To describe the map $\Phi$, we introduce   the forward move and  the backward move which are defined on the Gordon marking of a $\overline{B}_1$-overpartition and the reverse Gordon marking of a $\overline{B}_1$-overpartition.  A precise description of the first merging operation will be given later based on the restricted forward move and the restricted backward move and an explanation of the second merging operation will be provided by means of the $(k-1)$-insertion operation  and  the $(k-1)$-separation operation.

 \subsection{The forward move and the backward move}

\begin{defi}[The forward move]\label{forward-defi}For  $k>\lambda$ and $N\geq 1$, let $\pi$ be an overpartition satisfying {\rm(1)}, {\rm(2)} and  {\rm(3)}  in Definition \ref{defi-O-B}. Assume that there are $N$ parts marked with $k-1$ in  $RG(\pi)$, say  $\tilde{r}_1(\pi)>\tilde{r}_2(\pi)>\cdots>
\tilde{r}_N(\pi)$.  For $1\leq p\leq N$,  the forward  move  $\phi_p$ is defined as follows{\rm :} add $\eta$ to each of  $\tilde{r}_1(\pi),\,\tilde{r}_2(\pi),\ldots,\tilde{r}_{p}(\pi)$ and  rearrange the parts in non-increasing order to obtain a new overpartition, denoted $\phi_p(\pi)$.
\end{defi}

For example, let $\pi$ be  the overpartition given in  \eqref{mark-exa-1}.   Below is  the reverse Gordon marking   of $\pi$:
\begin{equation}  \label{rmark-exa-1-s-0000000001}
\begin{split}
&RG(\pi)=(\overline{80}_1,\overbrace{{ {80}_2,{80}_3,\overline{70}_1},
{ {70}_4}}^{{\{70\}_4}},\overline{69}_2,
\overbrace{{ \overline{60}_1,{60}_3,\overline{55}_2}, { \overline{51}_4}}^{{\{\overline{51}\}_4}},\overbrace{{{50}_1,\overline{49}_3,\overline{45}_2}, { \overline{41}_4}}^{{\{\overline{41}\}_4}},\overline{39}_1,\overline{35}_2,\\[5pt]
&\ \ \ \ \ \ \ \ \ \ \ \ \ \  \underbrace{{ \overline{29}_1,\overline{20}_2,} { {20}_3},{ {20}_4}}
_{\{20\}_4},
\underbrace{{ \overline{11}_1,\overline{10}_2,\overline{9}_3},
{ \overline{5}_4}}_{\{\overline{5}\}_4},\overline{1}_1).
\end{split}
\end{equation}
There are five $4$-marked parts in $RG(\pi)$. Applying the forward move $\phi_3$ to  $\pi$   in \eqref{mark-exa-1},    we obtain
 \begin{equation*}\label{exa-r-2}
\begin{split}
&\phi_3(\pi)=(\overline{80},{ {80}},{80},{80},\overline{70},\overline{69},{\overline{61}},\overline{60},
{60},\overline{55},{\overline{51}},{50},\overline{49},\overline{45},\overline{39},\overline{35},\\[3pt]
&\ \ \ \ \ \ \ \ \ \ \ \ \  \overline{29},\overline{20},{20},{20},\overline{11},
\overline{10},\overline{9},\overline{5},\overline{1}).
\end{split}
\end{equation*}

The following proposition gives several properties of  $\phi_p(\pi)$.

\begin{prop}\label{lem-for1}For  $k>\lambda$ and $N\geq 1$, let $\pi=(\pi_1,\ldots,\pi_\ell)$ be an overpartition satisfying {\rm(1)}, {\rm(2)} and  {\rm(3)}  in Definition \ref{defi-O-B}. Assume that there are $N$ parts marked with $k-1$ in $RG(\pi)$,  say    $\tilde{r}_1(\pi)>\tilde{r}_2(\pi)>\cdots>
\tilde{r}_N(\pi)$.   For $1\leq p\leq N$, let $\omega=(\omega_1,\ldots,\omega_\ell)=
\phi_p(\pi)$. Then
\begin{itemize}
  \item[{\rm (1)}] For $1\leq i\leq \ell$, $\omega_i\equiv0,\alpha_1,\ldots,\alpha_\lambda\pmod{\eta}${\rm{;}}

 \item[{\rm (2)}] Only  multiples of $\eta$ may be non-overlined in $\omega${\rm{;}}

    \item[{\rm (3)}] There are at most $k-1$ marks in $G(\omega)$ and  there are $N$ parts marked with $k-1$ in  $G(\omega)$,     say $\tilde{g}_1(\omega)>\tilde{g}_2(\omega)>\cdots>
\tilde{g}_N(\omega)${\rm{;}}

   \item[{\rm (4)}] Let $\{\tilde{r}_i(\pi)\}_{k-1}$ be  the $(k-1)$-band of $\pi$ induced by $\tilde{r}_i(\pi)$. Then  $ \tilde{g}_i(\omega)= \tilde{r}_i(\pi)+\eta$ for $1\leq i\leq p,$
   and $\tilde{g}_i(\omega) \in \{\tilde{r}_i(\pi)\}_{k-1}$ for $p< i\leq N$.
\end{itemize}
\end{prop}

For example, let $\pi$ be the overpartition in $\overline{\mathcal{B}}_1(1,5,9;10,5,4)$ given in    \eqref{mark-exa-1},  and let $\omega=\phi_3(\pi)$. Then
the Gordon marking of $\omega$ is given by
\begin{equation}\label{example-b}
\begin{split}
&G(\omega)=(\overline{80}_1,\overbrace{{ {80}_4},
{ {80}_3,
{80}_2,\overline{70}_1}}^{\{80\}_4},\overline{69}_2,\overbrace{{ \overline{61}_4},{ \overline{60}_3,{60}_1,\overline{55}_2}}^{\{\overline{61}\}_4},\overbrace{{ \overline{51}_4},{ {50}_3,\overline{49}_1,\overline{45}_2}}^{\{\overline{51}\}_4},\overline{39}_1,\overline{35}_2,\\[5pt]
&\ \ \ \ \ \ \ \ \ \ \ \
\overline{29}_1,\underbrace{{ \overline{20}_4},{ {20}_3,
{20}_2,\overline{11}_1}}_{\{\overline{20}\}_4},\underbrace{{ \overline{10}_4},
{ \overline{9}_3,\overline{5}_2,\overline{1}_1}}_{\{\overline{10}\}_4}).
\end{split}
\end{equation}
It can be readily checked that  $\omega$ satisfies  the properties (1)-(4) in   Proposition \ref{lem-for1}.

\noindent{\it Proof of Proposition \ref{lem-for1}.}  To prove   (1) and (2), it suffices to show that for $1\leq i\leq p$,  $\tilde{r}_i(\pi)+\eta$ cannot be repeated in $\omega$ if the part $\tilde{r}_i(\pi)$ is overlined. We now assume that  $\tilde{r}_i(\pi)$ is overlined, so the new generated part  $\tilde{r}_i(\pi)+\eta$ is also overlined. There are two cases.

Case 1: Assume that $\tilde{r}_i(\pi)+\eta$ is not a part of $\pi$. It is obvious that the generated overlined part $\tilde{r}_i(\pi)+\eta$ appears only once in  $\omega$.

Case 2: Assume that $\pi$ contains  an overlined part ${\pi}_t=\tilde{r}_i(\pi)+\eta$.  We claim that ${\pi}_t$ is marked with $k-1$ in $RG(\pi)$.  Since $\pi$ satisfies  the condition  {\rm(3)} in Definition \ref{defi-O-B},  the marks in  $RG(\pi)$ do not  exceed $k-1$.   Assume that $\tilde{r}_i(\pi)$ is the $r_i$-th part  of $\pi=(\pi_1,\pi_2,\ldots, \pi_\ell)$, that is, $\pi_{r_i}=\tilde{r}_i(\pi)$. Then  $\{\pi_{r_i-l}\}_{0\leq l\leq k-2}$ is   the $(k-1)$-band induced by $\tilde{r}_i(\pi)$. It follows that  the marks of the parts $\pi_{r_i-k+2},\ldots,\pi_{r_i}$  are distinct in $RG(\pi)$ and   $\pi_{r_i}$ is marked with $k-1$. This implies that the marks of $\pi_{r_i-k+2},\ldots,\pi_{r_i-1}$ in $RG(\pi)$ are distinct and less than $k-1$.  Suppose  to the contrary that the mark of $\pi_{t}$ in $RG(\pi)$ is less than $k-1$. Consequently, there is a part $\pi_{r_i-m}$ ($1\leq m\leq k-2$) in the $(k-1)$-band $\{\pi_{r_i-l}\}_{0\leq l\leq k-2}$  such that the mark of $\pi_{r_i-m}$ is the same as the mark of $\pi_{t}$.    Since  $\pi_{r_i}<\pi_{r_i-m}<\pi_{r_i}+\eta$ and $\pi_{t}=\pi_{r_i}+\eta$, we obtain   $\pi_{r_i-m}<\pi_{t}<\pi_{r_i-m}+\eta$, and so the marks of   $\pi_{t}$ and $\pi_{r_i-m}$ are distinct. But this is impossible under the prior assumption.  Therefore, the mark of $\pi_t$ is  $k-1$ in $RG(\pi)$, as claimed. In other words, we have $\tilde{r}_{i-1}(\pi)=\pi_t$. This enables us to employ the forward move to add $\eta$ to  $\tilde{r}_{i-1}(\pi)$. In the end,  the  generated overlined part $\tilde{r}_i(\pi)+\eta$ occurs only once in $\omega$.

We now turn to the properties (3) and (4). Assume that $\tilde{r}_i(\pi)$ is the $r_i$-th part  of $\pi=(\pi_1,\pi_2,\ldots, \pi_\ell)$, that is, $\pi_{r_i}=\tilde{r}_i(\pi)$. In fact, the forward move consists of two steps. First,  remove  the $(k-1)$-marked parts $\pi_{r_1},\,\pi_{r_2},\ldots, \pi_{r_p}$ from $\pi$ and denote the resulting overpartition  by $\pi^{(1)}$. Since   the largest mark in  $RG(\pi)$ is  $k-1$ and the parts removed from $\pi$ are marked with $k-1$ in $RG(\pi)$,  the marks of the remaining parts in $RG(\pi^{(1)})$ are the same as those in $RG(\pi)$. This implies that   the marks of the parts  of $\pi^{(1)}$ not less than $\pi_{r_p}$ do not exceed  $k-2$. Therefore, there are $N-p$ parts marked with $k-1$ in  $RG(\pi^{(1)})$, denoted $\tilde{r}_{1}(\pi^{(1)}),\ldots, \tilde{r}_{N-p}(\pi^{(1)})$, for which   $\tilde{r}_{i}(\pi^{(1)})=\tilde{r}_{i+p}(\pi) $ and $\{\tilde{r}_{i}(\pi^{(1)})\}_{k-1}=\{\tilde{r}_{i+p}(\pi)\}_{k-1}$ for $1\leq i\leq N-p$. In light of Proposition \ref{sequence-length}, we find that  there are also $N-p$ parts marked with $k-1$  in  $G(\pi^{(1)})$, denoted $\tilde{g}_{1}(\pi^{(1)}),\ldots, \tilde{g}_{N-p}(\pi^{(1)})$, for which  $\tilde{g}_{i}(\pi^{(1)}) \in \{\tilde{r}_{i}(\pi^{(1)})\}_{k-1}$. Meanwhile,  the marks of the parts not less than $\pi_{r_p}$   in $G(\pi^{(1)})$ do not exceed  $k-2$.  So we deduce that  $\tilde{g}_{i}(\pi^{(1)}) \in  \{\tilde{r}_{i+p}(\pi)\}_{k-1}$  for $1\leq i\leq N-p$.

The second step is to insert   $\pi_{r_1}+\eta,\,\pi_{r_2}+\eta,\ldots, \pi_{r_p}+\eta$   into $\pi^{(1)}$ and to rearrange the parts in non-increasing order to obtain $\omega$. We wish to show that for $1\leq i\leq p$, $\pi_{r_i}+\eta$  is marked with $k-1$ in $G(\omega)$.
We claim that $\pi_{r_i}>\pi_{r_i+1}$.  For $1\leq i\leq p$, since $\pi_{r_i}$ is the $(k-1)$-marked part in $RG(\pi)$, we know that $\{\pi_{r_i-l}\}_{0\leq l\leq k-2}$ is  the $(k-1)$-band of $\pi$ induced by $\pi_{r_i}$, which ensures that $\pi_{r_i-k+2}\leq \pi_{r_i}+\eta$ with strict inequality if $ \pi_{r_i-k+2}$  is overlined.
Under the assumption that $\pi$ satisfies  the condition  {\rm(3)} in Definition \ref{defi-O-B}, we have  $\pi_{r_i-k+2}\geq \pi_{r_i+1}+\eta$ with strict inequality if $ \pi_{r_i-k+2}$  is non-overlined.   But $\pi_{r_i-k+2}\leq \pi_{r_i}+\eta$  with strict inequality if $ \pi_{r_i-k+2}$  is overlined, so we conclude that  $\pi_{r_i}>\pi_{r_i+1}$, as claimed.

We continue to prove that $\pi_{r_i}+\eta$ is marked with $k-1$ in $G(\omega)$.  Based on the assumption  that $\pi$ satisfies the condition {\rm(3)} in Definition \ref{defi-O-B}, we obtain that
$\pi_{r_i-k+1}\geq \pi_{r_i}+\eta$ with strict inequality if $ \pi_{r_i-k+1}$  is non-overlined,  which implies that  for $1\leq i\leq p$,
\[ \pi_{r_i}+\eta, \text{ and } \pi_{r_i-k+2},   \ldots, \pi_{r_i-1}\]
are  parts of  $\omega$. Noting that $\pi_{r_i}>\pi_{r_i+1}$ and $\pi_{r_i-k+1}\geq \pi_{r_i}+\eta$ with strict inequality if $\pi_{r_i}$ is non-overlined, it is clear that   the mark of the part $\pi_{r_i}+\eta$ in $G(\omega)$ is  the smallest positive integer that is not used to mark $\pi_{r_i-k+2},   \ldots, \pi_{r_i-1}$.   Recalling that the marks of the parts of $\pi^{(1)}$ not less than $\pi_{r_p}$ in $G(\pi^{(1)})$ do not exceed  $k-2$,    the marks of the parts $\pi_{r_i-k+2},\,\pi_{r_i-k+3},\, \ldots ,\, \pi_{r_i-1}$  in $G(\omega)$ are less than $k-1$ for  $1\leq i\leq p$.  Thus, the mark of $\pi_{r_i}+\eta$ in $G(\omega)$
is $k-1$. Meanwhile, the marks of  remaining parts  in $G(\omega)$ are   the same as in $G(\pi^{(1)})$. Therefore,  we reach the conclusion that there are $N$ parts marked with $k-1$ in $G(\omega)$, and so the properties (3) and (4) are verified. This completes the proof.   \qed

For example, for  the  overpartition $\pi$ in $\overline{\mathcal{B}}_1(1,5,9;10,5,4)$ with the reverse Gordon marking given in  \eqref{rmark-exa-1-s-0000000001},   there are five $4$-marked parts in $RG(\pi)$.   Then $\omega=\phi_3(\pi)$ can be constructed via  two steps: The first step is to remove  $\tilde{r}_1(\pi)=70$, $\tilde{r}_2(\pi)=\overline{51}$, and  $\tilde{r}_3(\pi)=\overline{41}$   from $\pi$ to get $\pi^{(1)}$, whose reverse Gordon marking reads
 \[\begin{split}
&RG(\pi^{(1)})=(\overline{80}_1,{{80}_2,{80}_3,\overline{70}_1},\overline{69}_2,
{\overline{60}_1,{60}_3,\overline{55}_2},{{50}_1,\overline{49}_3,\overline{45}_2},\overline{39}_1,\overline{35}_2,\\[5pt]
&\ \ \ \ \ \ \ \ \ \ \ \ \ \ \ \ \  \underbrace{{ \overline{29}_1,\overline{20}_2,} { {20}_3},{ {20}_4}}
_{\{20\}_4},
\underbrace{{ \overline{11}_1,\overline{10}_2,\overline{9}_3},
{ \overline{5}_4}}_{\{\overline{5}\}_4},\overline{1}_1).
\end{split}\]
It can be checked that the marks of   parts in $RG(\pi^{(1)})$ are the same as those in  $RG(\pi)$.  On the other hand, below is the Gordon marking of $\pi^{(1)}$
 \[\begin{split}
&G(\pi^{(1)})=(\overline{80}_1,{{80}_3},
{{80}_2,
\overline{70}_1},\overline{69}_2,{\overline{60}_3},{{60}_1,\overline{55}_2},{{50}_3},{\overline{49}_1,\overline{45}_2},\overline{39}_1,\overline{35}_2,\\[5pt]
&\ \ \ \ \ \ \ \ \ \ \ \ \ \ \
\overline{29}_1,{ \overline{20}_4},{{20}_3,
{20}_2,\overline{11}_1},{ {\overline{10}_4}},
{\overline{9}_3,\overline{5}_2,\overline{1}_1}).
\end{split}\]
Evidently, the $4$-marked parts $\overline{20}_4$ and $\overline{10}_4$ in $G(\pi^{(1)})$  are in the $4$-bands $\{{ \overline{29},\overline{20},} { {20}},{ {20}}\}$ and $\{ \overline{11},\overline{10},\overline{9},
\overline{5}\}$ of $\pi^{(1)}$, respectively. In the second step, we insert $\tilde{r}_1(\pi)+10=80$, $\tilde{r}_2(\pi)+10=\overline{61}$,
and  $\tilde{r}_3(\pi)+10=\overline{51}$ into $\pi^{(1)}$ to get $\omega$, whose  Gordon marking   is displayed in \eqref{example-b}. As anticipated,  the marks of $\tilde{r}_1(\pi)+10=80$, $\tilde{r}_2(\pi)+10=\overline{61}$,
and  $\tilde{r}_3(\pi)+10=\overline{51}$  in $G(\omega)$
are $4$. Meanwhile the marks of  remaining parts  in $G(\omega)$ are   the same as in $G(\pi^{(1)})$. Therefore, $\omega$ satisfies  the properties (1)-(4) in   Proposition \ref{lem-for1}.

In parallel to the forward move, we now turn to the definition of the backward move relying on the Gordon marking of a $\overline{B}_1$-overpartition.

\begin{defi}[The backward move]\label{defi-backward}
For $k>\lambda$ and $N\geq p\geq 1$, let $\omega$ be an overpartition satisfying {\rm(1)}, {\rm(2)} and  {\rm(3)}  in Definition \ref{defi-O-B}.  Assume that there are $N$ parts marked with $k-1$ in  $G(\omega)$, denoted   $\tilde{g}_1(\omega)>\tilde{g}_2(\omega)>\cdots>
\tilde{g}_N(\omega)$, for which $\tilde{g}_p(\omega)\geq\overline{\eta+\alpha_1}$. The backward  move $\psi_p$ is defined as follows{\rm :} subtract $\eta$ from each  of $\tilde{g}_1(\omega),\,\tilde{g}_2(\omega),\ldots,
 \tilde{g}_{p}(\omega)$ and  rearrange the parts in  non-increasing order to obtain a new overpartition, denoted  $\psi_p(\omega)$.
\end{defi}

For example,  for the overpartition $\omega=\phi_3(\pi)$  with five $4$-marked parts in $G(\omega)$ as in  \eqref{example-b}, the backward move  $\psi_3$ transforms  $\omega$ back to $\pi$   in  \eqref{rmark-exa-1-s-0000000001}.

The backward move $\psi_p$ possesses the following properties with respect to certain overpartitions satisfying {\rm(1)}, {\rm(2)} and  {\rm(3)}  in Definition \ref{defi-O-B}.

\begin{prop}\label{lem-bac1} For $k>\lambda$ and $N\geq p\geq 1$, let $\omega$ be an overpartition satisfying {\rm(1)}, {\rm(2)} and  {\rm(3)}  in Definition \ref{defi-O-B}.  Assume that there are $N$ parts marked with $k-1$ in  $G(\omega)$, denoted  $\tilde{g}_1(\omega)>\tilde{g}_2(\omega)>\cdots>
\tilde{g}_N(\omega)$, for which $\tilde{g}_p(\omega)\geq \overline{\eta+\alpha_1}$, and assume that $\tilde{g}_p(\omega)$ is a part in any $(k-1)$-band  of $\omega$  belonging to  $I(\tilde{g}_p(\omega)-2\eta, \tilde{g}_p(\omega))$.
Let $\pi=(\pi_1,\ldots, \pi_\ell)=\psi_p(\omega)$. Then
 \begin{itemize}
     \item[{\rm (1)}] For $1\leq i\leq \ell$, $\pi_i\equiv0,\alpha_1,\ldots,\alpha_\lambda\pmod{\eta}${\rm{;}}

 \item[{\rm (2)}] Only  multiples of $\eta$ may be non-overlined in $\pi${\rm{;}}

     \item[{\rm{(3)}}] There are at most $k-1$ marks in  $RG(\pi)$ and there are $N$ parts marked with $k-1$ in $RG(\pi)$, say   $\tilde{r}_1(\pi)>\cdots
>\tilde{r}_N(\pi)${\rm{;}}

  \item[{\rm{(4)}}] Let   $\{\tilde{g}_i(\omega)\}_{k-1}$  be the $(k-1)$-band of $\omega$ induced by $\tilde{g}_i(\omega)$. Then $\tilde{r}_i(\pi)=\tilde{g}_i(\omega)
  -\eta$  for $1\leq i\leq p$,  and $ \tilde{r}_i(\pi)\in\{\tilde{g}_{i}(\omega)\}_{k-1}$  for $p< i\leq N$.
  \end{itemize}
\end{prop}
\pf  To prove   (1) and (2), it suffices to show that for  $1\leq i\leq p$, the   part   $\tilde{g}_i(\omega)-\eta$  occurs exactly once in $\pi$  if $\tilde{g}_i(\omega)$ is
overlined. We now assume that $\tilde{g}_i(\omega)$ is
overlined, so that the generated part  $\tilde{g}_i(\omega)-\eta$ is also overlined. There are two cases:

Case 1: Assume that the  part  $\tilde{g}_i(\omega)-\eta$ does not appear in  $\omega$. In this case, it is obvious that the generated part $\tilde{g}_i(\omega)-\eta$ occurs exactly once in  $\pi$.

Case 2: Assume that $\omega$ contains the overlined part $\omega_{t_i}=\tilde{g}_i(\omega)-\eta$, where $1\leq i\leq p$.   Using the same argument as in the proof of Proposition  \ref{lem-for1}, it can be shown that   $\omega_{t_i}$ is  marked with $k-1$ in $G(\omega)$, where $1\leq i\leq p$.

We proceed to show that  if $\omega$ satisfies the condition that $\tilde{g}_p(\omega)$ is a part in any $(k-1)$-band  of $\omega$  belonging to $I(\tilde{g}_p(\omega)-2\eta, \tilde{g}_p(\omega))$, then    $\omega$ does not contain the overlined part $\tilde{g}_p(\omega)-\eta$.   Suppose to the contrary that  $\omega$   contains the overlined part $\tilde{g}_p(\omega)-\eta$. Since  $\tilde{g}_p(\omega)-\eta$ is marked with $k-1$ in $G(\omega)$,  we have
\begin{equation}\label{back-pf0}
  \tilde{g}_{p+1}(\omega)=\tilde{g}_p(\omega)-\eta,
\end{equation}
where $\tilde{g}_{p+1}(\omega)$ is also overlined. Assume that $\tilde{g}_{p+1}(\omega)$ is the $g_{p+1}$-th part  of $\omega=(\omega_1,\omega_2,\ldots,\break $$ \omega_\ell)$, that is, $\omega_{g_{p+1}}=\tilde{g}_{p+1}(\omega)$. Then $\{\omega_{g_{p+1}+l}\}_{0\leq l\leq k-2}$ is the $(k-1)$-band induced by $\omega_{g_{p+1}}$, which, together with \eqref{back-pf0}, leads to
\[\omega_{g_{p+1}+k-2}>\omega_{g_{p+1}}-\eta=\tilde{g}_p(\omega)-2\eta.\]
By \eqref{back-pf0}, we have
\begin{equation*}\label{back-pf1}
 \omega_{g_{p+1}}=\tilde{g}_p(\omega)-\eta<\tilde{g}_p(\omega),
\end{equation*}
from which  $\{\omega_{g_{p+1}+l}\}_{0\leq l\leq k-2}$ is a $(k-1)$-band   of $\omega$  belonging to $I(\tilde{g}_p(\omega)-2\eta, \tilde{g}_p(\omega))$. But $\tilde{g}_p(\omega)$ is not a part in $\{\omega_{g_{p+1}+l}\}_{0\leq l\leq k-2}$, which contradicts the condition  that   $\tilde{g}_p(\omega)$ is a part in any $(k-1)$-band  of $\omega$  belonging to $I(\tilde{g}_p(\omega)-2\eta, \tilde{g}_p(\omega))$. This means that $\omega$ does not contain the overlined part $\tilde{g}_p(\omega)-\eta$. Using the fact that  $\tilde{g}_i(\omega)-\eta$ is marked with $k-1$ in $G(\omega)$, we have $\tilde{g}_{i+1}(\omega)=\tilde{g}_{i}(\omega)-\eta$ where $1\leq i<p$. Applying the backward move to $\omega$, we get the overpartition $\pi$ in which the part $\tilde{g}_i(\omega)-\eta$ appears exactly once. So we have verified   the properties (1) and (2).

We now turn to  the properties (3) and (4). Similarly, the backward move consists of two steps.  First,  remove the $(k-1)$-marked parts $\tilde{g}_{1}(\omega),\ldots, \tilde{g}_{p}(\omega)$ from $\omega$ and denote the resulting overpartition  by $\omega^{(1)}$. Along the same lines of reasoning as in  the proof of Proposition \ref{lem-for1}, we deduce that the marks of    the remaining parts   in $RG(\omega)$ are the same as those in $RG(\omega^{(1)})$. This implies that there are   $N-p$ parts marked with $k-1$  in $RG(\omega^{(1)})$, denoted $\tilde{r}_1(\omega^{(1)})>\cdots
>\tilde{r}_{N-p}(\omega^{(1)})$, for which $\tilde{r}_{i}(\omega^{(1)}) \in \{\tilde{g}_{i+p}(\omega)\}_{k-1}$.   We proceed to demonstrate  that  $\tilde{r}_{1}(\omega^{(1)})\leq\tilde{g}_{p}(\omega)-2\eta $ with strict inequality if $\tilde{g}_{p}(\omega)$ is non-overlined. Suppose to the contrary that $\tilde{r}_{1}(\omega^{(1)})\geq\tilde{g}_{p}(\omega)-2\eta $ with strict inequality if $\tilde{g}_{p}(\omega)$ is overlined. Since $\tilde{r}_{1}(\omega^{(1)}) \in \{\tilde{g}_{p+1}(\omega)\}_{k-1}$, we have $\tilde{r}_{1}(\omega^{(1)})\leq \tilde{g}_{p+1}(\omega)$. Note that $\tilde{g}_{p+1}(\omega)\leq\tilde{g}_{p}(\omega)-\eta$ with strict inequality if $\tilde{g}_{p+1}(\omega)$ is non-overlined. We obtain that $\tilde{r}_{1}(\omega^{(1)})\leq\tilde{g}_{p}(\omega)-\eta$ with strict inequality if $\tilde{r}_{1}(\omega^{(1)})$ is non-overlined. Assume that $\tilde{r}_{1}(\omega^{(1)})$ is the $r_{1}$-th part  of $\omega^{(1)}=(\omega^{(1)}_1,\omega^{(1)}_2,\ldots, \omega^{(1)}_\ell)$, that is, $\omega^{(1)}_{r_1}=\tilde{r}_{1}(\omega^{(1)})$. Then $\{\omega^{(1)}_{r_1-l}\}_{0\leq l\leq k-2}$ is the $(k-1)$-band induced by $\omega^{(1)}_{r_1}$, which implies that $\omega^{(1)}_{r_1-k+2}\leq \omega^{(1)}_{r_1}+\eta$ with strict inequality if $\omega^{(1)}_{r_1}$ is overlined. But, $\tilde{r}_{1}(\omega^{(1)})\leq\tilde{g}_{p}(\omega)-\eta$ with strict inequality if $\tilde{r}_{1}(\omega^{(1)})$ is non-overlined, it follows that $\omega^{(1)}_{r_1-k+2}<\tilde{g}_{p}(\omega)$.  Therefore, we conclude that $\{\omega^{(1)}_{r_1-l}\}_{0\leq l\leq k-2}$  is a $(k-1)$-band of $\omega^{(1)}$ belonging to $I(\tilde{g}_p(\omega)-2\eta, \tilde{g}_p(\omega))$. By the construction of $\omega^{(1)}$, we see that $\{\omega^{(1)}_{r_1-l}\}_{0\leq l\leq k-2}$ is also a $(k-1)$-band of $\omega$ belonging to $I(\tilde{g}_p(\omega)-2\eta, \tilde{g}_p(\omega))$. However, $\tilde{g}_p(\omega)$ is not in  $\{\omega^{(1)}_{r_1-l}\}_{0\leq l\leq k-2}$,  which contradicts the condition that  $\tilde{g}_p(\omega)$ is a part in any $(k-1)$-band  of $\omega$  belonging to $I(\tilde{g}_p(\omega)-2\eta, \tilde{g}_p(\omega))$. Hence $\tilde{r}_{1}(\omega^{(1)})\leq\tilde{g}_{p}(\omega)-2\eta$ with strict inequality if $\tilde{r}_{1}(\omega^{(1)})$ is non-overlined.

The second step is to insert   $\tilde{g}_{1}(\omega)-\eta,\ldots, \tilde{g}_{p}(\omega)-\eta$   into $\omega^{(1)}$ and to rearrange the parts in non-increasing order  to obtain  $\pi$.  It can be shown that the mark of  $ \tilde{g}_{i}(\omega)-\eta$   in $RG(\pi)$ is equal to $k-1$  for $1\leq i\leq p$. Furthermore, the remaining parts not less than $\tilde{g}_{p}(\omega)-\eta$ in $RG(\pi)$ are the same as in $RG(\omega^{(1)})$. We need to show that the marks of  remaining parts less than  $\tilde{g}_{p}(\omega)-\eta$ in $RG(\pi)$ are   the same as in $RG(\omega^{(1)})$. We first  verify that the marks of    the parts $\pi_i$ of $\pi$ such that   $\tilde{g}_{p}(\omega)-2\eta\leq \pi_i<\tilde{g}_{p}(\omega)-\eta$ with strict inequality if $\tilde{g}_{p}(\omega)$ is overlined   in $RG(\pi)$ are   the same as those in $RG(\omega^{(1)})$. Since $\tilde{r}_{1}(\omega^{(1)})\leq\tilde{g}_{p}(\omega)-2\eta $ with strict inequality if $\tilde{g}_{p}(\omega)$ is non-overlined,  the marks of   the parts $\pi_i$ of $\pi$ such that   $\tilde{g}_{p}(\omega)-2\eta\leq \pi_i<\tilde{g}_{p}(\omega)-\eta$ with strict inequality if $\tilde{g}_{p}(\omega)$ is overlined  are less than $k-1$ in $RG(\omega^{(1)})$. But the mark of $\tilde{g}_{p}(\omega)-\eta$ is $k-1$, we infer  that the marks of    the parts $\pi_i$ of $\pi$ such that   $\tilde{g}_{p}(\omega)-2\eta\leq \pi_i<\tilde{g}_{p}(\omega)-\eta$ with strict inequality if $\tilde{g}_{p}(\omega)$ is overlined   in $RG(\pi)$ are   the same as those in $RG(\omega^{(1)})$. Thus, the marks of the parts $\pi_i$ of $\pi$ such that   $\pi_i\leq \tilde{g}_{p}(\omega)-2\eta$ with strict inequality if $\tilde{g}_{p}(\omega)$ is non-overlined in $RG(\pi)$ are   the same as in $RG(\omega^{(1)})$.
Therefore,  the marks of  remaining parts less than  $\tilde{g}_{p}(\omega)-\eta$ in $RG(\pi)$ are   the same as in $RG(\omega^{(1)})$.  It follows that there are $N$ parts marked with $k-1$ in $RG(\pi)$, and so the properties (3) and (4) are verified. This completes the proof.  \qed

We now furnish an example to illustrate   Proposition \ref{lem-bac1}.  Let $\omega$ be the overpartition  in  $\overline{\mathcal{B}}_1(1,5,9;10,5,4)$  with the Gordon marking  given by \eqref{example-b}. There are five $4$-marked parts in $G(\omega)$, namely, $\tilde{g}_1(\omega)=80$, $\tilde{g}_2(\omega)=\overline{61}$,
 $\tilde{g}_3(\omega)=\overline{51}$, $\tilde{g}_4(\omega)=\overline{20}$ and $\tilde{g}_5(\omega)=\overline{10}$.   It can be checked that there are no  $4$-bands of $\omega$ belonging to $I(\tilde{g}_3(\omega)-2\eta, \tilde{g}_3(\omega))$.  The overpartition  $\pi=\psi_3(\omega)$ can be constructed as follows: First, remove  $\tilde{g}_1(\omega)=80$, $\tilde{g}_2(\omega)=\overline{61}$,
and  $\tilde{g}_3(\omega)=\overline{51}$ from $\omega$ to get $\omega^{(1)}$. We have
 \[\begin{split}
&G(\omega^{(1)})=(\overline{80}_1,{{80}_3},
{{80}_2,
\overline{70}_1},\overline{69}_2,{\overline{60}_3},{{60}_1,\overline{55}_2},{{50}_3},{\overline{49}_1,\overline{45}_2},\overline{39}_1,\overline{35}_2,\\[5pt]
&\ \ \ \ \ \ \ \ \ \ \ \ \ \ \
\overline{29}_1,\underbrace{{ \overline{20}_4},{ {20}_3,
{20}_2,\overline{11}_1}}_{\{\overline{20}\}_4},\underbrace{{ {\overline{10}_4}},
{ \overline{9}_3,\overline{5}_2,\overline{1}_1}}_{\{\overline{10}\}_4}).
\end{split}\]
It can be checked that the marks of   parts in $G(\omega^{(1)})$ are the same as those in  $G(\omega)$.  The reverse Gordon marking of $\omega^{(1)}$ is given by
 \[\begin{split}
&RG(\omega^{(1)})=(\overline{80}_1,{{80}_2,{80}_3,\overline{70}_1},\overline{69}_2,
{\overline{60}_1,{60}_3,\overline{55}_2},{{50}_1,\overline{49}_3,\overline{45}_2},\overline{39}_1,\overline{35}_2,\\[5pt]
&\ \ \ \ \ \ \ \ \ \ \ \ \ \ \ \ \  { \overline{29}_1,\overline{20}_2,} { {20}_3},{ {20}_4},
{ \overline{11}_1,\overline{10}_2,\overline{9}_3},
{ \overline{5}_4},\overline{1}_1),
\end{split}\]
from which we see that the $4$-marked parts  $20$ and
 $\overline{5}$ in $RG(\omega^{(1)})$ are also in  $\{\overline{20}\}_4$ and $\{\overline{10}\}_4$ of $\omega^{(1)}$, respectively.
Then $\pi$ can be obtained by inserting $\tilde{g}_1(\omega)-10=70$, $\tilde{g}_2(\omega)-10=\overline{51}$, and  $\tilde{g}_3(\omega)-10=\overline{41}$ into $\omega^{(1)}$. Below is the reverse Gordon marking of $\pi$
\begin{equation*}
\begin{split}
&RG(\pi)=(\overline{80}_1,{{ {80}_2,{80}_3,\overline{70}_1},
{ {70}_4}},\overline{69}_2,
{{ \overline{60}_1,{60}_3,\overline{55}_2}, { \overline{51}_4}},{{{50}_1,\overline{49}_3,\overline{45}_2}, { \overline{41}_4}},\overline{39}_1,\overline{35}_2,\\[5pt]
&\ \ \ \ \ \ \ \ \ \ \ \ \ \ \
{{ \overline{29}_1,\overline{20}_2,} { {20}_3},{ {20}_4}},
{{ \overline{11}_1,\overline{10}_2,\overline{9}_3},
{ \overline{5}_4}},\overline{1}_1).
\end{split}
\end{equation*}
Notice that  $\tilde{g}_1(\omega)-10=70$, $\tilde{g}_2(\omega)-10=\overline{51}$, and  $\tilde{g}_3(\omega)-10=\overline{41}$ are marked with $4$ in $RG(\pi)$. The marks of the remaining parts in  $RG(\pi)$ are the same as those in $RG(\omega^{(1)})$.

We remark that the condition in  Proposition \ref{lem-bac1} that $\tilde{g}_p(\omega)$ is a part in any $(k-1)$-band  of $\omega$  belonging to $I(\tilde{g}_p(\omega)-2\eta, \tilde{g}_p(\omega))$ is necessary. For example,  let ${\omega}$ be  an  overpartition  in ${\mathcal{B}}_1(1,5,9;10,5,4)$ having the Gordon marking
 \begin{equation*}
\begin{split}
&G({\omega})=(\overline{80}_1,{{ {80}_4},
{ {80}_3,
{80}_2,\overline{70}_1}},\overline{69}_2,{{ \overline{61}_4},{ \overline{60}_3,{60}_1,\overline{55}_2}},{{ \overline{51}_4},{ {\overline{50}_3},\overline{49}_1,\overline{45}_2}},{ \overline{40}_4},{40}_3,\overline{39}_1,\overline{35}_2,\\[5pt]
&\ \ \ \ \ \ \ \ \ \ \ \
\overline{29}_1,{{ \overline{20}_4},{ {20}_3,
{20}_2,\overline{11}_1}},{ {\overline{10}_4}},
{ \overline{9}_3,\overline{5}_2,\overline{1}_1}).
\end{split}
\end{equation*}
 There are six $4$-marked parts in $G(\omega)$, namely, $\tilde{g}_1(\omega)=80$, $\tilde{g}_2(\omega)=\overline{61}$,
 $\tilde{g}_3(\omega)=\overline{51}$,
     $\tilde{g}_4(\omega)=\overline{40}$, $\tilde{g}_5(\omega)=\overline{20}$ and
     $\tilde{g}_6(\omega)=\overline{10}$.
Furthermore,  ${\omega}$ has  three $4$-bands $\{\overline{49},\overline{45},\overline{40},{40}\}$, $\{\overline{45},\overline{40},{40},\overline{39}\}$ and  $\{\overline{40},{40},\overline{39},\overline{35}\}$  in  the interval $(\overline{31}, \overline{51})$.

The overpartition $\pi=\psi_3(\omega)$ can be obtained by subtracting $\eta=10$ from each of $\tilde{g}_1({\omega})=80$, $\tilde{g}_2({\omega})=\overline{61}$,
and  $\tilde{g}_3({\omega})=\overline{51}$, and so we get
   \[\begin{split}
&RG(\pi)=(\overline{80}_1,{{80}_2,{80}_3,\overline{70}_1},{ {70}_4},\overline{69}_2,
{\overline{60}_1,{60}_3,\overline{55}_2},{ \overline{51}_4},{{\overline{50}_1},
\overline{49}_3,\overline{45}_2},{ \overline{41}_4},{\overline{40}_1},
{ {{40}_5}},\overline{39}_3,\overline{35}_2,\\[5pt]
&\ \ \ \ \ \ \ \ \ \ \ \ \ \ \ { \overline{29}_1,\overline{20}_2,} { {20}_3},{ {{20}_4}},
{ \overline{11}_1,\overline{10}_2,\overline{9}_3},
{ {\overline{5}_4}},\overline{1}_1).
\end{split}\]
Because of the occurrence of the     $5$-marked part $40$ in $RG(\pi)$, the property (3) in Proposition \ref{lem-bac1} is violated.

\subsection{The restricted moves}

To describe the first step of the bijection $\Phi$ in Theorem \ref{lem-b-0-over}, we will  restrict the $(k-1)$-forward move and the $(k-1)$-backward move to two subsets of $\mathcal{B}_1(\alpha_1,\ldots,\alpha_\lambda;\eta,k, r)$. We assume that $\lambda,\ k$ and $r$ are integers such that $k\geq r\geq \lambda\geq0$ and $k>\lambda$. Recall that  $ {\mathcal{B}}_j(\alpha_1,\ldots,\alpha_\lambda;\eta,k,r)$  denotes the set of overpartitions   in $\overline{\mathcal{B}}_j(\alpha_1,\ldots,\alpha_\lambda;\eta,k,r)$  without overlined parts divisible by $\eta$.
We will be concerned with the following two subsets of  $\mathcal{B}_1(\alpha_1,\ldots,\alpha_\lambda;\eta,k, r)$.

 \begin{itemize}
\item  For $N\geq p\geq 1$, let $\mathcal{{B}}_e(\alpha_1,\ldots,\alpha_\lambda;\eta,k, r|N, p)$ denote the set of  overpartitions $\gamma$ in $\mathcal{B}_1(\alpha_1,\ldots,\alpha_\lambda;\eta,k, r)$ such that there are $N$ parts marked with $k-1$ in  $RG(\gamma)$, denoted $ \tilde{r}_1(\gamma)>  \tilde{r}_2(\gamma)>\cdots > \tilde{r}_N(\gamma)$, and for all $1\leq i\leq p$, the parity of $\{\tilde{r}_i(\gamma)\}_{k-1}$ is the same   as that of $\{\tilde{r}_{p+1}(\gamma)\}_{k-1}$.

\item For $N\geq p\geq 1$, let  $\mathcal{{B}}_d(\alpha_1,\ldots,\alpha_\lambda;\eta,k, r|N, p)$ denote the set of overpartitions $\gamma$ in $\mathcal{B}_1(\alpha_1,\ldots,\alpha_\lambda;\eta,k, r)$ such that there are $N$ parts marked with $k-1$ in $RG(\gamma)$, denoted   $ \tilde{r}_1(\gamma)>  \tilde{r}_2(\gamma)>\cdots > \tilde{r}_N(\gamma)$, and for all $1\leq i\leq p$,  the parity of $\{\tilde{r}_i(\gamma)\}_{k-1}$ is opposite from  the parity of $\{\tilde{r}_{p+1}(\gamma)\}_{k-1}$.
\end{itemize}
Notice that there are no $(k-1)$-bands $\{\tilde{r}_{N+1}(\gamma)\}_{k-1}$ in $\gamma$. In this case,  we  define the parity of the empty band  to be even, and so $\mathcal{{B}}_e(\alpha_1,\ldots,\alpha_\lambda;\eta,k, r|N, N)$ is a subset of $\mathcal{B}_0(\alpha_1,\ldots,\alpha_\lambda;\eta,k,r)$.

The following theorem  shows that the forward move  $\phi_p$ gives rise to  a bijection between $ \mathcal{B}_e(\alpha_1,\ldots,\alpha_\lambda;\eta,k,r|N, p)$ and $\mathcal{{B}}_d(\alpha_1,\ldots,\alpha_\lambda;\eta,k, r|N, p) $.

\begin{thm}\label{deltagammathm} For $N\geq p\geq 1$, the forward move  $\phi_p$  is a bijection between $\mathcal{{B}}_e(\alpha_1,\ldots,\alpha_\lambda;\eta,\break$
$k, r|N, p)$  and $\mathcal{{B}}_d(\alpha_1,\ldots,\alpha_\lambda;\eta,k, r|N, p)$. Moreover,  for $\gamma \in \mathcal{{B}}_e(\alpha_1,\ldots,\alpha_\lambda;\eta,k, r|N, p)$, let $\vartheta=\phi_p(\gamma)$, we have
$|\vartheta|=|\gamma|+p\eta.$
\end{thm}

For example,   let $\gamma$ be the overpartition in ${\mathcal{B}}_1(1,5,9;10,5,4)$, whose reverse Gordon marking reads
  \begin{equation*} \label{rmark-exa-1-s-gamma}
\begin{split}
&RG(\gamma)=(\overline{81}_1,\overbrace{{ {80}_2,{80}_3,\overline{71}_1},
{ {70}_4}}^{{\{70\}_4}},\overline{69}_2,
\overbrace{{ \overline{61}_1,{60}_3,\overline{59}_2}, { \overline{55}_4}}^{{\{\overline{55}\}_4}},\overbrace{{ {50}_1,\overline{49}_2,\overline{45}_3}, { \overline{41}_4}}^{{\{\overline{41}\}_4}},\\[5pt]
&\ \ \ \ \ \ \ \ \ \ \ \ \ \  \overline{39}_1,\overline{35}_2,\underbrace{{ \overline{29}_1,\overline{21}_2,} { {20}_3},{ {20}_4}}
_{\{20\}_4},
\underbrace{{ \overline{11}_1,{10}_2,\overline{9}_3},
{ \overline{5}_4}}_{\{\overline{5}\}_4}).
\end{split}
\end{equation*}
There are five  $4$-marked parts in $RG(\gamma)$. Moreover, it can be checked that  the $4$-bands induced by $\tilde{r}_1(\gamma)={70}$, $\tilde{r}_2(\gamma)=\overline{55}$, $\tilde{r}_3(\gamma)=\overline{41}$, $\tilde{r}_4(\gamma)=20$ are all even.
 Therefore,   $\gamma$  is an overpartition in $\mathcal{B}_e(1,5,9;10,5,4|5,3)$.  Let $\vartheta=\phi_3(\gamma)$. Then the reverse Gordon marking of $\vartheta$ is given by
 \begin{equation*}\label{example-b-delta}
\begin{split}
&RG(\vartheta)=(\overbrace{{ \overline{81}_1,{80}_2,
{80}_3},{
{80}_4}}^{\{{80}\}_4},\overline{71}_1,\overbrace{{ \overline{69}_2,\overline{65}_3,\overline{61}_1},{ {60}_4}}^{\{60\}_4}, \overline{59}_2,\overbrace{\overline{51}_1,{50}_3,{ \overline{49}_2},\overline{45}_4}^{\{\overline{45}\}_4},\\[5pt]
&\ \ \ \ \ \ \ \ \ \ \ \ \ \
\overline{39}_1,\overline{35}_2,\underbrace{{ \overline{29}_1,\overline{21}_2,{20}_3},
{ {20}_4}}_{\{20\}_4},\underbrace{{ \overline{11}_1,{10}_2,
\overline{9}_3},{ \overline{5}_4}}_{\{\overline{5}\}_4}).
\end{split}
\end{equation*}
There are five  $4$-marked parts in  $RG(\vartheta)$ and   the $4$-bands induced by  $\tilde{r}_1(\vartheta)={80}$, $\tilde{r}_2(\vartheta)={60}$ and $\tilde{r}_3(\vartheta)=\overline{45}$ are  odd, whereas  the $4$-band induced by $\tilde{r}_4(\vartheta)=20$ is even. This indicates  that  $\vartheta$   is an overpartition in $\mathcal{B}_d(1,5,9;10,5,4|5,3)$. Clearly, we have $|\vartheta|=|\gamma|+30$.

To prove Theorem \ref{deltagammathm}, we establish two lemmas. From now on,  we shall use $f_{\leq \eta}(\gamma)$  to denote the number of parts less than or equal to $\eta$ in an overpartition $\gamma$.

\begin{lem}\label{deltagamma} For $N\geq p\geq 1$,  let $\gamma$ be an overpartition in $\mathcal{{B}}_e(\alpha_1,\ldots,\alpha_\lambda;\eta,k, r|N, p)$ and let  $\vartheta=\phi_p(\gamma)$. Then  $\vartheta$ is an overpartition in $\mathcal{{B}}_d(\alpha_1,\ldots,\alpha_\lambda;\eta,k, r|N, p)$. Moreover,   $|\vartheta|=|\gamma|+p\eta$.
\end{lem}

\proof  Clearly,  $\gamma$ is an overpartition in $\mathcal{{B}}_1(\alpha_1,\ldots,\alpha_\lambda;\eta,k, r)$ with $N$ parts marked with $k-1$ in $RG(\gamma)$. In view of Proposition \ref{lem-for1}, we find that $\vartheta=\phi_p(\gamma)$ satisfies (1), (2) and (3) in Definition \ref{defi-O-B}. Furthermore, there are $N$ parts marked with $k-1$ in $RG(\vartheta)$. Thus, to prove that $\vartheta$ belongs to $\mathcal{{B}}_d(\alpha_1,\ldots,\alpha_\lambda;\eta,k, r|N, p)$, it suffices to verify that the following conditions hold:
\begin{itemize}
    \item[(A)] $f_{\leq \eta}(\vartheta)\leq r-1$;
    \item[(B)] For $1\leq i\leq p$, $\{\tilde{r}_{i}(\vartheta)\}_{k-1}$ and $\{\tilde{r}_{p+1}(\vartheta)\}_{k-1}$ have opposite  parities.
\end{itemize}

\noindent Condition (A).  It is readily seen that $f_{\leq \eta}(\vartheta)$ equals either $f_{\leq \eta}(\gamma)$ or $f_{\leq \eta}(\gamma)-1$. Under the condition   $f_{\leq \eta}(\gamma)\leq r-1$, we get $f_{\leq \eta}(\vartheta)\leq r-1$.

\noindent Condition (B). By the property  (3) in  Proposition \ref{lem-for1},  there are $N$ parts marked with $k-1$ in $G(\vartheta)$, denoted $ \tilde{g}_1(\vartheta)>  \tilde{g}_2(\vartheta)>\cdots > \tilde{g}_N(\vartheta)$. It follows from  Proposition \ref{sequence-length} that there are also $N$ parts marked with $k-1$ in $RG(\vartheta)$, denoted $ \tilde{r}_1(\vartheta)>  \tilde{r}_2(\vartheta)>\cdots > \tilde{r}_N(\vartheta)$ such that  $\tilde{g}_i(\vartheta)\in \{\tilde{r}_i(\vartheta)\}_{k-1}$  for $1\leq i \leq N$. This implies that $ \{\tilde{r}_i(\vartheta)\}_{k-1}$ is a $(k-1)$-band including $\tilde{g}_i(\vartheta)$. Utilizing Proposition  \ref{parity k-1 sequence-over}, we obtain  that for each $1\leq i \leq N$,   $\{\tilde{r}_{i}(\vartheta)\}_{k-1}$ and $\{\tilde{g}_i(\vartheta)\}_{k-1}$ have the same parity.  Therefore, to prove that   $\{\tilde{r}_{i}(\vartheta)\}_{k-1}$ and $\{\tilde{r}_{p+1}(\vartheta)\}_{k-1}$ have opposite  parities  for $1\leq i\leq p$, we are obliged to  show  that $\{\tilde{g}_{i}(\vartheta)\}_{k-1}$ and   $\{\tilde{g}_{p+1}(\vartheta)\}_{k-1}$  have opposite  parities for $1\leq i\leq p$.

For $1\leq i\leq N$, let
\[\tilde{r}_{i,1}(\gamma)  \geq  \cdots  \geq  \tilde{r}_{i,k-2}(\gamma)\geq \tilde{r}_{i}(\gamma)\]
be the parts in the $(k-1)$-band of $\gamma$ induced by $\tilde{r}_{i}(\gamma)$, and  let
\[ \tilde{g}_{i}(\vartheta)  \geq \tilde{g}_{i,2}(\vartheta)  \ge  \cdots  \geq  \tilde{g}_{i,k-1}(\vartheta)\]
be the parts in the $(k-1)$-band of $\vartheta$ induced by $\tilde{g}_{i}(\vartheta)$. Write
\[[\tilde{r}_{i,1}(\gamma)/\eta]+\cdots+
[\tilde{r}_{i,k-2}(\gamma)/\eta]+[\tilde{r}_{i}(\gamma)/\eta]\equiv a_i(\gamma)+\overline{V}_\gamma(\tilde{r}_{i,1}(\gamma))\pmod2,\]
and
\[[\tilde{g}_{i}(\vartheta)/\eta]+[\tilde{g}_{i,2}(\vartheta)/\eta]+\cdots+[\tilde{g}_{i,k-1}(\vartheta)/\eta]\equiv a_i(\vartheta)+\overline{V}_\vartheta(\tilde{g}_{i}(\vartheta))\pmod2,\]
where $a_i(\gamma)$ (resp. $a_i(\vartheta)$) either equals $r-1$ or $r$ for $1\leq i\leq N$ and $a_{N+1}(\gamma)$ (resp. $a_{N+1}(\vartheta)$) $=r-1$ with the  convention that  the empty band  is even.

Since  $\gamma$ belongs to $ {\mathcal{B}}_e(\alpha_1,\ldots,\alpha_\lambda;\eta,k,r|N,p)$, where $1\leq p\leq N$,
we have for $1\leq i\leq p$,
\begin{equation}\label{3-bij-lem1tttgamma}
  a_i(\gamma)=a_{p+1}(\gamma).
\end{equation}
We proceed to show that  $\{\tilde{g}_{i}(\vartheta)\}_{k-1}$ and   $\{\tilde{g}_{p+1}(\vartheta)\}_{k-1}$  have   opposite parities for $1\leq i\leq p$, or equivalently,   for $1\leq i\leq p$,
\begin{equation}\label{3-bij-lem1ttt}
a_{i}(\vartheta)\not = a_{p+1}(\vartheta).
\end{equation}

The proof of Proposition \ref{lem-for1} justifies  the following relation for $1\leq i\leq p$,
 \begin{equation}\label{gammadelta-r}
 \begin{array}{ccccccc}
  \tilde{r}_{i}(\gamma)+\eta&\geq& \tilde{r}_{i,1}(\gamma) &\geq& \cdots& \geq &\tilde{r}_{i,k-2}(\gamma)\\[5pt]
  \shortparallel& &\shortparallel & &  &  &\shortparallel\\[5pt]
   \tilde{g}_{i}(\vartheta)& \geq &\tilde{g}_{i,2}(\vartheta) &\geq& \cdots& \geq &\tilde{g}_{i,k-1}(\vartheta).
  \end{array}
\end{equation}
We claim that for $1\leq i\leq p$,
\begin{equation}\label{congrue-condition-r}
\overline{V}_{\vartheta}(\tilde{g}_{i}(\vartheta))=\overline{V}_{\gamma}(\tilde{r}_{i,1}(\gamma)).
\end{equation}
Recall that $\overline{V}_\pi(t)$ (resp.   $\overline{V}_\pi(\overline{t})$)  stands for the number of overlined parts not exceeding $t$ (resp.  $\overline{t}$) in $\pi$.

Owing to the relation \eqref{gammadelta-r}, we deduce that for $1\leq i\leq p$,
\begin{eqnarray}\label{congrue-condition-r1}
\overline{V}_\gamma(\tilde{r}_{i,1}(\gamma))-\overline{V}_\vartheta(\tilde{g}_{i,2}
(\vartheta))
=\begin{cases}
1, \quad \text{if} \  \tilde{r}_i(\gamma)\not\equiv 0 \bmod{\eta}, \\[5pt]
0, \quad \text{otherwise,} \quad
\end{cases}
\end{eqnarray}
and
\begin{eqnarray}\label{congrue-condition-r2}
\overline{V}_\vartheta(\tilde{g}_{i}(\vartheta))-\overline{V}_\vartheta(\tilde{g}_{i,2}(\vartheta))
=\begin{cases}
1, \quad \text{if} \  \tilde{g}_i(\vartheta)
\not\equiv 0 \bmod{\eta},  \\[5pt]
0, \quad \text{otherwise.} \quad
\end{cases}
\end{eqnarray}

By definition,  $\tilde{g}_i(\vartheta)=\tilde{r}_i(\gamma)+\eta$, and so $\tilde{g}_i(\vartheta)$ is divisible by $\eta$ if and only if $\tilde{r}_i(\gamma)$ is divisible by $\eta$. Therefore,
combining  \eqref{congrue-condition-r1}  and \eqref{congrue-condition-r2}   gives \eqref{congrue-condition-r}, and hence  the claim is proved.

Invoking the relation \eqref{gammadelta-r}, we find that for $1\leq i\leq p$,
\begin{equation*}\label{gammadelta-rrr-000}
\begin{split}
&[\tilde{g}_{i}(\vartheta)/\eta]+[\tilde{g}_{i,2}(\vartheta)/\eta]+\cdots+[\tilde{g}_{i,k-1}(\vartheta)/\eta]\\[5pt]
&\qquad=[(\tilde{r}_{i}(\gamma)+\eta)/\eta]+[\tilde{r}_{i,1}(\gamma)/\eta]+\cdots+[\tilde{r}_{i,k-2}(\gamma)/\eta]\\[5pt]
&\qquad=[\tilde{r}_{i}(\gamma)/\eta]+[\tilde{r}_{i,1}(\gamma)/\eta]+\cdots+
[\tilde{r}_{i,k-2}(\gamma)/\eta]+1\\[5pt]
&\qquad \equiv a_i(\gamma)+\overline{V}_\gamma(\tilde{r}_{i,1}(\gamma))+1\pmod2.
\end{split}
\end{equation*}
It follows from  \eqref{congrue-condition-r} that $a_{i}(\vartheta)\equiv a_{i}(\gamma)+1 \pmod{2}$ for $1\leq i\leq p$. In view of \eqref{3-bij-lem1tttgamma},
we obtain that for $1\leq i\leq p$,
\begin{equation}\label{3-bij-lem1}
a_{i}(\vartheta) \equiv a_{p+1}(\gamma) +1 \pmod{2}.
\end{equation}
We next show that
\begin{equation}\label{ttaam}
a_{p+1}(\gamma)= a_{p+1}(\vartheta).
\end{equation}
From Proposition  \ref{parity k-1 sequence-over}, we know that   the parity of $\{\tilde{r}_{p+1}(\vartheta)\}_{k-1}$ is the same as that of $\{\tilde{g}_{p+1}(\vartheta)\}_{k-1}$. On the other hand, the construction of the forward move $\phi_p$ indicates that the parity of $\{\tilde{r}_{p+1}(\gamma)\}_{k-1}$ is the same as that of $\{\tilde{r}_{p+1}(\vartheta)\}_{k-1}$, and
so the parity of $\{\tilde{r}_{p+1}(\gamma)\}_{k-1}$ agrees with
that of  $\{\tilde{g}_{p+1}(\vartheta)\}_{k-1}$. Thereby, we get
\eqref{ttaam}.   Combining     \eqref{3-bij-lem1} and \eqref{ttaam} gives \eqref{3-bij-lem1ttt}. It follows that $\{\tilde{g}_{i}(\vartheta)\}_{k-1}$  and  $\{\tilde{g}_{p+1}(\vartheta)\}_{k-1}$  have opposite  parities for $1\leq i\leq p$, and so  $\{\tilde{r}_{i}(\vartheta)\}_{k-1}$  and  $\{\tilde{r}_{p+1}(\vartheta)\}_{k-1}$  have opposite  parities for $1\leq i\leq p$. Hence the condition (B) is satisfied.

We have shown that     $\vartheta\in{\mathcal{B}}_d(\alpha_1,\ldots,\alpha_\lambda;\eta,k,r|
  N,p)$. It is routine to verify that $|\vartheta|=|\gamma|+p\eta$, and thus the proof is complete. \qed

\begin{lem}\label{deltagamma-reverse} For $N\geq p\geq 1$, let $\vartheta$ be an overpartition in $\mathcal{{B}}_d(\alpha_1,\ldots,\alpha_\lambda;\eta,k, r|N, p)$ and let $\gamma=\psi_p(\vartheta)$. Then $\gamma$ is an overpartition in $\mathcal{{B}}_e(\alpha_1,\ldots,\alpha_\lambda;\eta,k, r|N, p)$.   Furthermore,  $|\gamma|=|\vartheta|-p\eta$.
\end{lem}

\pf In order to show that $\gamma$ is an overpartition in $\mathcal{{B}}_e(\alpha_1,\ldots,\alpha_\lambda;\eta,k, r|N, p)$, we need to prove that $\gamma$ satisfies (1), (2) and (3) in Definition \ref{defi-O-B} and  there are $N$ parts marked with $k-1$ in $RG(\gamma)$, denoted $\tilde{r}_{1}(\gamma)>\tilde{r}_{2}(\gamma)>\cdots>\tilde{r}_{N}(\gamma)$. Moreover, the following conditions are also required:
\begin{itemize}
\item[(A)] $f_{\leq \eta}(\gamma)\leq r-1$;
\item[(B)] the parity of $\{\tilde{r}_{i}(\gamma)\}_{k-1}$ is the same as that of $\{\tilde{r}_{p+1}(\gamma)\}_{k-1}$ for $1\leq i\leq p$.
\end{itemize}

Now we consider (1), (2) and (3) in Definition  \ref{defi-O-B}.
Assume that $ \tilde{g}_1(\vartheta)>  \tilde{g}_2(\vartheta)>\cdots > \tilde{g}_N(\vartheta)$ are the $(k-1)$-marked parts in the Gordon marking of $\vartheta \in \mathcal{{B}}_d(\alpha_1,\ldots,\alpha_\lambda;\eta,k, r|N, p)$.  By   Proposition \ref{lem-bac1},  it is necessary to  prove that
\begin{itemize}
\item[(C)] $\tilde{g}_p(\vartheta)\geq\overline{\eta+\alpha_1}$;
\item[(D)] $\tilde{g}_p(\vartheta)$ is a part in any $(k-1)$-band   of $\vartheta$  belonging to $I(\tilde{g}_p(\vartheta)-2\eta, \tilde{g}_p(\vartheta))$.
\end{itemize}

 \noindent Condition (C). Given that $\vartheta \in \mathcal{{B}}_d(\alpha_1,\ldots,\alpha_\lambda;\eta,k, r|N, p)$,   combining Proposition \ref{sequence-length} and Proposition  \ref{parity k-1 sequence-over}, we realize  that  $ \{\tilde{g}_i(\vartheta)\}_{k-1}$ and  $ \{\tilde{g}_{p+1}(\vartheta)\}_{k-1}$ have opposite  parities  for $1\leq i\leq p$.  Suppose to the contrary that  $\tilde{g}_p(\vartheta)<\overline{\eta+\alpha_1}$, which means that $\tilde{g}_p(\vartheta)\leq \eta$.
In this case, we have $p=N$. Observing that $\tilde{g}_N(\vartheta)$ is marked with $k-1$ in $G(\vartheta)$,  so we get  $f_{\leq \eta}(\vartheta)=k-1$, that is, $r=k$.
Assume that
\[ \tilde{g}_{N}(\vartheta)  \geq \tilde{g}_{N,2}(\vartheta)  \ge  \cdots  \geq  \tilde{g}_{N,k-1}(\vartheta)\]
are the parts in the $(k-1)$-band of $\vartheta$ induced by $\tilde{g}_{N}(\vartheta)$. Under the condition that  $ \{\tilde{g}_N(\vartheta)\}_{k-1}$ and $ \{\tilde{g}_{N+1}(\vartheta)\}_{k-1}$ have opposite  parities   and the convention that  the empty band is even, we deduce that $\{\tilde{g}_N(\vartheta)\}_{k-1}$ is odd, that is,
\begin{equation}\label{delta-contra-333}
\begin{split}
&\ \ \ \ [\tilde{g}_{N}(\vartheta)/\eta]+[\tilde{g}_{N,2}(\vartheta)/\eta]+\cdots+[\tilde{g}_{N,k-1}(\vartheta)/\eta]\equiv r+\overline{V}_\vartheta(\tilde{g}_{N}(\vartheta))\pmod2.
\end{split}
\end{equation}
On the other hand, since $\tilde{g}_N(\vartheta)\leq \eta$, we obtain that
\begin{equation} \label{delta-contra-333t}
\begin{split}
&  [\tilde{g}_{N}(\vartheta)/\eta]+[\tilde{g}_{N,2}(\vartheta)/\eta]+\cdots
+[\tilde{g}_{N,k-1}(\vartheta)/\eta]\\[5pt]
&\quad =f_{\leq \eta}(\vartheta)-f_{<\eta}(\vartheta)=k-1-f_{<\eta}(\vartheta),
\end{split}
\end{equation}
where $f_{<\eta}(\vartheta)$    denotes the number of parts of $\vartheta$ less than $\eta$. Recall that $\overline{V}_\vartheta(\tilde{g}_{N}(\vartheta))$ counts   the number of overlined parts of $\vartheta$  not exceeding $\tilde{g}_{N}(\vartheta)$. Again,  under the assumption  $\tilde{g}_N(\vartheta)\leq \eta$, we have
$\overline{V}_\vartheta(\tilde{g}_{N}(\vartheta))=f_{<\eta}(\vartheta)$. Since $r=k$,  \eqref{delta-contra-333t} can be written as
\[ [\tilde{g}_{N}(\vartheta)/\eta]+[\tilde{g}_{N,2}(\vartheta)/\eta]+\cdots+[\tilde{g}_{N,k-1}(\vartheta)/\eta]=
r-1-\overline{V}_\vartheta(\tilde{g}_{N}(\vartheta)),\]
which contradicts  \eqref{delta-contra-333}. Hence $\tilde{g}_p(\vartheta)\geq \overline{\eta+\alpha_1}$.

 \noindent Condition (D).  Suppose that there is  a $(k-1)$-band belonging to $I(\tilde{g}_p(\vartheta)-2\eta, \tilde{g}_p(\vartheta))$ which does not contain $\tilde{g}_p(\vartheta)$ as a part, and let
\[ \vartheta_m\geq\vartheta_{m+1}\geq \cdots \geq \vartheta_{m+k-2}\]
be the parts in this $(k-1)$-band, that is,  $ \vartheta_{m}\leq \vartheta_{m+k-2}+\eta $ with strict inequality if $\vartheta_m$ is overlined, $\vartheta_m<\tilde{g}_p(\vartheta)$ and   $\vartheta_{m+k-2}\geq \tilde{g}_p(\vartheta)-2\eta$   with strict inequality if $\tilde{g}_p(\vartheta)$ is overlined. In view of   Lemma \ref{parity k-1 sequence-over-old}, we deduce that $\{\tilde{g}_p(\vartheta)\}_{k-1}$ and $\{\vartheta_{m+l}\}_{0\leq l\leq k-2}$ are of the same parity.

Now, since $\{\vartheta_{m+l}\}_{0\leq l\leq k-2}$ is a $(k-1)$-band of $\vartheta$, there is a part, say $\vartheta_{m+t}$ ($0\leq t\leq k-2$),  marked with $k-1$ in $G(\vartheta)$. But $\vartheta_m<\tilde{g}_p(\vartheta)$, so we get $\vartheta_{m+t}=\tilde{g}_{p+1}(\vartheta)$. This implies that $\{\vartheta_{m+l}\}_{0\leq l\leq k-2}$ is a $(k-1)$-band of $\vartheta$ including $\tilde{g}_{p+1}(\vartheta)$. By Proposition  \ref{parity k-1 sequence-over}, we deduce that $\{\tilde{g}_{p+1}(\vartheta)\}_{k-1}$ and  $\{\vartheta_{m+l}\}_{0\leq l\leq k-2}$ are of the same parity. It follows that  the parity of $\{\tilde{g}_{p+1}(\vartheta)\}_{k-1}$ is the same as that of  $\{\tilde{g}_{p}(\vartheta)\}_{k-1}$,   which contradicts  the condition that $\vartheta \in \mathcal{{B}}_d(\alpha_1,\ldots,\alpha_\lambda;\eta,k, r|N, p)$, that is, the parity of $\{\tilde{g}_{p+1}(\vartheta)\}_{k-1}$ is opposite from  the parity of $\{\tilde{g}_{p}(\vartheta)\}_{k-1}$. Therefore,   $\tilde{g}_p(\vartheta)$ is a part in any $(k-1)$-band   of $\vartheta$  belonging to $I(\tilde{g}_p(\vartheta)-2\eta, \tilde{g}_p(\vartheta))$.

Up to now, we have shown that $\vartheta$ satisfies the conditions (C) and (D). In view of  Proposition \ref{lem-bac1}, we see that $\gamma$ satisfies (1), (2) and (3) in Definition \ref{defi-O-B} and there are $N$ parts marked with $k-1$ in $RG(\gamma)$. We still have to show that $\gamma$ satisfies the conditions (A) and (B).

 \noindent Condition (A). By the condition (C), we have $\tilde{g}_p(\vartheta)\geq \overline{\eta+\alpha_1}$. We now consider two cases: (1) If $\tilde{g}_p(\vartheta)>2\eta$, then $f_{\leq \eta}(\gamma)=f_{\leq \eta}(\vartheta)\leq r-1$.
(2) If $\overline{\eta+\alpha_1}\leq \tilde{g}_p(\vartheta)\leq 2\eta$,  then  $f_{\leq \eta}(\gamma)=f_{\leq \eta}(\vartheta)+1$. In this case, we claim that $f_{\leq \eta}(\vartheta)< r-1$. Suppose to the contrary that $f_{\leq \eta}(\vartheta)= r-1$. Assume that
\[ \tilde{g}_{p}(\vartheta)  \geq \tilde{g}_{p,2}(\vartheta)  \ge  \cdots  \geq  \tilde{g}_{p,k-1}(\vartheta)\]
are the parts in the $(k-1)$-band of $\vartheta$ induced by $\tilde{g}_{p}(\vartheta)$. Then
\begin{equation*}\label{zp-re}
\begin{split}
&[\tilde{g}_{p}(\vartheta)/\eta]+[\tilde{g}_{p,2}(\vartheta)/\eta]+\cdots+[\tilde{g}_{p,k-1}(\vartheta)/\eta]\\
&\qquad \equiv \overline{V}_\vartheta(\tilde{g}_{p}(\vartheta))-f_{<\eta}(\vartheta)+f_{\eta}(\vartheta)\\
&\qquad \equiv f_{\leq \eta}(\vartheta)+\overline{V}_\vartheta(\tilde{g}_{p}(\vartheta))\pmod2,
\end{split}
\end{equation*}
which implies that $\{\tilde{g}_{p}(\vartheta)\}_{k-1}$ is even since $f_{\leq \eta}(\vartheta)=r-1$. Given that $\vartheta \in \mathcal{{B}}_d(\alpha_1,\ldots,\alpha_\lambda;\eta,$\break $k, r|N, p)$,  we see that the parity of $\{\tilde{g}_{p}(\vartheta)\}_{k-1}$ is opposite from  the parity of $\{\tilde{g}_{p+1}(\vartheta)\}_{k-1}$. It follows that $\{\tilde{g}_{p+1}(\vartheta)\}_{k-1}$ is odd, and so $\{\tilde{g}_{p+1}(\vartheta)\}_{k-1}$ is nonempty.   On the other hand, since  $\tilde{g}_{p+1}(\vartheta)\leq \tilde{g}_p(\vartheta)-\eta\leq\eta$, it is ensured by Lemma  \ref{parity k-1 sequence-over-old} that   the parity of $\{\tilde{g}_{p}(\vartheta)\}_{k-1}$ is the same as that of $\{\tilde{g}_{p+1}(\vartheta)\}_{k-1}$,   which leads to a contradiction. Hence  $f_{\leq \eta}(\vartheta)<r-1$ when $\overline{\eta+\alpha_1}\leq \tilde{g}_p(\vartheta)\leq 2\eta$, and so $f_{\leq \eta}(\gamma)\leq r-1$.

  \noindent Condition (B). Utilizing the property (4) in  Proposition \ref{lem-bac1},  we find that  for $1\leq i\leq p$, $\tilde{r}_i(\gamma)=\tilde{g}_i(\vartheta)-\eta$.  The reasoning in the proof of Lemma \ref{deltagamma} can  be adapted to deduce
  that the parity of $\{\tilde{r}_{i}(\gamma)\}_{k-1}$ is the same as that of  $\{\tilde{r}_{p+1}(\gamma)\}_{k-1}$ for $1\leq i\leq p$.

 Thus we conclude  that $\gamma$ is an overpartition in $\mathcal{{B}}_e(\alpha_1,\ldots,\alpha_\lambda;\eta,k, r|N, p)$. It is manifest from the construction of $\psi_p$  that
$ |\gamma|=|\vartheta|-p\eta$. This  completes the proof.     \qed

\noindent{\it Proof of Theorem \ref{deltagammathm}.}  Let $\gamma \in \mathcal{{B}}_e(\alpha_1,\ldots,\alpha_\lambda;\eta,k, r|N, p)$. Utilizing Lemma \ref{deltagamma}, we find that $\phi_p(\gamma)$ belongs to $\mathcal{{B}}_d(\alpha_1,\ldots,\alpha_\lambda;\eta,k, r|N, p)$. In view of the property (4) in Proposition \ref{lem-for1}, we deduce that $\psi_p(\phi_p(\gamma))=\gamma$.

Analogously, let $\vartheta \in \mathcal{{B}}_d(\alpha_1,\ldots,\alpha_\lambda;\eta,k, r|N, p)$. Invoking  Lemma \ref{deltagamma-reverse}, we get   $\psi_p(\vartheta) \in \mathcal{{B}}_e(\alpha_1,\ldots,\alpha_\lambda;\eta,k, r|N, p)$. By virtue of the property (4) in Proposition \ref{lem-bac1}, we obtain that  $\phi_p(\psi_p(\vartheta))=\vartheta$.

Thus, we arrive at the assertion  that the forward move $\phi_p$ is a bijection between $\mathcal{{B}}_e(\alpha_1,\ldots,\alpha_\lambda;\eta,k, r|N, p)$ and $\mathcal{{B}}_d(\alpha_1,\ldots,\alpha_\lambda;\eta,k, r|N, p)$. This completes the proof. \qed


 \subsection{The $(k-1)$-insertion  and the $(k-1)$-separation}
As mentioned before, a merging operation is needed in the construction  of the  bijection $\Phi$ between $\mathcal{D}_\eta\times{\mathcal{B}}_{0}(\alpha_1,\ldots,\alpha_\lambda;\eta,k,r)$
and $\mathcal{\overline{B}}_1(\alpha_1,\ldots,\alpha_\lambda;\eta,k,r)$. The main objective of this subsection is to present a
description  of this merging operation in terms of the $(k-1)$-
insertion operation and the $(k-1)$-separation operation. To be more specific,  the merging operation is  meant to  take the parts divisible by $\eta$ and the parts of the overpartition in ${\mathcal{B}}_1(\alpha_1,\ldots,\alpha_\lambda;\eta,k,r)$ to generate certain overlined parts divisible by $\eta$. As a result, we get an overpartition $\pi$  in $\mathcal{\overline{B}}_1(\alpha_1,\ldots,\alpha_\lambda;\eta,k,r)$.  To this end, we shall prepare two subsets of  $\overline{\mathcal{B}}_1(\alpha_1,\ldots,\alpha_\lambda;\eta,k,r)$. Assume that     $a=\eta$ or $\alpha_i$ for some $1\leq i\leq \lambda$.

 \begin{itemize}
 \item    For  $s\geq N\geq 0$, let $\overline{\mathcal{B}}^{\, a}_<(\alpha_1,\ldots,\alpha_\lambda;\eta,k,r|N,s)$ denote the set of   overpartitions $\tau$ in $\overline{\mathcal{B}}_1(\alpha_1,\ldots,\alpha_\lambda;\eta,k,r)$  satisfying
   \begin{itemize}
 \item[{\rm (1)}] There are $N$ parts marked with $k-1$ in   $RG(\tau)$, denoted $ \tilde{r}_1(\tau)>  \tilde{r}_2(\tau)>\cdots > \tilde{r}_N(\tau)${\rm{;}}

  \item[{\rm (2)}]  Assume that  $p$ is the smallest integer satisfying  $ \tilde{r}_{p+1}(\tau)+\eta \leq \overline{(s-p)\eta+a}$ with the convention that $\tilde{r}_{N+1}(\tau)=-\infty$. Then the largest overlined part $\equiv a\pmod\eta$ in $\tau$ is less than $\overline{(s-p)\eta+a}${\rm{;}}

 \item[{\rm (3)}]   If  $f_{\leq \eta}(\tau)=r-1$, $s=N\geq 1$ and  $a\neq\eta$, then  $\tilde{r}_{N}(\tau)\leq\eta${\rm{;}}

      \item[{\rm (4)}]  If   $s=N=0$ and  $a\neq \eta$, then $f_{\leq \eta}(\tau)<r-1$.

  \end{itemize}

  \item For $s\geq N\geq 0$, let $\overline{\mathcal{B}}^{\, a}_{=}(\alpha_1,\ldots,\alpha_\lambda;\eta,k,r|N,s)$ denote the set of overpartitions $\sigma$ in $\overline{\mathcal{B}}_1(\alpha_1,\ldots,\alpha_\lambda;\eta,k,r)$    subject to the following conditions:
     \begin{itemize}
  \item[{\rm (1)}] There exists  an  overlined part $\equiv a\pmod\eta$ in $\sigma$, and assume that the  largest overlined part $\equiv a\pmod\eta$  in  $\sigma$ is $\overline{t\eta+a}$;

    \item[{\rm (2)}] Let ${\hat{\sigma}}$ be the  overpartition obtained  by removing
  $\overline{t\eta+a}$  from $\sigma$.  Then there are $N$ parts marked with $k-1$ in  $G({\hat{\sigma}})$,  denoted $ \tilde{g}_1({\hat{\sigma}})>  \tilde{g}_2({\hat{\sigma}})>\cdots > \tilde{g}_N({\hat{\sigma}})$;

    \item[{\rm (3)}] Assume that  $p$  is the smallest integer such that  $ \tilde{g}_{p+1}({\hat{\sigma}})< \overline{t\eta+a}$ with the convention that $\tilde{g}_{N+1}({\hat{\sigma}})=-\infty$.  Then $s=p+t$.

      \end{itemize}

\end{itemize}

For example,  let $N=5$, $s=6$ and $a=10$ and  let $\tau$ be the overpartition in $\overline{\mathcal{B}}_1(1,5,9;10,5,4)$ with the reverse Gordon marking
  \begin{equation}\label{example-insertion-over-1}
\begin{split}
&RG(\tau)=(\overline{85}_1,{ {80}_2,{80}_3,\overline{75}_1},{ {70}_4}
,\overline{69}_2,
{ \overline{61}_1,{60}_3,\overline{59}_2}, { \overline{55}_4},{ {50}_1,\overline{49}_2,\overline{45}_3}, { \overline{41}_4},\\[5pt]
&\ \ \ \ \ \ \ \ \ \ \ \ \ \  \overline{39}_1,\overline{35}_2,{ \overline{29}_1,\overline{20}_2,} { {20}_3},{ {20}_4},
{ \overline{11}_1,\overline{10}_2,\overline{9}_3},
{ \overline{5}_4},\overline{1}_1).
\end{split}
\end{equation}
There are five  $4$-marked parts  $\tilde{r}_1(\tau)={70}$, $\tilde{r}_2(\tau)=\overline{55}$, $\tilde{r}_3(\tau)=\overline{41}$, $\tilde{r}_4(\tau)=20$ and $\tilde{r}_5(\tau)=\overline{5}$  in  $RG(\tau)$.   Then $p=3$ is the smallest integer such that $ 30= \tilde{r}_{p+1}(\tau)+\eta\leq \overline{(s-p)\eta+a}=\overline{40}$. Meanwhile,  the largest overlined part divisible by $10$ in $\tau$ is $\overline{20}$, which is less than $\overline{(s-p)\eta+a}=\overline{40}$.  So $\tau$ is an overpartition in $\overline{\mathcal{B}}^{\,10}_<(1,5,9;10,5,4|5,6)$.

Next example is   concerned   with  determining whether an overpartition in $\overline{\mathcal{B}}_1(\alpha_1,\ldots,\alpha_\lambda;\break $$\eta,k,r)$ belongs to   $\overline{\mathcal{B}}^{\, a}_{=}(\alpha_1,\ldots,\alpha_\lambda;\eta,k,r|N,s)$.  Let $N=5$, $s=6$, $a=10$ and let
  \begin{equation}\label{example-insertion-over-2}
\begin{split}
&\sigma=(\overline{85},{{80}}
,{ {80},{80},\overline{75}},\overline{69},
 {\overline{65}},{ \overline{61},{60},\overline{59}},{\overline{51}},{ {50},\overline{49},\overline{45}},{ \overline{40}},\\[5pt]
&\ \ \ \ \ \ \
\overline{39},\overline{35},{ \overline{29},\overline{20},} { {20}},{{20}},
{ \overline{11},\overline{10},\overline{9}},
{\overline{5}},\overline{1})
\end{split}
\end{equation}
  be an overpartition in $\overline{\mathcal{B}}_1(1,5,9;10,5,4)$.  The largest overlined part divisible by $10$ of $\sigma$ is  $\overline{40}$, and so $t=3$. Removing  $\overline{40}$  from $\sigma$, we get ${\hat{\sigma}}$ with the Gordon marking  \begin{equation}\label{example-insertion-over-3}
\begin{split}
&G({\hat{\sigma}})=(\overline{85}_2,{ {80}_4}
,{ {80}_3,{80}_1,\overline{75}_2},\overline{69}_1,
 { \overline{65}_4},{ \overline{61}_3,{60}_2,\overline{59}_1},{ \overline{51}_4},{ {50}_3,\overline{49}_1,\overline{45}_2},\\[5pt]
&\ \ \ \ \ \ \ \ \ \ \ \
\overline{39}_1,\overline{35}_2,{ \overline{29}_1,{ \overline{20}_4},} { {20}_3},{{20}_2},
{ \overline{11}_1,{ \overline{10}_4},\overline{9}_3},
{\overline{5}_2},\overline{1}_1).
\end{split}
\end{equation}
There are five  $4$-marked parts  $\tilde{g}_1({\hat{\sigma}})={80}$, $\tilde{g}_2({\hat{\sigma}})=\overline{65}$, $\tilde{g}_3({\hat{\sigma}})=\overline{51}$, $\tilde{g}_4({\hat{\sigma}})=\overline{20}$ and $\tilde{g}_5({\hat{\sigma}})=\overline{10}$  in  $G({\hat{\sigma}})$ and $p=3$ is the smallest integer such that
$\overline{20}=\tilde{g}_{p+1}({\hat{\sigma}})<\overline{40}$. Indeed, $p+t=s$ holds.   Thus, we conclude that $\sigma$ is an overpartition in $\overline{\mathcal{B}}^{\,10}_=(1,5,9;10,5,4|5,6)$.

We next give the definition of the $(k-1)$-insertion operation, which serves as a bijection between $\mathcal{\overline{B}}^{\, a}_<(\alpha_1,\ldots,\alpha_\lambda;\eta,k, r|N,s)$  and   $\mathcal{\overline{B}}^{\, a}_=(\alpha_1,\ldots,\alpha_\lambda;\eta,k, r|N,s)$.

\begin{defi}[The $(k-1)$-insertion]\label{the insertion} For  $s\geq N\geq 0$,  let $\tau$ be an overpartition in $\overline{\mathcal{{B}}}^{\,a}_<(\alpha_1,\ldots,\alpha_\lambda;\eta,k, r|N,s)$ with $N$ parts marked with $k-1$ in $RG(\tau)$, denoted  $\tilde{r}_1(\tau)>\cdots>\tilde{r}_N(\tau)$. Assume that  $p$ is the smallest integer such that   $0\leq p\leq N$ and $\overline{(s-p)\eta+a}\geq\tilde{r}_{p+1}(\tau)+\eta$. The $(k-1)$-insertion $I^a_s\colon \tau \rightarrow \sigma$ is defined  as follows:   first apply the forward move $\phi_p$ to $\tau$ to get ${{\tau}'}=\phi_p(\tau)$,  then insert $\overline{(s-p)\eta+a}$ into $\tau'$  as an overlined part of $\sigma$.
 \end{defi}
It should be understood that when $p=0$,  the forward move $\phi_p$ is considered as   the identity map, that is, $\phi_p(\tau)=\tau$. In this paper, we adopt the $(k-1)$-insertion with $a=\eta$. The case $a=\alpha_1$ will be used in our second paper \cite{He-Ji-Zhao}.

 For example, take  the overpartition $\tau$ in   $\overline{\mathcal{B}}^{\,10}_<(1,5,9;10,5,4|5,6)$ whose reverse Gordon marking is  given in \eqref{example-insertion-over-1}. In this case,  $p=3$   is the smallest integer such that $\overline{(s-p)\eta+a}=\overline{40}\geq 30= \tilde{r}_{p+1}(\tau)+\eta$, where $s=6$ and $a=10$.   Applying the forward move  $\phi_3$ to $\tau$, we  get
 \begin{equation*}
\begin{split}
&{\tau'}=(\overline{85},{{80}}
,{ {80},{80},\overline{75}},\overline{69},
 {\overline{65}},{ \overline{61},{60},\overline{59}},{\overline{51}},{ {50},\overline{49},\overline{45}},\\[5pt]
&\ \ \ \ \ \ \
\overline{39},\overline{35},{ \overline{29},\overline{20},} { {20}},{{20}},
{ \overline{11},\overline{10},\overline{9}},
{\overline{5}},\overline{1}),
\end{split}
\end{equation*}
 whose Gordon marking agrees with the one in \eqref{example-insertion-over-3}. Inserting $\overline{(s-p)\eta+a}=\overline{40}$  into  ${{\tau}'}$, we obtain $\sigma=I_s^a(\tau)$ as in \eqref{example-insertion-over-2}, which  belongs to $\overline{\mathcal{B}}^{10}_=(1,5,9;10,5,4|5,6)$. Clearly,  $|\sigma|=|\tau|+70$.

\begin{thm}\label{deltagammathmbb} For  $s\geq N\geq 0$,  the $(k-1)$-insertion  $ I^{a}_{s}$   is a bijection between   $\mathcal{\overline{B}}^{\, a}_<(\alpha_1,\ldots,\break$
$\alpha_\lambda;\eta,k, r|N,s)$  and   $\mathcal{\overline{B}}^{\, a}_=(\alpha_1,\ldots,\alpha_\lambda;\eta,k, r|N,s)$. Moreover, for $\tau \in \overline{\mathcal{{B}}}^{\,a}_<(\alpha_1,\ldots,\alpha_\lambda;\eta,k, r|N,\break$
$ s)$, let $\sigma=I^{a}_{s}(\tau)$, we have $|\sigma|=|\tau|+s\eta+a$.
 \end{thm}

The proof of the above theorem consists of three  parts.   Lemma \ref{iadd} shows that the $(k-1)$-insertion   is   a map from $\mathcal{\overline{B}}^{\, a}_<(\alpha_1,\ldots,\alpha_\lambda;\eta,k, r|N,s)$ to $\mathcal{\overline{B}}^{\, a}_=(\alpha_1,\ldots,\alpha_\lambda;\eta,k, r|N,s)$.  Lemma \ref{isub-lem} provides a map (that is, the $(k-1)$-separation) from   $\mathcal{\overline{B}}^{\, a}_=(\alpha_1,\ldots,\alpha_\lambda;\eta,k, r|N,s)$ to  $\mathcal{\overline{B}}^{\, a}_<(\alpha_1,\ldots,\alpha_\lambda;\eta,k, r|N,s)$.  Then we will finish the proof of Theorem \ref{deltagammathmbb} by confirming  that the $(k-1)$-insertion and the $(k-1)$-separation are inverses of each other.

\begin{lem}\label{iadd} For  $s\geq N\geq 0$, let $\tau$ be an overpartition in $\overline{\mathcal{{B}}}^{\,a}_<(\alpha_1,\ldots,\alpha_\lambda;\eta,k, r|N,s)$ and let  $\sigma=I^{a}_{s}(\tau)$. Then  $\sigma$ is an overpartition in $\overline{\mathcal{{B}}}^{\,a}_{=}
(\alpha_1,\ldots,\alpha_\lambda;\eta,k, r|N,s)$. Moreover, $|\sigma|=|\tau|+s\eta+a$.

\end{lem}

\pf To prove that $\sigma$ belongs  to  $\overline{\mathcal{{B}}}^{\,a}_{=}
(\alpha_1,\ldots,\alpha_\lambda;\eta,k, r|N,s)$, we must  verify the following conditions:
\begin{itemize}
\item[(A)]  There exists  an  overlined part $\equiv a\pmod\eta$ in $\sigma$, and assume that the  largest overlined part $\equiv a\pmod\eta$  in  $\sigma$ is $\overline{t\eta+a}$;

\item[{\rm (B)}] Let ${\hat{\sigma}}$ be the  overpartition obtained  by removing
  $\overline{t\eta+a}$  from $\sigma$. Then there are $N$ parts marked with $k-1$ in  $G({\hat{\sigma}})$,  denoted $ \tilde{g}_1({\hat{\sigma}})>  \tilde{g}_2({\hat{\sigma}})>\cdots > \tilde{g}_N({\hat{\sigma}})${\rm{;}}

    \item[(C)] Let $p$  be the smallest integer such that  $ \tilde{g}_{p+1}({\hat{\sigma}})< \overline{t\eta+a}$. Then we have $p+t=s${\rm{;}}

    \item[(D)]   $f_{\leq \eta}(\sigma)\leq r-1${\rm{;}}

    \item[(E)] The marks in $G(\sigma)$ do  not exceed $k-1$.
\end{itemize}

\noindent{  Condition (A).}  Let $\tilde{r}_1(\tau)>\cdots>\tilde{r}_N(\tau)$ be  the $(k-1)$-marked parts in $RG(\tau)$. Assume that $p$ is the smallest integer such that   $0\leq p\leq N$ and
\begin{equation}\label{iadd-lem-l}
\overline{(s-p)\eta+a}\geq\tilde{r}_{p+1}(\tau)+\eta.
\end{equation}
By the choice of $p$, we find that for $p\geq 1$ and $1\leq i\leq p$,
\begin{equation}\label{iadd-lem}
 \overline{(s-i+1)\eta+a}<\tilde{r}_{i}(\tau)+\eta.
\end{equation}
Since
$\tau\in \overline{\mathcal{{B}}}^{\,a}_<(\alpha_1,\ldots,\alpha_\lambda;\eta,k, r|N,s)$,  the largest overlined part $\equiv a\pmod\eta$ in $\tau$ is less than $\overline{(s-p)\eta+a}$. By the construction of  $I^{a}_{s}$, together with \eqref{iadd-lem}, we deduce that the largest overlined part $\equiv a\pmod\eta$ in $\sigma$ is  $\overline{(s-p)\eta+a}$, that is, $t=s-p$.

\noindent{ Condition (B).} Since ${\hat{\sigma}}$ is the  overpartition obtained  by removing
  $\overline{t\eta+a}$  from $\sigma$, by the construction of $I^{a}_{s}$, we find that ${\hat{\sigma}}=\phi_p(\tau)$. In view of  Proposition \ref{lem-for1}, we know that there are $N$ parts marked with $k-1$ in  $G({\hat{\sigma}})$, denoted  $\tilde{g}_1({\hat{\sigma}})>  \tilde{g}_2({\hat{\sigma}})>\cdots > \tilde{g}_N({\hat{\sigma}})$.

\noindent{ Condition (C).} From the proof  for the condition (A), we observe that  the largest overlined part $\equiv a\pmod\eta$ in $\sigma$ is  $\overline{(s-p)\eta+a}$, that is, $t=s-p$. We attempt  to show that $p$ is the smallest integer such that $\tilde{g}_{p+1}({\hat{\sigma}})<\overline{t\eta+a}$. By Proposition \ref{lem-for1}, we get
  \begin{equation*}\label{relation-ab-3}
 \tilde{g}_i(\hat{\sigma})= \tilde{r}_i(\tau)+\eta  \text{ for }1\leq i\leq p, \ \text{ and }\  \tilde{r}_{i}(\tau)\leq \tilde{g}_i({\hat{\sigma}})\leq  \tilde{r}_{i,1}(\tau)\text{ for }p< i\leq N,
\end{equation*}
where $\tilde{r}_{i,1}(\tau)\geq  \cdots \geq \tilde{r}_{i,{k-2}}(\tau) \geq \tilde{r}_{i}(\tau)$ are the parts in  the $(k-1)$-band of $\tau$ induced by $\tilde{r}_i(\tau)$.
It follows that for $1\leq i\leq p$,
\begin{equation}\label{relation-ab-3ttt}
    \tilde{g}_i({\hat{\sigma}})= \tilde{r}_i(\tau)+\eta> \overline{(s-i+1)\eta+a},
\end{equation}
and
\[\tilde{g}_{p+1}({\hat{\sigma}})\leq \tilde{r}_{p+1,1}(\tau)\leq \tilde{r}_{p+1}(\tau)+\eta,\]
with strict inequality if $ \tilde{r}_{p+1}(\tau)$ is overlined. Consequently, in view of \eqref{iadd-lem-l}, we deduce that
$\tilde{g}_{p+1}({\hat{\sigma}})<\overline{t\eta+a}.$ But, by \eqref{relation-ab-3ttt}, we find that $\tilde{g}_i({\hat{\sigma}})>\overline{t\eta+a}$  for $1\leq i\leq p$, from which it follows that $p$ is the smallest integer such that $\overline{t\eta+a}>\tilde{g}_{p+1}({\hat{\sigma}})$.

\noindent{Condition (D).} By the construction of $\hat{\sigma}$, we know that $f_{\leq \eta}(\hat{\sigma})\leq r-1$.
To show that $f_{\leq \eta}(\sigma)\leq r-1$, we consider the following two cases:

Case 1: If $\overline{(s-p)\eta+a}>\eta$, then   $f_{\leq \eta}(\sigma)=f_{\leq \eta}(\hat{\sigma})\leq r-1$.

Case 2: If $\overline{(s-p)\eta+a}\leq\eta$, then   $p=s$ and $a\neq \eta$.  Moreover, because of the choice of $p$,  we further have   $s=p=N$. We now encounter two subcases.

Subcase 2.1: If $f_{\leq \eta}(\tau)<r-1$, then $f_{\leq \eta}(\sigma)=f_{\leq \eta}(\hat{\sigma})+1\leq f_{\leq \eta}(\tau)+1\leq r-1$.

Subcase 2.2: If $f_{\leq \eta}(\tau)=r-1$,  then  $N\geq1$. Based on the fact that  $s=p=N\geq 1$ and the condition (3) in the definition of $ \overline{\mathcal{{B}}}^{\,a}_<(\alpha_1,\ldots,\alpha_\lambda;\eta,k, r|N,s)$, we see that   $\tilde{r}_N(\tau)\leq \eta$.  In the case $p=N$,   we may employ the forward move   to add $\eta$ to $\tilde{r}_N(\tau)$ in $\tau$ to get  $f_{\leq \eta}(\hat{\sigma})=f_{\leq \eta}(\tau)-1$.
Hence   $f_{\leq \eta}(\sigma)=f_{\leq \eta}(\hat{\sigma})+1= f_{\leq \eta}(\tau)= r-1$.

\noindent{Condition (E).} Recall that  the marks in $G(\hat{\sigma})$ do  not exceed $k-1$ and $\sigma$ is obtained   by inserting  $\overline{(s-p)\eta+a}$ into ${\hat{\sigma}}$. To   show that   the marks in $G(\sigma)$ do  not exceed $k-1$, it is enough to prove that    there are no  $(k-1)$-bands  of ${\hat{\sigma}}$ in $(\overline{(s-p-1)\eta+a}, \overline{(s-p+1)\eta+a})$. Suppose to the contrary  that there exists a $(k-1)$-band $\{{\hat{\sigma}}_{m+l}\}_{0\leq l\leq k-2}$  of ${\hat{\sigma}}$ in \break$(\overline{(s-p-1)\eta+a},\overline{(s-p+1)\eta+a})$, namely,
\[\overline{(s-p+1)\eta+a}> {\hat{\sigma}}_{m}\geq{\hat{\sigma}}_{m+1}\geq \cdots \geq {\hat{\sigma}}_{m+k-2}> \overline{(s-p-1)\eta+a}.\]
From  the construction of the $(k-1)$-insertion, we find that   $\{{\hat{\sigma}}_{m+l}\}_{0\leq l\leq k-2}$  is also a $(k-1)$-band of $\tau$. Hence there is a $(k-1)$-marked part $\hat{\sigma}_{m+t}$ ($0\leq t\leq k-2$) in $RG(\tau)$.

 Case 1: ${\hat{\sigma}}_{m}<\tilde{r}_{p}(\tau)$. In this case,  { $\tilde{r}_{p+1}(\tau)\geq\hat{\sigma}_{m+t}> \overline{(s-p-1)\eta+a}$, which contradicts \eqref{iadd-lem-l}.

  Case 2: ${\hat{\sigma}}_{m}>\tilde{r}_{p}(\tau)$.  Setting $i=p$ in \eqref{iadd-lem} gives  $\overline{(s-p+1)\eta+a}<\tilde{r}_{p}(\tau)+\eta$, whence  ${\hat{\sigma}}_{m}<\overline{(s-p+1)\eta+a}<\tilde{r}_{p}(\tau)+\eta$. Consequently,   $\tilde{r}_{p}(\tau)$ is a part of $\tau$  in $({\hat{\sigma}}_m-\eta, {\hat{\sigma}}_m)$, that is, ${\hat{\sigma}}_m-\eta<\tilde{r}_{p}(\tau)<{\hat{\sigma}}_m$.  Since $\{{\hat{\sigma}}_{m+l}\}_{0\leq l\leq k-2}$ is a $(k-1)$-band  of ${\hat{\sigma}}$, there are exactly $k-2$ parts of ${\hat{\sigma}}$  after ${\hat{\sigma}}_m$ belonging to $I({\hat{\sigma}}_m-\eta, {\hat{\sigma}}_m)$. Recalling that  $\tilde{r}_{p}(\tau)$ does not appear in $\hat{\sigma}$,  we infer  that   there are exactly $k-1$ parts of $\tau$  after ${\hat{\sigma}}_m$ belonging to $I({\hat{\sigma}}_m-\eta, {\hat{\sigma}}_m)$. This implies that there is one  part belonging to $I({\hat{\sigma}}_m-\eta, {\hat{\sigma}}_m)$   marked with $k$ in $RG(\tau)$, which is again a contradiction since  the  marks  in $RG(\tau)$ are supposed not to exceed $k-1$.

 Therefore,  there are no  $(k-1)$-bands  of ${\hat{\sigma}}$  in $(\overline{(s-p-1)\eta+a},\overline{(s-p+1)\eta+a})$.
That is to say, the marks in $G(\sigma)$ do  not exceed $k-1$ after inserting $\overline{(s-p)\eta+a}$ into ${\hat{\sigma}}$, and so the condition (E) is verified.

Thus, we have shown that  $\sigma$ is an overpartition in $\overline{\mathcal{B}}_1(\alpha_1,\ldots,\alpha_\lambda;\eta,k,r)$. Clearly,
$|\sigma|=|\tau|+s\eta+a$.  This completes the proof. \qed

We now define the $(k-1)$-separation, which plays the role of the inverse map of the $(k-1)$-insertion.

\begin{defi}[The $(k-1)$-separation]\label{isub} For  $s\geq N\geq 0$,
let $\sigma$ be an overpartition in $\overline{\mathcal{{B}}}^{\,a}_{=}(\alpha_1,\ldots,\alpha_\lambda;\eta,k, r|N,s)$ with   the largest overlined part $\equiv a\pmod \eta$ being $\overline{t\eta+a}$.   The $(k-1)$-separation  ${J}^{\, a}_{s}\colon \sigma \rightarrow \tau$ is defined as follows{\rm :} First remove  $\overline{t\eta+a}$ from $\sigma$ to produce ${\hat{\sigma}}$, and then apply the backward move $\psi_{s-t}$ to ${\hat{\sigma}}$  to obtain $\tau$.
\end{defi}

  The following lemma states that the $(k-1)$-separation has the specified image set.

\begin{lem}\label{isub-lem} For  $s\geq N\geq 0$,
let $\sigma$ be an overpartition in $\overline{\mathcal{{B}}}^{\,a}_{=}(\alpha_1,\ldots,\alpha_\lambda;\eta,k, r|N,s)$, and let $\tau={J}^{\, a}_{s}(\sigma)$. Then  $\tau$ is an overpartition in $\overline{\mathcal{{B}}}^{\,a}_<(\alpha_1,\ldots,\alpha_\lambda;\eta,k, r|N,s)$. Moreover,  $|\tau|=|\sigma|-s\eta-a$.
\end{lem}

\pf   To prove that $\tau$ belongs  to  $\overline{\mathcal{{B}}}^{\,a}_{<}
(\alpha_1,\ldots,\alpha_\lambda;\eta,k, r|N,s)$, we need to check the following conditions:
\begin{itemize}
\item[(A)]  There are $N$ parts marked with $k-1$ in  $RG(\tau)$, denoted  $\tilde{r}_1(\tau)>  \tilde{r}_2(\tau)>\cdots > \tilde{r}_N(\tau)$;

\item[(B)] Assume that $p$ is the smallest integer such that $\tilde{r}_{p+1}(\tau)+\eta\leq \overline{(s-p)\eta+a}$. Then the largest overlined part $\equiv a\pmod\eta$ in $\tau$ is less than $\overline{(s-p)\eta+a}${\rm{;}}

\item[(C)]  $f_{\leq \eta}(\tau)\leq r-1$;

\item[(D)]  If $s=N=0$ and $a\neq \eta$, then $f_{\leq \eta}(\tau)< r-1$;

\item[(E)] If  $f_{\leq \eta}(\tau)=r-1$, $s=N\geq 1$ and  $a\neq\eta$, then $\tilde{r}_N(\tau)\leq\eta$.

\end{itemize}

\noindent{ Condition (A).}  Assume that the largest overlined part $\equiv a \pmod{\eta}$ in $\sigma$ is $\overline{t\eta+a}$. Let ${\hat{\sigma}}$ be the  overpartition obtained  by removing $\overline{t\eta+a}$  from $\sigma$. By definition,   there are $N$ parts marked with $k-1$ in  $G({\hat{\sigma}})$, denoted $\tilde{g}_1({\hat{\sigma}})>\tilde{g}_2({\hat{\sigma}})>\cdots>
\tilde{g}_N({\hat{\sigma}})$. Assume that $p$ is the smallest integer such that $\overline{t\eta+a}>\tilde{g}_{p+1}({\hat{\sigma}})$. Since $\sigma$ is an overpartition in $\overline{\mathcal{{B}}}^{\,a}_=(\alpha_1,\ldots,\alpha_\lambda;\eta,k, r|N,s)$,   we have $p+t=s$.   To show that  the condition (A) holds, in view of Proposition  \ref{lem-bac1},  it suffices to verify the following statements:
\begin{itemize}
    \item[(A1)] $\tilde{g}_{p}({\hat{\sigma}})\geq\overline{\eta+\alpha_1}${\rm{;}}
    \item[(A2)]  $\tilde{g}_p({\hat{\sigma}})$ is a part in any $(k-1)$-band   of ${\hat{\sigma}}$  belonging to $I(\tilde{g}_p({\hat{\sigma}})-2\eta, \tilde{g}_p({\hat{\sigma}}))$.
\end{itemize}

\noindent{Condition (A1).}  Since  $\overline{t\eta+a}$ does not appear in ${\hat{\sigma}}$, the minimality  of $p$ implies that  $\overline{t\eta+a}< \tilde{g}_{p}({\hat{\sigma}})$.
Under the condition that the marks in $G(\sigma)$ do not exceed $k-1$,  it is obvious that there are no  $(k-1)$-bands  of ${\hat{\sigma}}$ in $(\overline{(t-1)\eta+a},\overline{(t+1)\eta+a})$. Let
\[\tilde{g}_{p}(\hat{\sigma})\geq \tilde{g}_{p,2}(\hat{\sigma})\geq \cdots\geq \tilde{g}_{p,k-2}(\hat{\sigma})\]
be the parts in the $(k-1)$-band    of $\hat{\sigma}$ induced by $\tilde{g}_{p}(\hat{\sigma})$. Then
\[\tilde{g}_{p,k-2}(\hat{\sigma})\geq \tilde{g}_{p}(\hat{\sigma})-\eta >\overline{(t-1)\eta+a}.\]
Consequently,
$\tilde{g}_{p}({\hat{\sigma}})\geq \overline{(t+1)\eta+a}$ since  there are no  $(k-1)$-bands  of ${\hat{\sigma}}$ in $(\overline{(t-1)\eta+a},\break$$\overline{(t+1)\eta+a})$. Using the fact that the largest overlined  part $\equiv a \pmod{\eta}$ in $\hat{\sigma}$ is less than $\overline{t\eta+a}$,  we are led to the strict inequality
\begin{equation}\label{condition(a)-1}
\tilde{g}_{p}({\hat{\sigma}})> \overline{(t+1)\eta+a},
\end{equation}
which yields (A1).

\noindent{Condition (A2).}
Suppose to the contrary that there is  a $(k-1)$-band of $\hat{\sigma}$ belonging to $I(\tilde{g}_p({\hat{\sigma}})-2\eta, \tilde{g}_p({\hat{\sigma}}))$ that does not contain $\tilde{g}_p({\hat{\sigma}})$ as a part, and let
\[\hat{\sigma}_{m}\geq \cdots\geq {\hat{\sigma}}_{m+k-2}\]
be the parts in this $(k-1)$-band.  We  assume that ${\hat{\sigma}}_{m+l}$ ($0\leq l\leq k-2$)  is  a part in this $(k-1)$-band marked with $k-1$ in $G({\hat{\sigma}})$. Evidently,  ${\hat{\sigma}}_{m+l}\leq {\hat{\sigma}}_{m}<\tilde{g}_p({\hat{\sigma}})$.
We   claim that ${\hat{\sigma}}_{m+l}>\tilde{g}_{p+1}({\hat{\sigma}})$. According to \eqref{condition(a)-1}, we get
 \begin{equation}\label{relation-ab-3stta}
 \overline{(t-1)\eta+a}< \tilde{g}_p({\hat{\sigma}})-2\eta\leq {\hat{\sigma}}_{m+k-2}\leq \cdots\leq {\hat{\sigma}}_{m}.
 \end{equation}
As mentioned before, there are no  $(k-1)$-bands  of ${\hat{\sigma}}$ in $(\overline{(t-1)\eta+a},\overline{(t+1)\eta+a})$, and so  it follows
from \eqref{relation-ab-3stta} that ${\hat{\sigma}}_{m}\geq \overline{(t+1)\eta+a}$. In fact,  we attain the strict inequality ${\hat{\sigma}}_{m}> \overline{(t+1)\eta+a}$ owing to the fact that the largest overlined part $\equiv a\pmod\eta$ in ${\hat{\sigma}}$ is less than $\overline{t\eta+a}$.  The assumption  that $\{\hat{\sigma}_{m+l}\}_{0\leq l\leq k-2}$ is a $(k-1)$-band ensures that ${\hat{\sigma}}_{m+k-2}\geq {\hat{\sigma}}_{m}-\eta$. Noting that ${\hat{\sigma}}_{m}> \overline{(t+1)\eta+a}$, we obtain that ${\hat{\sigma}}_{m+k-2}> \overline{t\eta+a}$. But  $\tilde{g}_{p+1}({\hat{\sigma}})< \overline{t\eta+a}$, we arrive at
\[{\hat{\sigma}}_{m+l}\geq {\hat{\sigma}}_{m+k-2}>\overline{t\eta+a}> \tilde{g}_{p+1}({\hat{\sigma}}),\]
 as claimed. Thus, we conclude that $\tilde{g}_{p+1}({\hat{\sigma}})<{\hat{\sigma}}_{m+l}<\tilde{g}_p({\hat{\sigma}})$. However, ${\hat{\sigma}}_{m+l}$ is marked with $k-1$ in $G(\hat{\sigma})$, which leads to a contradiction since   there are no  $(k-1)$-marked parts in  $G({\hat{\sigma}})$ between $\tilde{g}_{p}({\hat{\sigma}})$ and $\tilde{g}_{p+1}({\hat{\sigma}})$. This confirms the condition (A2).

With the conditions (A1) and (A2) in hand,  Proposition \ref{lem-bac1} guarantees  that  there are $N$ parts marked with $k-1$ in  $RG(\tau)$.  In addition, it gives that
  \begin{equation}\label{relation-ab-3s}
 \tilde{r}_i(\tau)= \tilde{g}_i({\hat{\sigma}})-\eta  \text{ for }1\leq i\leq p, \ \text{ and }\  \tilde{g}_{i,k-1}({\hat{\sigma}})\leq \tilde{r}_i(\tau)\leq  \tilde{g}_{i}({\hat{\sigma}})\text{ for }p< i\leq N.
\end{equation}
Thus, we have proved that $\tau$ satisfies the condition (A).

\noindent{  Condition (B).}  We aim to show that $p$ is also the smallest integer   such that $\tilde{r}_{p+1}(\tau)+\eta\leq \overline{(s-p)\eta+a}$, where   $p$ is defined to be the smallest integer such that $\tilde{g}_{p+1}({\hat{\sigma}})<\overline{t\eta+a}$. Applying   \eqref{relation-ab-3s} with $i=p+1$ yields that $\tilde{r}_{p+1}(\tau)\leq \tilde{g}_{p+1}({\hat{\sigma}})< \overline{t\eta+a}$. Let
\[\tilde{r}_{p+1,1}(\tau)\geq \cdots\geq \tilde{r}_{p+1,k-2}(\tau)\geq \tilde{r}_{p+1}(\tau)\]
be the parts in the $(k-1)$-band   of $\tau$ induced by $\tilde{r}_{p+1}(\tau)$. Then
\[\tilde{r}_{p+1,1}(\tau)\leq \tilde{r}_{p+1}(\tau)+\eta< \overline{(t+1)\eta+a}.\]
Since $\{\tilde{r}_{p+1}(\tau)\}_{k-1}$ is also a  $(k-1)$-band  of ${\hat{\sigma}}$ and  there are no  $(k-1)$-bands  of ${\hat{\sigma}}$ in $(\overline{(t-1)\eta+a},\overline{(t+1)\eta+a})$, we deduce that
    \begin{equation}\label{relation-ab-3bb}
 \tilde{r}_{p+1}(\tau)\leq \overline{(t-1)\eta+a}=\overline{(s-p-1)\eta+a}.
 \end{equation}
Combining \eqref{condition(a)-1} and \eqref{relation-ab-3s}, we find that
  \begin{equation*}
     \tilde{r}_p(\tau)= \tilde{g}_{p}({\hat{\sigma}})-\eta > \overline{t\eta+a}=\overline{(s-p)\eta+a}.
  \end{equation*}
Hence  for $1\leq i<p$,
\begin{equation}\label{relation-ab-3cc}
    \tilde{r}_i(\tau)\geq \tilde{r}_{i+1}(\tau)+\eta\geq\cdots\geq \tilde{r}_{p}(\tau)+(p-i)\eta>\overline{(s-i)\eta+a}.
\end{equation}
By inspection of \eqref{relation-ab-3bb}  and \eqref{relation-ab-3cc},   we conclude that $p$ is the smallest integer such that $\tilde{r}_{p+1}(\tau)+\eta\leq \overline{(s-p)\eta+a}$. On the other hand, by the definition of $J_s^a$, we obtain that  the largest overlined part $\equiv a\pmod\eta$ in $\tau$ is less than  $\overline{(s-p)\eta+a}$, and so the condition (B) is justified.

\noindent{  Condition (C).} To show that  $f_{\leq\eta}(\tau)\leq r-1$, we consider  three cases:

Case c1:  $p=0$. In this case, $\tau=\hat{\sigma}$,  so  $f_{\leq\eta}(\tau)\leq f_{\leq\eta}(\hat{\sigma})\leq r-1$.

Case c2:    $p\geq 1$ and $\tilde{g}_p({\hat{\sigma}})-\eta>\eta$. In this case,  $f_{\leq\eta}(\tau)=f_{\leq\eta}(\hat{\sigma})\leq f_{\leq\eta}(\sigma)\leq r-1$.

Case c3:  $p\geq 1$  and $\tilde{g}_p({\hat{\sigma}})-\eta\leq\eta$. In this case, $f_{\leq\eta}(\tau)=f_{\leq\eta}(\hat{\sigma})+1$. It follows from \eqref{condition(a)-1} that  $\eta\geq  \tilde{g}_{p}({\hat{\sigma}})-\eta>\overline{t\eta+a}$, and so  $f_{\leq\eta}(\hat{\sigma})=f_{\leq\eta}(\sigma)-1$. Hence $f_{\leq\eta}(\tau)=f_{\leq\eta}(\hat{\sigma})+1=f_{\leq\eta}(\sigma)\leq r-1$.

\noindent{ Condition (D).}  If $s=N=0$ and $a\neq \eta$, then    $\tau$ is obtained by removing $\overline{a}$ from $\sigma$. This implies that $f_{\leq\eta}(\tau)=f_{\leq\eta}(\sigma)-1<r-1$.

\noindent{  Condition (E).}  There are two cases.

Case e1: If $\tilde{g}_p(\hat{\sigma})-\eta\leq\eta$, then $p=N$ and  $\tilde{r}_N(\tau)=\tilde{g}_N(\hat{\sigma})-\eta\leq\eta$.

Case e2: If $\tilde{g}_p(\hat{\sigma})-\eta>\eta$, then  $f_{\leq\eta}(\hat{\sigma})=f_{\leq\eta}(\tau)$. The condition     $f_{\leq\eta}(\tau)=r-1$ implies that $f_{\leq\eta}(\hat{\sigma})=r-1$. We claim that  $p<N$ in this case. Suppose to the contrary that $p=N$.  If so, we have $t=s-p=N-p=0$. This implies that $\hat{\sigma}$ is obtained by removing $\overline{a}$ from $\sigma$, and thus   $f_{\leq\eta}(\hat{\sigma})=f_{\leq\eta}(\sigma)-1<r-1$, which contradicts  $f_{\leq\eta}(\hat{\sigma})=r-1$. Hence we have $p<N$. In light of \eqref{relation-ab-3bb}, we obtain that
\[\tilde{r}_N(\pi)\leq \tilde{r}_{N-1}(\pi)-\eta\leq\cdots\leq  \tilde{r}_{p+1}(\pi)-(N-p-1)\eta\leq \overline{(s-N)\eta+a}= \overline{a}<\eta.\]
Therefore, we have proved that the condition (E) is fulfilled.

So far we   have accomplished the task of showing that  $\tau$ is an overpartition in $\overline{\mathcal{{B}}}^{\,a}_<(\alpha_1,\ldots,\break$$\alpha_\lambda;\eta,k,r|N,s)$. Evidently,  $|\tau|=|\sigma|-s\eta-a$.  This completes the proof. \qed

 We are now ready to give a proof of Theorem \ref{deltagammathmbb} based on Lemma \ref{iadd} and Lemma \ref{isub-lem}.

{\noindent \it Proof of Theorem \ref{deltagammathmbb}.}  Let $\tau \in\overline{\mathcal{{B}}}^{\,a}_<(\alpha_1,\ldots,\alpha_\lambda;\eta,k, r|N,s)$. Utilizing Lemma \ref{iadd}, we find that $I^{a}_s(\tau)$ belongs to $\overline{\mathcal{{B}}}^{\,a}_=(\alpha_1,\ldots,\alpha_\lambda;\eta,k, r|N,s)$. Appealing to the condition (C) in the proof of Lemma \ref{iadd} and the property (4) in Proposition \ref{lem-for1}, we deduce that ${J}^{\, a}_s(I^{a}_s(\tau))=\tau$.

Conversely, let $\gamma \in \overline{\mathcal{{B}}}^{\,a}_=(\alpha_1,\ldots,\alpha_\lambda;\eta,k, r|N,s)$. Invoking  Lemma \ref{isub-lem}, we know that   ${J}^{\, a}_s(\gamma) \in \overline{\mathcal{{B}}}^{\,a}_<(\alpha_1,\ldots,\alpha_\lambda;\eta,k, r|N,s).$ By virtue of the condition (B) in the proof of Lemma \ref{isub-lem} and the property (4) in Proposition \ref{lem-bac1}, we obtain that  $I^{a}_s({J}^{\, a}_s(\gamma))=\gamma$.

Therefore,   the map $I^{a}_s$ is a bijection between $\overline{\mathcal{{B}}}^{\,a}_<(\alpha_1,\ldots,\alpha_\lambda;\eta,k, r|N,s)$ and $\overline{\mathcal{{B}}}^{\,a}_=(\alpha_1,\ldots,\break$
$\alpha_\lambda;\eta,k, r|N,s)$. This completes the proof. \qed

The following theorem gives a   criterion   to determine whether an overpartition in $\overline{\mathcal{{B}}}^{\,a}_{=}(\alpha_1,\ldots,\alpha_\lambda;\eta,k, r|N,s)$ is also an overpartition in $\overline{\mathcal{{B}}}^{\,a}_{<}(\alpha_1,\ldots,\alpha_\lambda;\eta,k, r|N',s')$, which involves  the successive application of the $(k-1)$-insertion operations.

\begin{thm}\label{ssins}
For $s\geq N\geq 0$, let $\sigma$ be an overpartition in $\overline{\mathcal{{B}}}^{\,a}_{=}(\alpha_1,\ldots,\alpha_\lambda;\eta,k, r|N,s)$. Assume that there are $N'$ parts marked with $k-1$ in the reverse Gordon marking of $\sigma$. Then  $\sigma$ is also an overpartition in $\overline{\mathcal{{B}}}^{\,a}_<(\alpha_1,\ldots,\alpha_\lambda;\eta,k, r|N',s')$ if and only if $s'>s$.

\end{thm}
\pf   We  first show that if $s'>s$, then $\sigma$ is  in $\overline{\mathcal{{B}}}^{\,a}_<(\alpha_1,\ldots,\alpha_\lambda;\eta,k, r|N',s')$.  Let $\tilde{r}_1({{\sigma}})>\cdots>\tilde{r}_{N'}({{\sigma}})$ be the   $(k-1)$-marked parts in $RG(\sigma)$. We are required to prove that $\sigma$ satisfies the following conditions:
\begin{itemize}
    \item[(A)]  If $p'$ is the smallest integer such that
\begin{equation*}
 \overline{(s'-p')\eta+a}\geq\tilde{r}_{p'+1}(\sigma)+\eta,
\end{equation*}
then the largest overlined part $\equiv a\pmod\eta$ in $\sigma$ is  less than $\overline{(s'-p')\eta+a}$;
 \item[(B)] If $f_{\leq\eta}(\sigma)=r-1$, $s'=N'\geq 1$ and $a\neq \eta$, then $\tilde{r}_{N'}(\sigma)\leq \eta$.

\end{itemize}

\noindent{Condition (A).} Assume that $\overline{t\eta+a}$ is the largest overlined part $\equiv a\pmod \eta$ in $\sigma$. Let ${\hat{\sigma}}$  be the   overpartition   obtained  from $\sigma$ by removing $ {\overline{t\eta+a}}$.  By definition, there are  $N$ parts marked with $k-1$ in $G({\hat{\sigma}})$, denoted $\tilde{g}_1({\hat{\sigma}})>\cdots>
\tilde{g}_N({\hat{\sigma}})$. Let  $p$ be the smallest integer such that $ \tilde{g}_{p+1}({\hat{\sigma}})<\overline{t\eta+a}$. Since $\sigma \in \overline{\mathcal{{B}}}^{\,a}_{=}(\alpha_1,\ldots,\alpha_\lambda;\eta,k, r|N,s)$, we have  $p=s-t$. Using Proposition \ref{sequence-length}, we find that there are also $N$ parts marked with $k-1$ in $RG({\hat{\sigma}})$, denoted  $\tilde{r}_1({\hat{\sigma}})>\cdots>\tilde{r}_{N}({\hat{\sigma}})$ and  that $\tilde{r}_{i}({\hat{\sigma}})\leq\tilde{g}_{i}({\hat{\sigma}})$ for $1\leq i\leq N$.
In particular, $\tilde{r}_{p+1}({\hat{\sigma}})\leq\tilde{g}_{p+1}({\hat{\sigma}})$. But $ \tilde{g}_{p+1}({\hat{\sigma}})<\overline{t\eta+a}$, so we get  $\tilde{r}_{p+1}({\hat{\sigma}})<\overline{t\eta+a}$.

 We now attempt to show that $\tilde{r}_{p+1}(\sigma)\leq \overline{t\eta+a}.$  Suppose to the contrary that $\tilde{r}_{p+1}(\sigma)> \overline{t\eta+a}$. Since ${\hat{\sigma}}$  is the   overpartition   obtained  from $\sigma$ by removing $ {\overline{t\eta+a}}$, we find that $\tilde{r}_{p+1}({\hat{\sigma}})=\tilde{r}_{p+1}(\sigma)$,  which implies  $\tilde{r}_{p+1}({\hat{\sigma}})>\overline{t\eta+a}$, contradicting  the preceding
 assertion that  $\tilde{r}_{p+1}({\hat{\sigma}})< \overline{t\eta+a}$. This proves  $\tilde{r}_{p+1}(\sigma)\leq \overline{t\eta+a}.$

Examining  the construction of  ${\hat{\sigma}}$, we notice that $N'$  equals either $N$  or $N+1$. Under the condition that $s'>s$, we get $s'\geq s+1\geq N+1\geq N'$.  Let $p'$ be the smallest integer such that
\[ \overline{(s'-p')\eta+a}\geq\tilde{r}_{p'+1}(\sigma)+\eta.\]
Since $s'-p>s-p=t$ and   $\tilde{r}_{p+1}(\sigma)\leq \overline{t\eta+a}$,   we find that
  \[\overline{(s'-p)\eta+a}\geq \overline{(t+1)\eta+a}=\overline{t\eta+a}+\eta\geq \tilde{r}_{p+1}(\sigma)+\eta.\]
Hence the choice of $p'$ implies that
  \[p'\leq p\leq N\leq N'.\]
This leads to
  \[ \overline{(s'-p')\eta+a}> \overline{(s-p)\eta+a}= \overline{t\eta+a}.\]
 This proves the condition (A)  because   $\overline{t\eta+a}$ is the largest overlined part $\equiv a\pmod\eta$ of $\sigma$.

\noindent{Condition (B).}  As we know,  $\tilde{r}_{p+1}(\sigma)\leq \overline{t\eta+a}$, so that
\[ \overline{t\eta+a}\geq\tilde{r}_{p+1}(\sigma)\geq\tilde{r}_{p+2}(\sigma)+\eta\geq\cdots\geq \tilde{r}_{N+1}(\sigma)+(N-p)\eta.\]
Observing that $s'\geq s+1 \geq N+1\geq N'$, we find that $N'=N+1\geq 1$ and $s=N$ when  $s'=N'$.  It follows that $\tilde{r}_{N'}(\sigma)=\tilde{r}_{N+1}(\sigma)\leq \overline{t\eta+a}-(N-p)\eta=\overline{a}< \eta$ because  $p+t=s=N$. Thus, we have proved the condition (B) is valid.

This completes the proof of the sufficiency. Conversely, assume that $\sigma$ is  in  both  $\overline{\mathcal{{B}}}^{\,a}_=(\alpha_1,\ldots,\alpha_\lambda;\eta,k, r|N,s)$ and    $\overline{\mathcal{{B}}}^{\,a}_<(\alpha_1,\ldots,\alpha_\lambda;\eta,k, r|N',s')$, we intend to show that  $s'>s$.

Given that $ \sigma$  belongs to $\overline{\mathcal{{B}}}^{\,a}_{=}(\alpha_1,\ldots,\alpha_\lambda;\eta,k, r|N,s)$, we may assume that  $\overline{t\eta+a}$ is the largest overlined part $\equiv a\pmod\eta$ in $ \sigma$. Let  $\hat{ \sigma}$ be the   overpartition  obtained  from $ \sigma$ by removing $\overline{t\eta+a}$. Then there are $N$ parts marked with $k-1$ in $G(\hat{ \sigma})$, denote $\tilde{g}_1(\hat{ \sigma})>\cdots>\tilde{g}_{N}(\hat{ \sigma})$. Let $p$ be the smallest integer such that $\tilde{g}_{p+1}(\hat{ \sigma})<\overline{t\eta+a}$. Since $ \sigma \in \overline{\mathcal{{B}}}^{\,a}_{=}(\alpha_1,\ldots,\alpha_\lambda;\eta,k, r|N,s)$, we have $p=s-t$.    By the reasoning in the proof of Lemma \ref{isub-lem}, we  establish that
  \begin{equation}\label{lemma4.8aa}
 \tilde{g}_{p}(\hat{ \sigma})>\overline{(t+1)\eta+a}.
 \end{equation}

On the other hand,  since $ \sigma$ is also  in ${\mathcal{{B}}}^{\,a}_<(\alpha_1,\ldots,\alpha_\lambda;\eta,k, r|N',s')$,   there are $N'$ parts marked with $k-1$ in $RG( \sigma)$, denote $\tilde{r}_1( \sigma)>\cdots>
\tilde{r}_{N'}( \sigma)$. Assume that $p'$ is the smallest integer such that    \begin{equation}\label{thm4-10a}
\overline{(s'-p')\eta+a}\geq\tilde{r}_{p'+1}( \sigma)+\eta.
\end{equation}
The condition that $ \sigma \in {\mathcal{{B}}}^{\,a}_<(\alpha_1,\ldots,\alpha_\lambda;\eta,k, r|N',s')$ ensures that the largest overlined part $\equiv a\pmod \eta$ in $ \sigma$ is less than $\overline{(s'-p')\eta+a}$. But the largest overlined part $\equiv a\pmod \eta$ in $ \sigma$ is supposed  to be $\overline{t\eta+a}$, so we get  $s'-p'>t=s-p$.

Our final goal is to show that   $s'>s$.  Suppose to the contrary that $s'\leq s$. This implies that  $p'<p$ since $s'-p'>s-p$. Let $\tilde{r}_1(\hat{ \sigma})>\cdots>\tilde{r}_{N}(\hat{ \sigma})$ be the $(k-1)$-marked parts in $RG(\hat{ \sigma})$. In view of Proposition \ref{sequence-length}, we find that $\tilde{r}_{p}(\hat{ \sigma})\geq \tilde{g}_{p}(\hat{ \sigma})-\eta$. Comparison with \eqref{lemma4.8aa} yields
 \[
 \tilde{r}_{p}(\hat{ \sigma})\geq \tilde{g}_{p}(\hat{ \sigma})-\eta>\overline{t\eta+a},
 \]
so that
 \begin{equation}\label{lemma4.8aaa}
   \tilde{r}_{p}( \sigma)=\tilde{r}_{p}(\hat{ \sigma})>\overline{t\eta+a}.
 \end{equation}
Using the fact that $p'<p$, we obtain that
  \begin{equation}\label{lemma4.8aaab}
    \tilde{r}_{p'+1}( \sigma)\geq \tilde{r}_{p'+2}( \sigma)+\eta\geq \cdots\geq \tilde{r}_{p}( \sigma)+(p-p'-1)\eta.
  \end{equation}
  Substituting \eqref{lemma4.8aaa} into \eqref{lemma4.8aaab}, we arrive at
    \[\tilde{r}_{p'+1}( \sigma)>\overline{t\eta+a}+(p-p'-1)\eta=\overline{(s-p')\eta+a}-\eta\geq \overline{(s'-p')\eta+a}-\eta,\]
that is,
\[\tilde{r}_{p'+1}( \sigma)+\eta> \overline{(s'-p')\eta+a},\]
which is  in contradiction to  \eqref{thm4-10a}.
Thus, we have shown $s'>s$. This completes the proof.   \qed

\subsection{Proof of Theorem \ref{lem-b-0-over}}

In this subsection, we will give a proof of Theorem \ref{lem-b-0-over}  by successively  applying
the forward move  and  the $(k-1)$-insertion with $a=\eta$.

{\it \noindent Proof of Theorem \ref{lem-b-0-over}.} Let $\mu$ be an overpartition in $ \mathcal{{B}}_0(\alpha_1,\ldots,\alpha_\lambda;\eta,k,r)$ and let $\zeta$ be a partition with distinct  parts divisible by $\eta$. We wish to construct an overpartition $\pi=\Phi(\zeta,\mu)$ in $\overline{\mathcal{{B}}}_1(\alpha_1,\ldots,\alpha_\lambda;\eta,k,r)$ such that $|\pi|=|\zeta|+|\mu|$. There are two cases:

\noindent Case 1:  $\zeta=\emptyset$. Set $\pi=\mu$. Obviously, $\pi \in  \overline{\mathcal{{B}}}_1(\alpha_1,\ldots,\alpha_\lambda;\eta,k,r)$ and $|\pi|=|\zeta|+|\mu|$.

\noindent Case 2:  $\zeta\neq\emptyset$. Assume that there  are $N$ parts marked with $k-1$ in $G(\mu)$,  and set $\zeta=(\eta\zeta_1,\ldots,\eta\zeta_{c},\eta\zeta_{c+1}\ldots,
\eta\zeta_{c+m})$, where $\zeta_1>\cdots>\zeta_c>N\geq \zeta_{c+1}>\cdots>\zeta_{c+m}>0.$  We first  merge $\eta\zeta_{c+1},\ldots,\eta\zeta_{c+m}$ and $\mu$ by successively applying the forward move. Then, we will merge  $\eta\zeta_1,\ldots,\eta\zeta_c$ and $\mu$ by applying the $(k-1)$-insertion with $a=\eta$  to generate $c$ overlined parts  divisible by $\eta$.

\noindent{ Step 1.} Let $ \tilde{g}_1(\mu)> \tilde{g}_2(\mu)>\cdots >\tilde{g}_N(\mu)$ be the $(k-1)$-marked parts in  $G(\mu)$. Note that  $\mu\in\mathcal{{B}}_0(\alpha_1,\ldots,\alpha_\lambda;\eta,k,r)$, we see that  $\{\tilde{g}_{i}(\mu)\}_{k-1}$ are even for $1\leq i\leq N$. We first merge $\eta\zeta_{c+1},\ldots,\eta\zeta_{c+m}$ into $\mu$ by successively  applying the forward move. Denote the intermediate overpartitions by $\mu^{(0)},\mu^{(1)},\ldots,\mu^{(m)}$ with $\mu^{(0)}=\mu$.

 Since $\zeta_{c+1}\leq N$, we find that   $\mu$ is  an overpartition in ${\mathcal{B}}_e(\alpha_1,\ldots,\alpha_\lambda;\eta,k,r|N,
 \zeta_{c+1})$.  Set $b=0$ and repeat the following procedure until $b=m$:

\begin{itemize}

\item[(A)]  Merge $\eta \zeta_{c+b+1}$ into $\mu^{(b)}$.  Apply the forward move $\phi_{\zeta_{c+b+1}}$   to $\mu^{(b)}$ to obtain  $\mu^{(b+1)}$,  that is,     \[\mu^{(b+1)}=\phi_{\zeta_{c+b+1}}(\mu^{(b)}).\]
Since
\[\mu^{(b)}\in {\mathcal{B}}_e(\alpha_1,\ldots,\alpha_\lambda;\eta,k,r|N,
   \zeta_{c+b+1}),\]
in view of Lemma \ref{deltagamma}, we deduce that
 \begin{equation*}
\mu^{(b+1)}\in {\mathcal{B}}_d(\alpha_1,\ldots,\alpha_\lambda;\eta,k,r|N,
 \zeta_{c+b+1}),
 \end{equation*}
 and
 \[
  |\mu^{(b+1)}|=|\mu^{(b)}|+\eta\zeta_{c+b+1}.\]

   \item[(B)] Replace $b$ by $b+1$.  If $b=m$, then we are done. If $b<m$, then  we have
\[\mu^{(b)}\in {\mathcal{B}}_e(\alpha_1,\ldots,\alpha_\lambda;\eta,k,r|N,
   \zeta_{c+b+1}),\]
   since    $\zeta_{c+b+1}<\zeta_{c+b}\leq N$. Go back to (A).

     \end{itemize}
 Eventually, the above procedure yields $\mu^{(m)}\in {\mathcal{B}}_d(\alpha_1,\ldots,\alpha_\lambda;\eta,k,r
    |N,\zeta_{c+m})$
such that
    \begin{equation} \label{wei-par1}
      |\mu^{(m)}|=|\mu^{(0)}|+\eta\zeta_{c+1}+\cdots+\eta\zeta_{c+m}.
    \end{equation}

{\noindent Step 2.} We continue to merge $\eta\zeta_c,\ldots,\eta\zeta_1$    into $\mu^{(m)}$ by successively applying the $(k-1)$-insertion with $a=\eta$. Denote the intermediate overpartitions by $\mu^{(m)},\mu^{(m+1)},\ldots,\mu^{(m+c)}$ and set $\pi=\mu^{(m+c)}$.
Assume that there are $N(\mu^{(i)})$ parts marked with $k-1$ in $RG(\mu^{(i)})$, where $m\leq i\leq m+c$ and  $N(\mu^{(m)})=N$.

Assume that  $p$ is the smallest integer such that $0\leq p\leq N$ and $\overline{(\zeta_c-p)\eta}\geq \tilde{r}_{p+1}(\mu^{(m)})+\eta$. Such an integer $p$ exists  because
  $\zeta_c>N$ and $\overline{(\zeta_c-N)\eta}>0\geq-\infty=\tilde{r}_{N+1}(\mu^{(m)})+\eta$.  Since $\mu^{(m)}\in {\mathcal{B}}_d(\alpha_1,\ldots,\alpha_\lambda;\eta,k,r
    |N,\zeta_{c+m})$,  there are no overlined parts divisible by $\eta$ in $\mu^{(m)}$. Hence the largest overlined part divisible by $\eta$ in $\mu^{(m)}$ is less than  $\overline{(\zeta_c-p)\eta}$. It follows that
\[\mu^{(m)}\in {\overline{\mathcal{B}}}^{\,\eta}_<(\alpha_1,\ldots,\alpha_\lambda;\eta,k,r|N(\mu^{(m)}),\zeta_{c}-1).\]

Merging $\eta\zeta_c,\ldots,\eta\zeta_1$    into $\mu^{(m)}$, the following procedure  generates $c$ overlined parts divisible by $\eta$. We start with setting $b=0$.
 \begin{itemize}
\item[(A)] Merge $\eta \zeta_{c-b}$ into $\mu^{(m+b)}$ to generate an overlined part divisible by $\eta$.  More precisely,  applying the $(k-1)$-insertion $I^{\,\eta}_{\zeta_{c-b}-1}$ to $\mu^{(m+b)}$, we  obtain
\[\mu^{(m+b+1)}=I^{\,\eta}_{\zeta_{c-b}-1}(\mu^{(m+b)}).\]
Since
\begin{equation*}
 \mu^{(m+b)} \in {\overline{\mathcal{B}}}^{\,\eta}_<(\alpha_1,\ldots,\alpha_\lambda;
 \eta,k,r|N(\mu^{(m+b)}),\zeta_{c-b}-1),
 \end{equation*}
in view of Lemma \ref{iadd}, we find that
 \begin{equation*}
 \mu^{(m+b+1)} \in {\overline{\mathcal{B}}}^{\,\eta}_=(\alpha_1,\ldots,\alpha_\lambda;
 \eta,k,r|N(\mu^{(m+b)}),\zeta_{c-b}-1),
 \end{equation*}
 and
 \[
 |\mu^{(m+b+1)}|=|\mu^{(m+b)}|+\eta\zeta_{c-b}.\]

\item[(B)] Replace $b$ by $b+1$.  If $b=c$, then we are done. If $b<c$, since $\zeta_{c-b}>\zeta_{c-b+1}$,  it follows from Theorem \ref{ssins}   that
 \[\mu^{(m+b)}\in {\overline{\mathcal{B}}}^{\,\eta}_<(\alpha_1,\ldots,\alpha_\lambda;\eta,k,r
 |N(\mu^{(m+b)}),\zeta_{c-b}-1).\]
Go back to (A).

 \end{itemize}
The above procedure generates an overpartition
$\pi=\mu^{(m+c)}\in \overline{\mathcal{B}}^{\,\eta}_{=} (\alpha_1,\ldots,\alpha_\lambda;\eta,k,r\break$
$|N(\mu^{(m+c-1)}),  \zeta_{1}-1)$
     such that
     \begin{equation} \label{wei-par2}
         |\mu^{(m+c)}|=|\mu^{(m)}|+\eta\zeta_{c}
         +\cdots+\eta\zeta_{1}.
    \end{equation}
From the construction of the $(k-1)$-insertion with $a=\eta$, it can be seen that $\pi$ is an overpartition in $\overline{\mathcal{B}}_1(\alpha_1,\ldots,\alpha_\lambda;\eta,k,r)$ with $c$ overlined parts divisible by $\eta$. Furthermore, combining \eqref{wei-par1} with \eqref{wei-par2}, we find that $|\pi|=|\mu|+|\zeta|$. Therefore, $\Phi$ is a desired map from $\mathcal{D}_\eta\times{\mathcal{B}}_{0}(\alpha_1,\ldots,\alpha_\lambda;\eta,k,r)$ to $\mathcal{\overline{B}}_1(\alpha_1,\ldots,\alpha_\lambda;\eta,k,r)$.

To prove that $\Phi$ is a bijection, we shall define the inverse map $\Psi$ of $\Phi$ from $\mathcal{\overline{B}}_1(\alpha_1,\ldots,\break$$\alpha_\lambda;\eta,k,r)$ to  $\mathcal{D}_\eta\times{\mathcal{B}}_{0}(\alpha_1,\ldots,\alpha_\lambda;\eta,k,r)$ by successively applying the $(k-1)$-separation with $a=\eta$ and the backward move.  Let $\pi$ be  an overpartition in $\overline{\mathcal{B}}_1(\alpha_1,\ldots,\alpha_\lambda;\eta,k,r)$.  We shall construct a  pair of overpartitions  $(\zeta,\mu)$, that is,  $\Psi(\pi)=(\zeta,\mu)$, such that $|\zeta|+|\mu|=|\pi|$, where $\zeta \in \mathcal{D}_\eta$ and $\mu \in {\mathcal{B}}_{0}(\alpha_1,\ldots,\alpha_\lambda;\eta,k,r)$.

There are two steps in the construction of $(\zeta,\mu)$ from $\pi$. In the first step,   we   eliminate all overlined parts of $\pi$ divisible by $\eta$  by successively applying the $(k-1)$-separation with $a=\eta$. In the second step,  we  successively apply the backward move to the resulting overpartition  in the first step  so that all $(k-1)$-bands  of the obtained overpartition are  even.

{\noindent Step 1.} Assume that there are $c\geq 0$ overlined parts divisible by $\eta$ in $\pi$. We eliminate the $c$ overlined parts divisible by $\eta$ from $\pi$ by applying  the $(k-1)$-separation with $a=\eta$. Denote the intermediate pairs by  $(\zeta^{(0)},\pi^{(0)}),\ldots,(\zeta^{(c)},\pi^{(c)})$, with  $(\zeta^{(0)},\pi^{(0)})=(\emptyset,\pi)$.  There are two cases:

{\noindent Case 1:}  $c=0$. Then set $\zeta^{(c)}=\emptyset$ and $\pi^{(c)}=\pi$.

{\noindent  Case 2:} $c\geq 1$. Assume that  $\overline{\eta t_0}>\overline{\eta t_1}>\cdots>\overline{\eta t_{c-1}}$ are the overlined parts of $\pi$ divisible by $\eta$. Set $b=0$  and carry out the following procedure.

 \begin{itemize}

\item[(A)]  Let  $\hat{\pi}^{(b)}$ be the  overpartition  obtained  from $\pi^{(b)}$ by removing the overlined part $\overline{\eta t_{b}}$. Assume that $\tilde{g}_1(\hat{\pi}^{(b)})>\cdots>\tilde{g}_{N(\hat{\pi}^{(b)})}
    (\hat{\pi}^{(b)})$  are the $N(\hat{\pi}^{(b)})$ parts marked with $k-1$ in $G(\hat{\pi}^{(b)})$, and $p_{b}$ is the smallest integer such that $ \tilde{g}_{p_{b}+1}(\hat{\pi}^{(b)})<\overline{\eta t_{b}}$. Let $s^{(b)}=p_b+t_b$. By definition,
 \[\pi^{(b)} \in \overline{\mathcal{B}}^{\,\eta}_=(\alpha_1,\ldots,\alpha_\lambda;\eta,k,r|
 N(\hat{\pi}^{(b)}),s^{(b)}-1).\]

Apply the $(k-1)$-separation ${J}^{\,\eta}_{s^{(b)}-1}$  to $\pi^{(b)}$ to get $\pi^{(b+1)}$, that is,
 \[\pi^{(b+1)}={J}^{\,\eta}_{s^{(b)}-1}(\pi^{(b)}).\]
By means of Lemma \ref{isub-lem}, we find that
 \[\pi^{(b+1)}\in \overline{\mathcal{B}}^{\,\eta}_<(\alpha_1,\ldots,\alpha_\lambda;\eta,k,r|
 N(\hat{\pi}^{(b)}),s^{(b)}-1),\]
 and
 \[  \left|\pi^{(b+1)}\right|=|\pi^{(b)}|-\eta s^{(b)}.\]
 Then insert $\eta s^{(b)}$ into $\zeta^{(b)}$ as a part to obtain $\zeta^{(b+1)}$.

 \item[(B)] Replace $b$ by $b+1$. If  $b=c$, then we are done. Otherwise, go back to (A).

 \end{itemize}

Observe that for $0\leq b\leq c$,   there are $c-b$ overlined parts divisible by $\eta$ in $\pi^{(b)}$.
 Theorem  \ref{ssins} reveals  that for $0\leq b< c-1$,
  \begin{equation}\label{ttt1}
  s^{(b)}>s^{(b+1)}>N(\hat{\pi}^{(b+1)}).
  \end{equation}
Therefore,  there are   no overlined parts divisible by $\eta$ in $\pi^{(c)}$ and   $\zeta^{(c)}=(\eta s^{(0)},\ldots, \eta s^{(c-1)})$  is a partition with distinct parts divisible by $\eta$.
 Moreover, we have
  \begin{equation}\label{ttt3}
 |\pi|=|\pi^{(c)}|+|\zeta^{(c)}|.
 \end{equation}

Let us now move on to the second step.

{\noindent Step 2.}  Applying the  backward move  successively to $\pi^{(c)}$, we are led to a pair of overpartitions $(\zeta,\mu)\in\mathcal{D}_\eta\times{\mathcal{B}}_{0}(\alpha_1,\ldots,\alpha_\lambda;\eta,k,r)$ such that $|\pi^{(c)}|=|\mu|+|\zeta|.$ Let $N$ be the number of the $(k-1)$-marked parts in $G(\pi^{(c)})$ and let $\tilde{g}_1(\pi^{(c)})>\cdots>\tilde{g}_{N}(\pi^{(c)})$ be  the $(k-1)$-marked parts in $G(\pi^{(c)})$.  There are two cases:

{\noindent Case 1:} All the $(k-1)$-bands  $\{\tilde{g}_i(\pi^{(c)})\}_{k-1}$  of $\pi^{(c)}$ are even. In view of Theorem \ref{parity k-1 sequence-over-g}, we have  $\pi^{(c)} \in \mathcal{{B}}_0(\alpha_1,\ldots,\alpha_\lambda;\eta,k,r)$. Set $\mu=\pi^{(c)}$ and $\zeta=\zeta^{(c)}$. Then $(\zeta,\mu)\in\mathcal{D}_\eta\times{\mathcal{B}}_{0}(\alpha_1,\ldots,\alpha_\lambda;\eta,k,r)$ and  $|\pi|=|\mu|+|\zeta|.$

{\noindent Case 2:} There exists $i$ such that $1\leq i\leq N$ and $\{\tilde{g}_i(\pi^{(c)})\}_{k-1}$ is odd.

In this case, we set $b=0$ and execute  the following procedure. Denote the intermediate pairs by  $(\zeta^{(c)},\pi^{(c)}),(\zeta^{(c+1)},\pi^{(c+1)})$, and so on.

\begin{itemize}
\item[(A)]   Let $\tilde{g}_1(\pi^{(c+b)})>\cdots>\tilde{g}_{N}(\pi^{(c+b)})$ be the $(k-1)$-marked parts in $G(\pi^{(c+b)})$ and let $1\leq p_{c+b} \leq N$   be the smallest  integer such that $\{\tilde{g}_{p_{c+b}}(\pi^{(c+b)})\}_{k-1}$ and  $\{\tilde{g}_{p_{c+b}+1}(\pi^{(c+b)})\}_{k-1}$ have opposite   parities. By definition,  we get
\[\pi^{(c+b)}\in {\mathcal{B}}_d(\alpha_1,\ldots,\alpha_\lambda;\eta,k,r|N,p_{c+b}).\]
Apply the backward move $\psi_{p_{c+b}}$ to $\pi^{(c+b)}$ to get  $\pi^{(c+b+1)}$, that is,
\[\pi^{(c+b+1)}=\psi_{p_{c+b}}(\pi^{(c+b)}).
\]

By Lemma  \ref{deltagamma-reverse}, we obtain  that
\[
\pi^{(c+b+1)}\in {\mathcal{B}}_e(\alpha_1,\ldots,\alpha_\lambda;\eta,k,r|N,p_{c+b}),
\]
and
\begin{equation}\label{ttt32}
|\pi^{(c+b+1)}|=|\pi^{(c+b)}|-\eta p_{c+b}.
\end{equation}
Then  insert  $\eta p_{c+b}$ into $\zeta^{(c+b)}$ as a part to get a partition $\zeta^{(c+b+1)}$.
 \item[(B)]   Replace $b$ by $b+1$. If   all the $(k-1)$-bands  $\{\tilde{g}_i(\pi^{(c+b)})\}_{k-1}$ of $\pi^{(c+b)}$   are even, then we are done. Otherwise, go back to (A).

\end{itemize}

We claim that during the above procedure,   we have
\begin{equation}\label{budengshi}
N\geq p_{c+b+1}>p_{c+b}.
\end{equation}
Given $b\geq 0$, since $\pi^{(c+b+1)}\in {\mathcal{B}}_e(\alpha_1,\ldots,\alpha_\lambda;\eta,k,r|N,p_{c+b})$, we know that $p_{c+b}$ is the least  integer such that $\{\tilde{g}_{i}(\pi^{(c+b+1)})\}_{k-1}$ have the same parity  for  $1\leq i\leq p_{c+b}+1$.  Whereas $\pi^{(c+b+1)}$ is  in $ {\mathcal{B}}_d(\alpha_1,\ldots,\alpha_\lambda;\eta,k,r|N,p_{c+b+1})$, so that  $p_{c+b+1}$ is the least  integer such that $\{\tilde{g}_{p_{c+b+1}}(\pi^{(c+b+1)})\}_{k-1}$ and $\{\tilde{g}_{p_{c+b+1}+1}(\pi^{(c+b+1)})\}_{k-1}$ have opposite parities. Hence we obtain \eqref{budengshi}, and this proves  the claim.

The relation \eqref{budengshi} ensures that the above procedure terminates after at most $N$ iterations. Assume that it terminates with  $b=m$, that is,   all the $(k-1)$-bands   $\{\tilde{g}_{i}(\pi^{(c+m)})\}_{k-1}$  are even  for $1\leq i\leq N$. Set
\[\mu=\pi^{(c+m)}\quad \text{and} \quad \zeta=\zeta^{(c+m)}=(\eta s_0,\ldots, \eta s_{c-1},\eta p_{c+m-1},\ldots, \eta p_{c}).\]
Utilizing  Theorem \ref{parity k-1 sequence-over-g}, we find that $\mu$ is an overpartition in $ \mathcal{{B}}_0(\alpha_1,\ldots,\alpha_\lambda;\eta,k,r)$.
Observe that $N=N(\hat{\pi}^{(c-1)})$ when $c\geq 1$. In light of  \eqref{ttt1} and \eqref{budengshi},  we conclude that $\zeta$ is a partition with distinct parts divisible by $\eta$. Combining \eqref{ttt3} and \eqref{ttt32}, we have $|\pi|=|\mu|+|\zeta|.$ Therefore, $\Psi$ is a   map from $\mathcal{\overline{B}}_1(\alpha_1,\ldots,\alpha_\lambda;\eta,k,r)$ to  $\mathcal{D}_\eta\times{\mathcal{B}}_{0}(\alpha_1,\ldots,\alpha_\lambda;\eta,k,r)$.

Combining Theorem \ref{deltagammathm} and Theorem \ref{deltagammathmbb}, we obtain that $\Psi(\Phi(\zeta,\mu))=(\zeta, \mu)$ for all $(\zeta, \mu) \in \mathcal{D}_\eta\times{\mathcal{B}}_0(\alpha_1,\ldots,\alpha_\lambda;
\eta,k,r)$ and $\Phi(\Psi(\pi))=\pi$ for all $\pi \in \overline{\mathcal{B}}_1(\alpha_1,\ldots,\alpha_\lambda;
\eta,k,r)$. Hence $\Phi$ is a   bijection between $\mathcal{D}_\eta\times{\mathcal{B}}_0(\alpha_1,\ldots,\alpha_\lambda;
\eta,k,r)$ and $\overline{\mathcal{B}}_1(\alpha_1,\ldots,\alpha_\lambda;\eta,k,r)$. This completes the proof.  \qed

\subsection{An example}

We provide an example to illustrate  the bijection $\Phi$ in Theorem \ref{lem-b-0-over}. Let
\begin{equation*}\label{example-dist}
\zeta=(100,80,50,40,20)
\end{equation*}
  be a partition in $\mathcal{D}_{10}$, and let $\mu$ be an overpartition in ${\mathcal{B}}_0(3,7;10,4,3)$ with the reverse Gordon marking
\begin{equation*}\label{example-mu}
\begin{split}
&RG(\mu)=(\overbrace{{ {\overline{87}}_1,{80}_2},
{ {80}_3}}^{{\{80\}_3}},\overline{67}_1,\overline{63}_2,
\overbrace{{ \overline{57}_1,{50}_2}, { {50}_3}}^{{\{{50}\}_3}},\overline{43}_1,\overbrace{{ \overline{37}_2,\overline{33}_1},{ {30}_3}}
^{{\{30\}_3}},\\[5pt]
&\ \ \ \ \ \ \ \ \ \ \ \ \ \underbrace{{ {20}_1,{20}_2},{ \overline{13}_3}}
_{\{\overline{13}\}_3},\overline{7}_1,\overline{3}_2).
\end{split}
\end{equation*}

The overpartition $\pi=\Phi(\zeta,\mu)$ is obtained by successively applying the forward move and the $3$-insertion with $a=10$. Observe that there are four $3$-marked parts in $RG(\mu)$,
 that is,  $N=4$. We first merge $40$ and $20$ of $\zeta$ into $\mu$ by  successively applying the forward move and then  merge $50$, $80$ and $100$ of $\zeta$ into the resulting overpartition by successively applying the $3$-insertion with $a=10$.

\noindent{   Step 1. Merge $40$ and $20$ of $\zeta$ into $\mu$ by  successively applying the forward move.}

Note that   $\{80\}_3,\,\{50\}_3,\,\{30\}_3,\,\{\overline{13}\}_3$ in $RG(\mu)$  are all even,  so
$\mu \in{\mathcal{B}}_e(3,7;10,4,3|4,4).$

\begin{itemize}

\item Set $\mu^{(0)}=\mu$, and merge $40$ into $\mu^{(0)}$.

  Apply the forward move $\phi_4$ to $\mu^{(0)}$ to get $\mu^{(1)}$, namely, add $\eta=10$ to each of the $3$-marked parts $80,50,30$ and $\overline{13}$ in $RG(\mu^{(0)})$ respectively and rearrange the parts in non-increasing  order to obtain $\mu^{(1)}=( {90},\overline{87},{80},\overline{67}_1,\overline{63}, {60},\overline{57},{50},\overline{43},{40},\overline{37},
  \break$$\overline{33},\overline{23},{20},{20},\overline{7},\overline{3})$. The reverse Gordon marking of $\mu^{(1)}$ is given by
\begin{equation*}
\begin{split}
&RG(\mu^{(1)})=(\overbrace{{ {{90}_1,\overline{87}}_2},
{ {80}_3}}^{{\{80\}_3}},
\overbrace{{ \overline{67}_1,\overline{63}_2}, { {60}_3}}^{{\{{60}\}_3}},\overline{57}_1,\overbrace{{ {50}_2,\overline{43}_1},{ {40}_3}}
^{{\{40\}_3}},\\[5pt]
&\ \ \ \ \ \ \ \ \ \ \ \ \ \ \ \ \  \overline{37}_2,\overline{33}_1,\underbrace{{ \overline{23}_1,{20}_2},{ {20}_3}}
_{\{{20}\}_3},\overline{7}_1,\overline{3}_2).
\end{split}
\end{equation*}

By Lemma \ref{deltagamma}, we deduce that $\mu^{(1)} \in{\mathcal{B}}_d(3,7;10,4,3|4,4)$. Indeed,  $\{80\}_3,\,\{60\}_3,\,\break$$\{40\}_3,\,\{20\}_3$ in $RG(\mu^{(1)})$ are odd. This implies that $\mu^{(1)}\in{\mathcal{B}}_e(3,7;10,4,3|4,2).$

\item Merge $20$ into $\mu^{(1)}$.

Apply the forward move $\phi_2$ to $\mu^{(1)}$ to obtain $\mu^{(2)}$, namely, add $\eta=10$ to each of the $3$-marked parts $80$ and $60$ in $RG(\mu^{(1)})$. We get
\begin{equation*}
\begin{split}
&RG(\mu^{(2)})=(\overbrace{ {90}_1,{90}_2,
{ \overline{87}_3}}^{\{\overline{87}\}_3},
\overbrace{{70}_1,\overline{67}_2, { \overline{63}_3}}^{\{\overline{63}\}_3},\overline{57}_1,\overbrace{{ {50}_2,\overline{43}_1},{ {40}_3}}^{\{40\}_3}
,\\[5pt]
&\ \ \ \ \ \ \ \ \ \ \ \ \ \ \ \ \  \overline{37}_2,\overline{33}_1,\underbrace{{ \overline{23}_1,{20}_2},{ {20}_3}}_{\{20\}_3},\overline{7}_1,\overline{3}_2).
\end{split}
\end{equation*}
Again, it follows from Lemma \ref{deltagamma} that
$\mu^{(2)} \in{\mathcal{B}}_d(3,7;10,4,3|4,2)$. In fact, $\{\overline{87}\}_3$ and $\{\overline{63}\}_3$ in $RG(\mu^{(2)})$ are even, but $\{40\}_3$ in $RG(\mu^{(2)})$ is odd.
\end{itemize}

\noindent{  Step 2. Successively employ  the  $3$-insertion with $a=10$ to merge $50$, $80$ and   $100$ of $\zeta$ into $\mu^{(2)}$.}

 \begin{itemize}

 \item We start with merging $50$ into $\mu^{(2)}$, and  set $s=4$.

There are four 3-marked parts in $RG(\mu^{(2)})$, which are  $\tilde{r}_1(\mu^{(2)})=\overline{87}$, $\tilde{r}_2(\mu^{(2)})=\overline{63}$, $\tilde{r}_3(\mu^{(2)})={40}$ and $\tilde{r}_4(\mu^{(2)})={20}$.
In this occasion,  $p=4$ is the smallest integer such that $\overline{(4-p)\cdot 10+10}=\overline{10}\geq \tilde{r}_{p+1}(\mu^{(2)})+10=-\infty$ and  there are no  overlined parts divisible by $10$ in $\mu^{(2)}$. Hence  $\mu^{(2)}\in\overline{\mathcal{B}}^{\,10}_<(3,7;10,4,3|4,4).$

Apply the $3$-insertion $I^{10}_{4}$ to $\mu^{(2)}$ to get $\mu^{(3)}$. More precisely,  add $\eta=10$ to each of the $3$-marked parts $\overline{87}$, $\overline{63}$, $40$ and $20$ in $RG(\mu^{(2)})$  and then insert $\overline{10}$ into the resulting overpartition as an overlined part.  The resulting reverse Gordon marking reads
\begin{equation*}
\begin{split}
&RG(\mu^{(3)})=({{ {\overline{97}_1,{90}}_2},
{ {90}_3}},
{{ \overline{73}_1,{70}_2}, { \overline{67}_3}},{{ \overline{57}_1,{50}_2},{ {50}_3}},\overline{43}_1,\\[5pt]
&\ \ \ \ \ \ \ \ \ \ \ \ \ \ \ \ \ {{ \overline{37}_2,\overline{33}_1},{ {30}_3}},\overline{23}_1,{20}_2,{{ \overline{10}_1,\overline{7}_2},{ \overline{3}_3}}).
\end{split}
\end{equation*}
Utilizing  Lemma \ref{iadd} gives
$ \mu^{(3)} \in \overline{\mathcal{B}}^{\,10}_=(3,7;10,4,3|4,4)$.

\item   Merge $80$ into $\mu^{(3)}$ and   set $s=7$.

There are five 3-marked parts in $RG(\mu^{(3)})$, to wit,  $\tilde{r}_1(\mu^{(3)})=90$, $\tilde{r}_2(\mu^{(3)})=\overline{67}$, $\tilde{r}_3(\mu^{(3)})={50}$, $\tilde{r}_4(\mu^{(3)})={30}$ and $\tilde{r}_5(\mu^{(3)})=\overline{3}$.  Moreover,   $p=2$ is the smallest integer such that $\overline{(7-p)\cdot 10+10}=\overline{60}\geq \tilde{r}_{p+1}(\mu^{(3)})+10=60$.
Given that $ \mu^{(3)} \in \overline{\mathcal{B}}^{\,10}_=(3,7;10,4,3|4,4)$,    Theorem \ref{ssins} yields  that $ \mu^{(3)} \in \overline{\mathcal{B}}^{\,10}_<(3,7;10,4,3|5,7)$.

Apply the $3$-insertion $I^{10}_{7}$ to $\mu^{(3)}$ to get $\mu^{(4)}$, that is,  add $\eta=10$ to each of the $3$-marked parts ${90}$ and $\overline{67}$ in $RG(\mu^{(3)})$ and then insert $\overline{60}$ into the resulting overpartition as an overlined part. We are led to
\begin{equation*}
\begin{split}
&RG(\mu^{(4)})=({{ {{100}}_1,\overline{97}_2},
{ {90}_3}},
{{ \overline{77}_1,\overline{73}_2}, { {70}_3}},\overline{60}_1,{{ \overline{57}_2,{50}_1},{ {50}_3}},\overline{43}_2,\\[5pt]
&\ \ \ \ \ \ \ \ \ \ \ \ \ \ \ \ \ {{ \overline{37}_1,\overline{33}_2},{ {30}_3}},\overline{23}_1,{20}_2,{{ \overline{10}_1,\overline{7}_2},{ \overline{3}_3}}).
\end{split}
\end{equation*}
As asserted by Lemma \ref{iadd}, we have $ \mu^{(4)} \in \overline{\mathcal{B}}^{\,10}_=(3,7;10,4,3|5,7)$.

\item  Finally,  merge $100$ into $\mu^{(4)}$, and   set $s=9$.

There  are five $3$-marked parts in  $RG(\mu^{(4)})$, namely, $\tilde{r}_1(\mu^{(4)})=90$, $\tilde{r}_2(\mu^{(4)})=70$, $\tilde{r}_3(\mu^{(4)})={50}$, $\tilde{r}_4(\mu^{(4)})={30}$ and $\tilde{r}_5(\mu^{(4)})=\overline{3}$.  Moreover,  $p=0$ is the smallest integer such that $\overline{(9-p)\cdot 10+10}=\overline{100}\geq \tilde{r}_1(\mu^{(4)})+10=100$.  Knowing that $ \mu^{(4)} \in \overline{\mathcal{B}}^{\,10}_=(3,7;10,4,3|5,7)$,    Theorem \ref{ssins} indicates that $ \mu^{(4)} \in \overline{\mathcal{B}}^{\,10}_<(3,7;10,4,3|5,9)$.

Apply the $3$-insertion $I^{10}_{9}$ to $\mu^{(4)}$  to get $\mu^{(5)}$. In other words,   insert $\overline{100}$ into  $\mu^{(4)}$ as an overlined part to generate
 \begin{equation}\label{rrr}
\begin{split}
&RG(\mu^{(5)})=({{ \overline{{100}}_1,{{100}}_2},
{ \overline{97}_3}},{90}_1,
{{ \overline{77}_1,\overline{73}_2}, { {70}_3}},\overline{60}_1,
{{ \overline{57}_2,{50}_1},{ {50}_3}},
\overline{43}_2,\\[5pt]
&\ \ \ \ \ \ \ \ \ \ \ \ \ \ \ \ \ {{ \overline{37}_1,\overline{33}_2},
{ {30}_3}},\overline{23}_1,{20}_2,{{ \overline{10}_1,
\overline{7}_2},{ \overline{3}_3}}
).
\end{split}
\end{equation}
Using Lemma \ref{iadd} again, we conclude that
$ \mu^{(5)} \in \overline{\mathcal{B}}^{\,10}_=(3,7;10,4,3|5,9)$.

\end{itemize}

Set $\pi=\mu^{(5)}$. Clearly,   $\pi$ is an overpartition in $\overline{\mathcal{B}}_1(3,7;10,4,3)$ such that $|\pi|=|\mu|+|\zeta|$.

Conversely, let $\pi$ be   an overpartition in $\overline{\mathcal{B}}_1(3,7;10,4,3)$ whose reverse Gordon marking is given by \eqref{rrr}. The pair of overpartitions $\Psi(\pi)=(\zeta,\mu)$ can be recovered  by successively applying the $3$-separation with $a=10$ and the backward move. There are three overlined parts in $\pi$ divisible by $10$, as identified by  $\overline{100}$, $\overline{60}$ and $\overline{10}$.

\noindent{  Step 1. Eliminate $\overline{100}$, $\overline{60}$ and $\overline{10}$ from $\pi$ by successively  using the $3$-separation with $a=10$.}

\begin{itemize}
\item  Eliminate   $\overline{100}$  from $\pi$, and set $t_0=10$.

    Set $\pi^{(0)}=\pi$ and $\zeta^{(0)}=\emptyset$. Let $\hat{\pi}^{(0)}$ be the   overpartition obtained  from $\pi^{(0)}$ by removing $\overline{100}$, which has the Gordon marking
\begin{equation}\label{example-pi-conv}
\begin{split}
&G(\hat{\pi}^{(0)})=({{ {{100}}_3},
{ \overline{97}_2,{90}_1}},
{{ \overline{77}_3}, { \overline{73}_2,{70}_1}},\overline{60}_2,\overline{57}_1,{{ {50}_3},{ {50}_2,\overline{43}_1}},\overline{37}_2,\overline{33}_1\\[5pt]
&\ \ \ \ \ \ \ \ \ \ \ \ \ \ \ {{ {30}_3},{ \overline{23}_2,{20}_1}},{{ \overline{10}_3},{ \overline{7}_2,\overline{3}_1}}).
\end{split}
\end{equation}
There are five $3$-marked parts in $G(\hat{\pi}^{(0)})$, namely, $\tilde{g}_1(\hat{\pi}^{(0)})=100$, $\tilde{g}_2(\hat{\pi}^{(0)})=\overline{77}$, $\tilde{g}_3(\hat{\pi}^{(0)})=50$, $\tilde{g}_4(\hat{\pi}^{(0)})=30$ and $\tilde{g}_5(\hat{\pi}^{(0)})=\overline{10}$. Moreover, $p_0=0$ is the smallest integer such that $\overline{10\cdot t_0}=\overline{100}>\tilde{g}_{p_0+1}(\hat{\pi}^{(0)})
=100$. Set $s_0=p_0+t_0=10$. Then $\pi^{(0)}\in\overline{\mathcal{B}}^{\,10}_=(3,7;10,4,3|5,9).$

Set  $\zeta^{(1)}=(100)$. Apply the $3$-separation $J^{10}_{9}$  to $\pi^{(0)}$ to get $\pi^{(1)}$. In other words,  $\pi^{(1)}$ is obtained from $\pi^{(0)}$ by  removing  $\overline{100}$, which means that $\pi^{(1)}=\hat{\pi}^{(0)}$ and the Gordon marking of ${\pi}^{(1)}$ is given by \eqref{example-pi-conv}.
Appealing to Lemma \ref{isub-lem}, we deduce that $\pi^{(1)}\in \overline{\mathcal{B}}^{\,10}_<(3,7;10,4,3|5,9)$.

 \item Eliminate   $\overline{60}$  from $\pi^{(1)}$,  and set $t_1=6$.

     Let $\hat{\pi}^{(1)}$ be the   overpartition  obtained  from $\pi^{(1)}$  by removing $\overline{60}$. We have
\begin{equation*}
\begin{split}
&G(\hat{\pi}^{(1)})=({{ {{100}}_3},
{ \overline{97}_2,{90}_1}},
{{ \overline{77}_3}, { \overline{73}_2,{70}_1}},\overline{57}_1,{{ {50}_3},{ {50}_2,\overline{43}_1}},\overline{37}_2,\overline{33}_1\\[5pt]
&\ \ \ \ \ \ \ \ \ \ \ \ \ \ \ {{ {30}_3},{ \overline{23}_2,{20}_1}},{{ \overline{10}_3},{ \overline{7}_2,\overline{3}_1}}).
\end{split}
\end{equation*}
There are five $3$-marked parts in  $G(\hat{\pi}^{(1)})$, which are $\tilde{g}_1(\hat{\pi}^{(1)})=100$, $\tilde{g}_2(\hat{\pi}^{(1)})=\overline{77}$, $\tilde{g}_3(\hat{\pi}^{(1)})=50$, $\tilde{g}_4(\hat{\pi}^{(1)})=30$ and $\tilde{g}_5(\hat{\pi}^{(1)})=\overline{10}$. Now, $p_1=2$ is the smallest integer such that  $\overline{10\cdot t_1}=\overline{60}>\tilde{g}_{p_1+1}(\hat{\pi}^{(1)})
=50$. Set $s_1=p_1+t_1=8$, and we get
$\pi^{(1)}\in\overline{\mathcal{B}}^{\,10}_=(3,7;10,4,3|5,7).$  Clearly,  $s_0>s_1$, in agreement with Theorem \ref{ssins}.

Set $\zeta^{(2)}=(100,80).$ Apply the $3$-separation $J^{10}_{7}$  to $\pi^{(1)}$ to get $\pi^{(2)}$, namely,   remove $\overline{60}$ from $\pi^{(1)}$ to get $\hat{\pi}^{(1)}$, and then subtract $\eta=10$ from each of the $3$-marked parts $100$ and $\overline{77}$ in $G(\hat{\pi}^{(1)})$  to get $\pi^{(2)}$.  The Gordon marking of ${\pi}^{(2)}$ is given below:
\begin{equation*}
\begin{split}
&G(\pi^{(2)})=({{{ \overline{97}_3}},
{ {90}_2,{90}_1}},
{{ \overline{73}_3}, { {70}_2,\overline{67}_1}},\overline{57}_1,{{ {50}_3},{ {50}_2,\overline{43}_1}},\overline{37}_2,\overline{33}_1\\[5pt]
&\ \ \ \ \ \ \ \ \ \ \ \ \ \ {{ {30}_3},{ \overline{23}_2,{20}_1}},{{ \overline{10}_3},{ \overline{7}_2,\overline{3}_1}}).
\end{split}
\end{equation*}
We now have $\pi^{(2)}\in \overline{\mathcal{B}}^{\,10}_<(3,7;10,4,3|5,7)$, as expected by Lemma \ref{isub-lem}.

\item Finally, eliminate   $\overline{10}$ from $\pi^{(2)}$, and set $t_2=1$.

Let $\hat{\pi}^{(2)}$ be the  overpartition  obtained  from $\pi^{(2)}$ by removing $\overline{10}$, so that
\begin{equation*}
\begin{split}
&G(\hat{\pi}^{(2)})=({{{ \overline{97}_3}},
{ {90}_2,{90}_1}},
{{ \overline{73}_3}, { {70}_2,\overline{67}_1}},\overline{57}_1,{{ {50}_3},{ {50}_2,\overline{43}_1}},\overline{37}_2,\overline{33}_1\\[5pt]
&\ \ \ \ \ \ \ \ \ \ \ \ \ \ {{ {30}_3},{ \overline{23}_2,{20}_1}},\overline{7}_2,\overline{3}_1).
\end{split}
\end{equation*}
There are four $3$-marked parts in  $G(\hat{\pi}^{(2)})$, namely, $\tilde{g}_1(\hat{\pi}^{(2)})=\overline{97}$, $\tilde{g}_2(\hat{\pi}^{(2)})=\overline{73}$, $\tilde{g}_3(\hat{\pi}^{(2)})=50$ and $\tilde{g}_4(\hat{\pi}^{(2)})=30$. Meanwhile, $p_2=4$ is the smallest integer such that  $\overline{10\cdot t_2}=\overline{10}>\tilde{g}_{p_2+1}(\hat{\pi}^{(1)})
=-\infty$. Set $s_2=t_2+p_2=5$. Then
$\pi^{(2)}\in\overline{\mathcal{B}}^{\,10}_=(3,7;10,4,3|4,4).$  In accordance with Theorem   \ref{ssins}, we have $s_1>s_2$.

 Set $\zeta^{(3)}=(100,80,50).$ Apply the $3$-separation $J^{10}_{4}$ to $\pi^{(2)}$ to get $\pi^{(3)}$. To wit,   remove $\overline{10}$ from $\pi^{(2)}$ to get $\hat{\pi}^{(2)}$, then  subtract $\eta=10$ from each of  the $3$-marked parts $\overline{97}$, $\overline{73}$, $50$ and $30$ in $G(\hat{\pi}^{(2)})$  to obtain $\pi^{(3)}$.  We get
\begin{equation*}
\begin{split}
&G(\pi^{(3)})=(\overbrace{{{ {90}_3}},
{ {90}_2,\overline{87}_1}}^{{\{{90}\}_3}},
\overbrace{{ {70}_3}, { \overline{67}_1,\overline{63}_2}}^{{\{{70}\}_3}},\overline{57}_1,{50}_2,\overline{43}_1,\overbrace{{ {40}_3},{ \overline{37}_2,\overline{33}_1}}
^{{\{40\}_3}},\\[5pt]
&\ \ \ \ \ \ \ \ \ \ \ \ \ \ \ \underbrace{{ \overline{23}_3},{ {20}_2,{20}_1}}_{\{\overline{23}\}_3},\overline{7}_2,\overline{3}_1).
\end{split}
\end{equation*}
Using Lemma \ref{isub-lem}, we have
$\pi^{(3)}\in \overline{\mathcal{B}}^{\,10}_=(3,7;10,4,3|4,4)$.

\end{itemize}
There are no  overlined parts  divisible by $10$ in $\pi^{(3)}$.  The fact that $\zeta^{(3)}=(100,80,50)$ is a partition with distinct parts reflects the claim of Theorem  \ref{ssins}.

\noindent{ Step 2. Successively apply the backward move to $\pi^{(3)}$ to derive a pair of overpartitions $(\zeta,\mu)$ in  ${\mathcal{D}}_\eta \times {\mathcal{B}}_0(3,7;10,4,3)$.}

There are four  $3$-marked parts in $G(\pi^{(3)})$, namely, $\tilde{g}_1(\pi^{(3)})={90}$, $\tilde{g}_2(\pi^{(3)})={70}$, $\tilde{g}_3(\pi^{(3)})=40$ and $\tilde{g}_4(\pi^{(3)})=\overline{23}$. Moreover,     $\{90\}_3$ and $\{70\}_3$  are even and $\{40\}_3$ and $\{\overline{23}\}_3$   are odd, whereas $p_3=2$ is the smallest integer such that  $\{\tilde{g}_{p_3}(\pi^{(3)})\}_{3}$ and $\{\tilde{g}_{p_3+1}(\pi^{(3)})\}_{3}$ have opposite  parities. Hence $\pi^{(3)} \in{\mathcal{B}}_d(3,7;10,4,3|4,2).$

 \begin{itemize}
\item Set $\zeta^{(4)}=(100,80,50,20).$ Apply the backward move $\psi_2$ to $\pi^{(3)}$ to produce $\pi^{(4)}$. Strictly speaking,  subtract $\eta=10$ from each of the $3$-marked parts ${90}$ and ${70}$ in $G(\pi^{(3)})$ to get $\pi^{(4)}$.  The Gordon marking of $\pi^{(4)}$ is
\begin{equation*}
\begin{split}
&G(\pi^{(4)})=(\overbrace{{{ {90}_3}},
{ \overline{87}_2,{80}_1}}^{{\{{90}\}_3}},\overline{67}_1,\overline{63}_2,
\overbrace{{ {60}_3}, { \overline{57}_1,{50}_2}}^{{\{{60}\}_3}},\overline{43}_1,\overbrace{{ {40}_3},{ \overline{37}_2,\overline{33}_1}}
^{{\{40\}_3}},\\[5pt]
&\ \ \ \ \ \ \ \ \ \ \ \ \ \ \ \underbrace{{ \overline{23}_3},
{ {20}_2,{20}_1}}_{\{\overline{23}\}_3},
\overline{7}_2,\overline{3}_1).
\end{split}
\end{equation*}
In view of Lemma  \ref{deltagamma-reverse}, we may say that $\pi^{(4)}\in{\mathcal{B}}_e(3,7;10,4,3|4,2)$. To be more specific, $\{90\}_3$, $\{60\}_3$ $\{40\}_3$ and $\{\overline{23}\}_3$ are all odd. Hence $p_4=4$ is the smallest integer such that  $\{\tilde{g}_{p_4}(\pi^{(4)})\}_{3}$ and $\{\tilde{g}_{p_4+1}(\pi^{(4)})\}_{3}$ have opposite parities. It follows that $\pi^{(4)} \in{\mathcal{B}}_d(3,7;10,4,3|4,4).$ Obviously, $N\geq p_4>p_3$.

\item Set  $\zeta^{(5)}=(100,80,50,40,20).$ Apply the backward move $\psi_4$ to $\pi^{(4)}$ to obtain $\pi^{(5)}$,  namely, subtract $\eta=10$ from each of the $3$-marked parts ${90}$, ${60}$, $40$ and $\overline{23}$ in $G(\pi^{(4)})$. We get
\begin{equation*}
\begin{split}
&G(\pi^{(5)})=(\overbrace{{ \overline{87}_3},
{80}_2,{80}_1}^{\{\overline{87}\}_3},\overline{67}_1,\overline{63}_2,
\overline{57}_1,\overbrace{{ {50}_3},  \overline{50}_2,\overline{43}_1}^{\{50\}_3},\overbrace{{ \overline{37}_3}, \overline{33}_2,{30}_1}^{\{\overline{37}\}_3},\\[5pt]
&\ \ \ \ \ \ \ \ \ \ \ \ \ \ \underbrace{{ {20}_3},{20}_2,
\overline{13}_1}_{\{20\}_3},
\overline{7}_2,\overline{3}_1).
\end{split}
\end{equation*}
By Lemma  \ref{deltagamma-reverse}, we see that $\pi^{(5)} \in{\mathcal{B}}_e(3,7;10,4,3|4,4)$. More precisely,  $\{\overline{87}\}_3,\,\{50\}_3,\,\break$
$\{\overline{37}\}_3,\,\{20\}_3$ in $G(\pi^{(5)})$   are even. Resorting to  Theorem \ref{parity k-1 sequence-over-g}, we arrive at  $\pi^{(5)}\in{\mathcal{B}}_0(3,7;10,4,3)$.
\end{itemize}

In conclusion, set   $\zeta=\zeta^{(5)}$ and $\mu=\pi^{(5)}$.  Then $(\zeta,\mu)\in \mathcal{D}_{10}\times {\mathcal{B}}_0(3,7;10,4,3)$ and $|\pi|=|\mu|+|\zeta|$.

\section{Proof of Theorem \ref{rel-over2} }

The goal of this section  is to give a proof of Theorem \ref{rel-over2}, which can be restated in  purely    combinatorial terms. Here we use the common notation $\delta_{r,k}=1$ if $r=k$, and $\delta_{r,k}=0$ otherwise.

\begin{thm}\label{lem-b-1-over} Let $k$, $r$  and $\lambda$ be integers  such that $k\geq r\geq\lambda\geq0$ and $k-1>\lambda$.
There is a bijection $\Theta$ between  $\mathcal{\overline{B}}_0(\alpha_1,\ldots,\alpha_\lambda;\eta,k,r)$ and $\mathcal{D}_\eta\times{\mathcal{B}}_{1}(\alpha_1,\ldots,\alpha_\lambda;\eta,k-1,r-\delta_{r,k})$, namely, for an overpartition $\nu \in\mathcal{\overline{B}}_0(\alpha_1,\ldots,\alpha_\lambda;\eta,k,r)$, we have $\Theta(\nu)=(\zeta,\omega) \in \mathcal{D}_\eta\times\mathcal{{B}}_{1}(\alpha_1,\ldots,\alpha_\lambda;\eta,k-1,r-\delta_{r,k}) $ such that  $|\nu|=|\zeta|+|\omega|$ and $\ell(\nu)=\ell(\zeta)+\ell(\omega)$.
\end{thm}

 Since there are no  overlined parts divisible by $\eta$ in $\omega$ and there are no $(k-1)$-marked parts in   $RG(\omega)$, in order to obtain $(\zeta,\omega)$, we need to remove all overlined parts divisible by $\eta$ and certain non-overlined parts divisible by $\eta$ from $\nu$ to generate $\omega$, and use the removed parts to generate $\zeta$.  To this end, we shall define the  $(k-1)$-reduction operation  and   the $(k-1)$-augmentation operation, which are the main ingredients in the construction of $\Theta$.

\subsection{The $(k-1)$-reduction and the $(k-1)$-augmentation }

 The definitions of the  $(k-1)$-reduction and the $(k-1)$-augmentation are based on two subsets of $\overline{\mathcal{B}}_0(\alpha_1,\ldots,\alpha_\lambda;\eta,k,r)$. To describe these two subsets, we need
to introduce the following notation.  Define $ol(\nu)$ to be the largest overlined part divisible by $\eta$ in $\nu$ with the convention that $ol(\nu)=\overline{0}$ if there are no   overlined parts  divisible by $\eta$ in $\nu$. Define $\tilde{r}_1(\nu)$ to be the largest $(k-1)$-marked part  in   $RG(\nu)$ with the convention that  $\tilde{r}_1(\nu)=-\infty$ if there  are no $(k-1)$-marked parts  in   $RG(\nu)$.

 We now assume that $k$, $r$  and $\lambda$ are integers  such that $k\geq r\geq\lambda\geq0$ and $k-1>\lambda$.

\begin{itemize}

\item  For $t\geq 1$, let $\overline{\mathcal{B}}^{\,=}_0(\alpha_1,\ldots,\alpha_\lambda;\eta,k,r|t)$ denote the set of   overpartitions $\nu$ in $\overline{\mathcal{B}}_0(\alpha_1,\ldots,\break $$\alpha_\lambda;\eta,k,r)$ such that either $ol(\nu)=\overline{t\eta}$ and $\tilde{r}_1(\nu)\leq\overline{t\eta}$, or $ol(\nu)<\overline{t\eta}$ and  $\overline{(t-1)\eta}<\tilde{r}_1(\nu)\leq{t\eta}$.

\item For $t\geq 1$, let $\overline{\mathcal{B}}^{\,<}_0(\alpha_1,\ldots,\alpha_\lambda;\eta,k,r|t)$ denote the set of   overpartitions $\nu$ in $\overline{\mathcal{B}}_0(\alpha_1,\ldots,\break$$\alpha_\lambda;\eta,k,r)$ such that
 $ol(\nu)<\overline{t\eta}$ and $\tilde{r}_1(\nu)\leq \overline{(t-1)\eta}$.
\end{itemize}

With  the above two subsets in hand, we are ready to give the definition of the $(k-1)$-reduction operation.

\begin{defi}[The $(k-1)$-reduction]\label{defi-division} For  $t\geq 1$,   let $\nu$ be an overpartition in $\overline{\mathcal{B}}^{\,=}_0(\alpha_1,\ldots,\break$
$\alpha_\lambda;\eta,k,r|t)$. Define the $(k-1)$-reduction $D_t\colon \nu \rightarrow \omega$ as follows{\rm{:}}    If $ol(\nu)=\overline{t\eta}$, then $\omega$ is obtained  from $\nu$ by removing the overlined part $\overline{t\eta}$. Otherwise, $\omega$ is obtained  from $\nu$ by removing a non-overlined part   ${t\eta}$.
\end{defi}

The following proposition guarantees that the $(k-1)$-reduction is well defined.

\begin{prop}\label{daxiao-0-5}
For $t\geq1$, let $\nu$ be an overpartition in $\overline{\mathcal{B}}_0(\alpha_1,\ldots,\alpha_\lambda;\eta,k,r)$ such that $ol(\nu)<\overline{t\eta}$ and  $\overline{(t-1)\eta}<\tilde{r}_1(\nu)\leq{t\eta}$. Then $\nu$ contains a non-overlined part  $t\eta$.
\end{prop}

\pf Assume that $\tilde{r}_1(\nu)$ is the $r_1$-th part of $\nu=(\nu_1,\nu_2,\ldots,\nu_\ell)$ in $\overline{\mathcal{B}}_0(\alpha_1,\ldots,\alpha_\lambda;\eta,k,r)$, that is, $\tilde{r}_1(\nu)=\nu_{r_1}$. Since $\nu_{r_1}$ is a $(k-1)$-marked part in $RG(\nu)$, there is a unique $(k-1)$-band of $\nu$ induced by $\nu_{r_1}$. Assume that
\[\nu_{r_1-k+2}\geq \nu_{r_1-k+3}\geq\cdots\geq\nu_{r_1}\]
are the parts in the $(k-1)$-band induced by $\nu_{r_1}$, where $\nu_{r_1-k+2}\leq\nu_{r_1}+\eta$ with strict inequality if $\nu_{r_1}$ is overlined. Under the condition  $\overline{(t-1)\eta}<\nu_{r_1}\leq{t\eta}$, we deduce that $\nu_{r_1-k+2}\leq (t+1)\eta$, and so
\begin{equation*}
(t+1)\eta\geq\nu_{r_1-k+2}\geq \nu_{r_1-k+3}\geq\cdots\geq\nu_{r_1}> \overline{(t-1)\eta}.
\end{equation*}
Moreover, we may assume that  $m$ is the smallest integer such that $r_1-k+2\leq m\leq r_1$ and $\nu_{m}\leq t\eta$. This implies that $\nu_{l}\leq t\eta$ for $m\leq l\leq  r_1$ and $\nu_{l}\geq \overline{t\eta}$ for $r_1-k+2\leq  l< m$.  We claim that  $\nu_{m}=t\eta$. Suppose to the contrary that  $\nu_{m}< t\eta$. In this case, we have $\overline{(t-1)\eta}<\nu_{r_1}\leq \nu_{m}<t\eta$, so we can write  $\nu_{r_1}=\overline{(t-1)\eta+\alpha_i}$,
where $1\leq i\leq \lambda$. Then  we have
\begin{equation}\label{pro-examaa}
 \overline{(t-1)\eta+\alpha_i}=\nu_{r_1}\leq \cdots\leq \nu_{m}\leq \overline{(t-1)\eta+\alpha_\lambda},
\end{equation}
and  $\nu_{r_1-k+2}<\nu_{r_1}+\eta=\overline{t\eta+\alpha_i}$. The condition   $ol(\nu)<\overline{t\eta}$ implies that  $\nu_{l}> \overline{t\eta}$ for $r_1-k+2\leq  l< m$. Hence we have
\begin{equation}\label{pro-exambb}
 \overline{t\eta}<\nu_{m-1}\leq \cdots\leq \nu_{r_1-k+2}< \overline{t\eta+\alpha_i}.
\end{equation}
Combining \eqref{pro-examaa} and \eqref{pro-exambb}, we deduce that $k-1\leq (\lambda-i+1)+(i-1)=\lambda$,  which contradicts the assumption that $k-1>\lambda$. Hence  $\nu_m=t\eta$. This completes the proof. \qed

The following theorem says that the $(k-1)$-reduction operation is  indeed a bijection.

\begin{thm}\label{theorem4.4}  For  $t\geq 1$,  the $(k-1)$-reduction $D_{t}$ is a bijection between   $\overline{\mathcal{B}}^{\,=}_0(\alpha_1,\ldots,
\alpha_\lambda;\break$$\eta,k,r|t)$ and  $\overline{\mathcal{B}}^{\,<}_0(\alpha_1,\ldots,\alpha_\lambda;\eta,k,r|t)$. Moreover, for $\nu \in \overline{\mathcal{B}}^{\,=}_0(\alpha_1,\ldots,\alpha_\lambda;\eta,k, r|t)$ and $\omega=D_t(\nu) $, we have   $|\omega|=|\nu|-t\eta$ and $\ell(\omega)=\ell(\nu)-1$.
\end{thm}

The proof of Theorem \ref{theorem4.4} consists of three parts.
In  Lemma \ref{division}, we show that  the $(k-1)$-reduction is a map from $\overline{\mathcal{B}}^{\,=}_0(\alpha_1,\ldots,
\alpha_\lambda;\eta,k,r|t)$ to  $\overline{\mathcal{B}}^{\,<}_0(\alpha_1,\ldots,\alpha_\lambda;\eta,k,r|t)$. Lemma \ref{combinatio} exhibits
 the $(k-1)$-augmentation map from $\overline{\mathcal{B}}^{\,<}_0(\alpha_1,\ldots,
 \alpha_\lambda;\eta,k,r|t)$ to $\overline{\mathcal{B}}^{\,=}_0(\alpha_1,\ldots,
 \alpha_\lambda;\eta,k,\break$$r|t)$. Then we show
 that the $(k-1)$-reduction and the $(k-1)$-augmentation are inverses of each other.

 \begin{lem}\label{division} For  $t\geq 1$,  let $\nu$ be an overpartition in $\overline{\mathcal{B}}^{\,=}_0(\alpha_1,\ldots,\alpha_\lambda;\eta,k, r|t)$ and let $\omega=D_t(\nu)$. Then $\omega$ is overpartition in $\overline{\mathcal{B}}^{\,<}_0(\alpha_1,\ldots,\alpha_\lambda;\eta,k,r|t)$. Furthermore, $|\omega|=|\nu|-t\eta$ and $\ell(\omega)=\ell(\nu)-1$.
\end{lem}
\pf  By definition, we wish to show that $\omega$ satisfies the following conditions:
\begin{itemize}
    \item[(A)] $\omega$ is an overpartition in $\overline{\mathcal{B}}_1(\alpha_1,\ldots,\alpha_\lambda;\eta,k,r)$;

    \item [(B)] $ol(\omega)<\overline{t\eta}$ and $\tilde{r}_1(\omega)\leq \overline{(t-1)\eta}$;

    \item [(C)] All the $(k-1)$-bands  of $\omega$ induced by the $(k-1)$-marked parts in $RG(\omega)$ are even.
\end{itemize}

{\noindent Condition (A).} Given  the precondition  $\nu \in \overline{\mathcal{B}}^{\,=}_0(\alpha_1,\ldots,\alpha_\lambda;\eta,k, r|t)$, it is immediate from the construction of $\omega$   that it satisfies (1)-(4) in the definition of $\overline{\mathcal{B}}_1(\alpha_1,\ldots,\alpha_\lambda;\eta,k,r)$. That is to say,  $\omega$ is an overpartition in $\overline{\mathcal{B}}_1(\alpha_1,\ldots,\alpha_\lambda;\eta,k,r)$.

{\noindent Condition (B).} Since $ol(\nu)\leq \overline{t\eta}$,  $\tilde{r}_1(\nu)\leq\overline{t\eta}$ and $\omega$ is obtained  from $\nu$ by removing an overlined part $\overline{t\eta}$ or a non-overlined part $t\eta$, we obtain that  $ol(\omega)< \overline{t\eta}$ and  $\tilde{r}_1(\omega)\leq{t\eta}$.

We further show that $\tilde{r}_1(\omega)\leq \overline{(t-1)\eta}$. Suppose to the contrary  that
$\tilde{r}_1(\omega)>\overline{(t-1)\eta}$. In this case, we have $\overline{(t-1)\eta}<\tilde{r}_1(\omega)\leq {t\eta}$.
Since $\tilde{r}_1(\omega)$ is the largest  $(k-1)$-marked part in $RG(\omega)$,  there are exactly $k-2$ parts of $\omega$ appearing before $\tilde{r}_1(\omega)$  in the interval $I(\tilde{r}_1(\omega),\tilde{r}_1(\omega)+\eta)$.   The assumption   $\overline{(t-1)\eta}<\tilde{r}_1(\omega)\leq {t\eta}$ implies that $\overline{t\eta}<\tilde{r}_1(\omega)+\eta\leq {(t+1)\eta}$.   Hence the removed part of $\nu$ (that is, ${t\eta}$ or $\overline{t\eta}$) is also in the interval $I(\tilde{r}_1(\omega),\tilde{r}_1(\omega)+\eta)$. It follows that  there are exactly $k-1$ parts  of $\nu$  appearing  before $\tilde{r}_1(\omega)$ in the interval $I(\tilde{r}_1(\omega),\tilde{r}_1(\omega)+\eta)$. This means that  there exists a part of $\nu$  marked with $k$ in $RG(\nu)$, which is impossible because  $\nu \in \overline{\mathcal{B}}^{\,=}_0(\alpha_1,\ldots,\alpha_\lambda;\eta,k, r|t)$.

{\noindent Condition (C).}   Given a $(k-1)$-marked part $\omega_i$ in $RG(\omega)$,   assume that
 $\{\omega_{i-l}\}_{0\leq l\leq k-2}$ is the $(k-1)$-band induced by the $(k-1)$-marked part $\omega_i$. We aim to show that $\{\omega_{i-l}\}_{0\leq l\leq k-2}$ is  even in $\omega$.  Using the condition (B), we know that $\tilde{r}_1(\omega)\leq \overline{(t-1)\eta}$, and so $\omega_{i}\leq \overline{(t-1)\eta}$. The assumption that $\{\omega_{i-l}\}_{0\leq l\leq k-2}$ is a $(k-1)$-band yields $\omega_{i-k+2}<\overline{t\eta}$, more precisely,
    \begin{equation}\label{interttcc}
    \overline{t\eta}>\omega_{i-k+2}\geq \omega_{i-k+3}\geq \cdots \geq \omega_{i}.
    \end{equation}
  It follows that $\{\omega_{i-l}\}_{0\leq l\leq k-2}$ is  also a $(k-1)$-band in $\nu$. Since $\nu\in \overline{\mathcal{B}}_0(\alpha_1,\ldots,\alpha_\lambda;\eta,k,r)$, we find that  the $(k-1)$-band $\{\omega_{i-l}\}_{0\leq l\leq k-2}$  is even in $\nu$, that is,
  \begin{equation}\label{div-pro}
 \left[\omega_{i-k+2}/\eta\right]+\cdots+\left[\omega_{i}/\eta\right]\equiv r-1+\overline{V}_{\nu}(\omega_{i-k+2})\pmod{2}.
 \end{equation}
Noting that $\omega$ is obtained  from $\nu$ by removing an overlined part $\overline{t\eta}$ or a non-overlined part $t\eta$, by   \eqref{interttcc}, we get $\overline{V}_\nu(\omega_{i-k+2})=\overline{V}_\omega(\omega_{i-k+2})$.
  Therefore, it is immediate from   \eqref{div-pro} that $\{\omega_{i-l}\}_{0\leq l\leq k-2}$ is also even in $\omega$, and so the condition (C) is justified.

In conclusion,  we have shown that   $\omega$ is an overpartition in  $\overline{\mathcal{B}}^{\,<}_0(\alpha_1,\ldots,\alpha_\lambda;\eta,k,r|t)$.  Clearly,   $|\omega|=|\nu|-t\eta$ and $\ell(\omega)=\ell(\nu)-1$. This completes the proof.  \qed

We now turn to  the $(k-1)$-augmentation operation, which will be shown to be the inverse map of
the $(k-1)$-reduction operation.

\begin{defi}[The $(k-1)$-augmentation]\label{insertion} For $t\geq 1$, let $\omega$ be an overpartition in $\overline{\mathcal{B}}^{\,<}_0(\alpha_1,\ldots,\alpha_\lambda;\eta,k,r|t)$. We say that $\omega$ satisfies the   condition U if   there exist $k-2$ consecutive parts $\omega_i,\ldots,\omega_{i+k-3}$ of $\omega$ such that
 \begin{itemize}
\item[{\rm(1)}] $(t+1)\eta\geq \omega_i\geq \cdots \geq \omega_{i+k-3}>\overline{(t-1)\eta}${\rm{;}}

\item[{\rm (2)}] $\omega_i\leq \omega_{i+k-3}+\eta$ with strict inequality if $\omega_i$ is overlined{\rm{;}}

\item[{\rm(3)}] $[\omega_{i}/\eta]+\cdots+[\omega_{i+k-3}/\eta]\equiv t+r-1+\overline{V}_{\omega}(\omega_{i})\pmod2.$
\end{itemize}
The $(k-1)$-augmentation  $C_t\colon \omega \rightarrow \nu$ is defined as follows{\rm{:}} If $\omega$ satisfies the condition U,   then $\nu$ is obtained by inserting $t\eta$  into $\omega$  as a non-overlined part. Otherwise, we say that
$\omega$ satisfies the condition O and $\nu$ is obtained by inserting $\overline{t\eta}$   into $\omega$  as an overlined part.

\end{defi}

The following lemma says that the   $(k-1)$-augmentation is a  map from $\overline{\mathcal{B}}^{\,<}_0(\alpha_1,\ldots,
\alpha_\lambda;\eta,\break$$ k,r|t)$ to $\overline{\mathcal{B}}^{\,=}_0(\alpha_1,\ldots,\alpha_\lambda;\eta, k,r|t)$.

\begin{lem}\label{combinatio}    For  $t\geq 1$,  let $\omega$ be an overpartition in $\overline{\mathcal{B}}^{\,<}_0(\alpha_1,\ldots,\alpha_\lambda;\eta,k, r|t)$ and let $\nu=C_t(\omega)$. Then   $\nu$ is an overpartition in $\overline{\mathcal{B}}^{\,=}_0(\alpha_1,\ldots,\alpha_\lambda;\eta,k,r|t)$ such that $|\nu|=|\omega|+t\eta$ and $\ell(\nu)=\ell(\omega)+1$.

\end{lem}

\pf  To prove that $\nu$ is an overpartition in $\overline{\mathcal{B}}^{\,=}_0(\alpha_1,\ldots,\alpha_\lambda;\eta,k,r|t)$, we need to verify that $\nu$ satisfies the following conditions:

\begin{itemize}

\item[(A)] $f_{\leq \eta}(\nu)\leq r-1$;

\item[(B)] For $1\leq i\leq \ell(\nu)-k+1$, $\nu_i\geq \nu_{i+k-1}+\eta$ with strict inequality if $\nu_i$ is non-overlined;

\item[(C)] $ol(\nu)=\overline{t\eta}$ and $\tilde{r}_1(\nu)\leq\overline{t\eta}$, or $ol(\nu)<\overline{t\eta}$ and  $\overline{(t-1)\eta}<\tilde{r}_1(\nu)\leq{t\eta}$;

\item[(D)] All $(k-1)$-bands  of $\nu$  are even.
\end{itemize}

{\noindent Condition (A).}  It is clear that  $f_{\leq \eta}(\omega)\leq r-1$ since $\omega  \in \overline{\mathcal{B}}^{\,<}_0(\alpha_1,\ldots,\alpha_\lambda;\eta,k, r|t)$. To prove that $f_{\leq \eta}(\nu)\leq r-1$, we consider three cases:

{\noindent Case  1.}  $t\geq 2$. In this case,   $f_{\leq \eta}(\nu)=f_{\leq \eta}(\omega)\leq r-1$.

{\noindent Case  2.} $t=1$ and $f_{\leq \eta}(\omega)<r-1$. We have $f_{\leq \eta}(\nu)\leq f_{\leq \eta}(\omega)+1\leq r-1$.

{\noindent Case  3.} $t=1$ and $f_{\leq \eta}(\omega)=r-1$. We claim that  $\omega$  satisfies the condition O.  Assume that
\begin{equation}\label{construction-comb-2tt}
 2\eta\geq \omega_i\geq \cdots \geq \omega_{i+k-3}>0
\end{equation}
are the $k-2$ consecutive parts of $\omega$  such that $\omega_i\leq \omega_{i+k-3}+\eta$  with strict inequality whenever $\omega_i$ is overlined.  Since $\tilde{r}_1(\omega)\leq\overline{(t-1)\eta}=\overline{0}$, there are no $(k-1)$-marked parts in $RG(\omega)$, that is, there are no $(k-1)$-bands  of $\omega$, which implies that $f_{\leq \eta}(\omega)\leq k-2$. Therefore, all parts of $\omega$ not exceeding  $\eta$ are after $\omega_{i-1}$. Hence, by \eqref{construction-comb-2tt}, we obtain that
\[\begin{split}
&[\omega_i/\eta]+\cdots+[\omega_{i+k-3}/\eta]\\[5pt]
&\qquad \equiv f_{=\eta}(\omega)+(\overline{V}_\omega(\omega_i)-f_{<\eta}(\omega))\\[5pt]
&\qquad \equiv f_{\leq \eta}(\omega)+\overline{V}_\omega(\omega_i) \\[5pt]
&\qquad =r-1+\overline{V}_\omega(\omega_i)\pmod{2}.
\end{split}
\]
So the claim is confirmed, and hence $\nu$ is obtained by inserting $\overline{\eta}$ into $\omega$  as an overlined part,  from which we get $f_{\leq \eta}(\nu)=f_{\leq \eta}(\omega)=r-1$.

{\noindent Condition (B).}  Suppose to the contrary that   there exists  $1\leq c\leq \ell(\nu)-k+1$  such that
\begin{equation}\label{lema-ctd}
 \nu_c\leq \nu_{c+k-1}+\eta  \text{ with strict inequality if } \nu_c  \text{ is overlined.}
\end{equation}
 Assume that the $m$-th part of $\nu$ is the inserted part of the $(k-1)$-augmentation operation, that is, $\nu_m=\overline{t\eta}$ or $t\eta$. Since $\omega\in\overline{\mathcal{B}}^{\,<}_0(\alpha_1,\ldots,\alpha_\lambda;\eta,k, r|t)$, we have  $\nu_c\geq \nu_m\geq \nu_{c+k-1}$. Comparing with \eqref{lema-ctd}, we get   $\nu_c\leq \nu_{m}+\eta=(t+1)\eta$ and $\nu_{c+k-1}\geq \nu_m-\eta=(t-1)\eta$  with strict inequality if $\nu_m$ is overlined. Thus, we arrive at
\[(t+1)\eta\geq\nu_c\geq \nu_{c+1}\geq \cdots \geq \nu_{c+k-1}\geq(t-1)\eta.\]
The condition  $\tilde{r}_1(\omega)\leq\overline{(t-1)\eta}$ implies that there are no  $(k-1)$-bands  of $\omega$ in   $(\overline{(t-1)\eta},(t+1)\eta]$.   It follows that  $\nu_{c+k-1}=\overline{(t-1)\eta}$ or $(t-1)\eta$,  and so $\nu_c\leq t\eta$. But $\nu_c\geq \nu_m\geq t\eta$, we obtain that
\begin{equation}\label{lema-ctdaadd}
 \nu_c=\nu_m=t\eta.
\end{equation}
Hence,    $\omega_{c+i-1}=\nu_{c+i}$, where $1\leq i\leq k-1$.  More precisely,
\[t\eta\geq\omega_c\geq \omega_{c+1}\geq \cdots \geq \omega_{c+k-2}\geq (t-1)\eta,\]
and
\begin{equation}\label{lema-ctdaa}
 \omega_{c+k-2}=\nu_{c+k-1}=\overline{(t-1) \eta} \text{ or } (t-1)\eta.
\end{equation}
Therefore,
$\{\omega_{c+l}\}_{0\leq l\leq k-2}$ is   a $(k-1)$-band of $\omega$ in $[(t-1)\eta,t\eta]$.
  Using the condition   $\omega\in\overline{\mathcal{B}}^{\,<}_0(\alpha_1,\ldots,\alpha_\lambda;\eta,k, r|t)$, we find that
\begin{equation}\label{construction-comb-2}
  [\omega_{c}/\eta]+\cdots+[\omega_{c+k-2}/\eta]\equiv r-1+\overline{V}_\omega(\omega_c).
\end{equation}
It is clear from \eqref{lema-ctdaadd} that $\nu$ is obtained  by inserting a non-overlined part $t\eta$ into   $\omega$. So  $\omega$ satisfies the   condition U, which means that there exist $k-2$ consecutive parts in $(\overline{(t-1)\eta},(t+1)\eta]$, say
\[(t+1)\eta\geq \omega_i\geq \cdots \geq \omega_{i+k-3}>\overline{(t-1)\eta},\]
satisfying $\omega_i\leq \omega_{i+k-3}+\eta$ with strict inequality if $\omega_i$ is overlined, and \begin{equation}\label{construction-comb-3}
[\omega_{i}/\eta]+\cdots+[\omega_{i+k-3}/\eta]\equiv t+r-1+\overline{V}_{\omega}(\omega_{i})\pmod2.
\end{equation}
Now, \eqref{lema-ctdaa} yields that $\omega_{i+k-3}>\overline{(t-1)\eta}\geq \omega_{c+k-2}$, which implies that $i\leq c$.  Set $c=i+t$, where $t\geq0$. Then we have
\begin{equation}\label{construction-comb-3aa}
    (t+1)\eta\geq \omega_i\geq \cdots \geq \omega_{i+t-1}>\overline{t\eta},
\end{equation}
and
\begin{equation}\label{construction-comb-3bbb}
t\eta>\omega_{c+k-t-2}\geq \cdots \geq \omega_{c+k-3}\geq (t-1)
\eta.
\end{equation}
Combining \eqref{construction-comb-3aa} and \eqref{construction-comb-3bbb}, we obtain that
\[[\omega_{i}/\eta]+\cdots+[\omega_{i+t-1}/\eta]\equiv[\omega_{c+k-t-2}/\eta]+\cdots+[\omega_{c+k-3}/\eta]+\overline{V}_{\omega}(\omega_{i})-\overline{V}_{\omega}(\omega_{c})\pmod2,\]
which can be rewritten as
\begin{equation}\label{construction-comb-3aaee}
[\omega_{i}/\eta]+\cdots+[\omega_{i+k-3}/\eta]\equiv[\omega_{c}/\eta]+\cdots+[\omega_{c+k-3}/\eta]+\overline{V}_{\omega}(\omega_{i})-\overline{V}_{\omega}(\omega_{c})\pmod2.
\end{equation}
Substituting \eqref{lema-ctdaa} and \eqref{construction-comb-2}  into \eqref{construction-comb-3aaee} gives
\[
\begin{split}
&[\omega_{i}/\eta]+\cdots+[\omega_{i+k-3}/\eta]\\[5pt]
&\qquad =[\omega_{c}/\eta]+\cdots+[\omega_{c+k-3}/\eta]+[\omega_{c+k-2}/\eta]-(t-1)+\overline{V}_{\omega}(\omega_{i})-\overline{V}_{\omega}(\omega_{c})\\[5pt]
&\qquad \equiv r-1+\overline{V}_\omega(\omega_c)-(t-1)+\overline{V}_{\omega}(\omega_{i})-\overline{V}_{\omega}(\omega_{c})\\[5pt]
&\qquad \equiv t+r+\overline{V}_{\omega}(\omega_{i})\pmod{2},
\end{split}
\]
which contradicts \eqref{construction-comb-3}. Hence the condition (B) holds. Together with the condition (A), we conclude  that $\nu$ is an overpartition in $\overline{\mathcal{B}}_1(\alpha_1,\ldots,\alpha_\lambda;\eta,k,r)$.

{\noindent Condition (C).}  We consider the following two cases.

{\noindent Case 1.}  $\omega$ satisfies the   condition O  in Definition \ref{insertion}. In this circumstance, $\nu$ is obtained from $\omega$ by inserting $\overline{t\eta}$    as an overlined part. Obviously,   $ol(\nu)=\overline{t\eta}$ and  $\tilde{r}_1(\nu)\leq \overline{t\eta}$.

{\noindent Case 2.}  $\omega$ satisfies the   condition U  in Definition \ref{insertion}.  If so,
 $\nu$ is obtained from $\omega$ by inserting   ${t\eta}$ as a non-overlined part. Under the condition that $ol(\omega)<\overline{t\eta}$ and $\tilde{r}_1(\omega)\leq \overline{(t-1)\eta}$, we deduce that $ol(\nu)=ol(\omega)<\overline{t\eta}$ and $\tilde{r}_1(\nu)\leq {t\eta}$.

To prove that $\tilde{r}_1(\nu)>\overline{(t-1)\eta}$, it suffices to show that there is a $(k-1)$-band in $(\overline{(t-1)\eta}, {(t+1)\eta}]$. With the assumption that $\omega$ satisfies  the  condition U, we know that there exist $k-2$ consecutive parts   of $\omega$, say
\[
 (t+1)\eta\geq \omega_s\geq \cdots\geq \omega_{s+k-3}>\overline{(t-1)\eta},
\]
satisfying $\omega_s\leq \omega_{s+k-3}+\eta$ with strict inequality if $\omega_s$ is overlined, and
 \begin{equation}\label{even-moduli-over-0}
 [\omega_{s}/\eta]+\cdots+[\omega_{s+k-3}/\eta]\equiv t+r-1+\overline{V}_{\omega}(\omega_{s})\pmod2.
 \end{equation}
 Since $\tilde{r}_1(\omega)\leq \overline{(t-1)\eta}$, there are no $(k-1)$-bands of $\omega$ in $(\overline{(t-1)\eta}, (t+1)\eta]$. It follows that $\omega_{s-1}>\overline{t\eta}$ and $\omega_{s+k-2}<t\eta$. Assume that   $c$ is the smallest integer such that $\nu_{c}\leq t\eta$. Then
we have $s\leq c\leq s+k-2$,  $\nu_{l}=\omega_{l}$   for $s\leq l\leq c-1$,  and   $\nu_{l+1}=\omega_{l}$ for $c\leq l\leq s+k-3$, namely,
\begin{equation}\label{even-moduli-over-011}
  (t+1)\eta\geq \nu_s\geq\cdots\geq \nu_{c-1}> t\eta \geq \nu_{c+1}  \geq\cdots\geq \nu_{s+k-2}>\overline{(t-1)\eta}
\end{equation}
are the $k-1$ parts of $\nu$ such that      $\nu_s\leq \nu_{s+k-2}+\eta$ with strict inequality if $\nu_s$ is overlined.
Hence $\{\nu_{s+l}\}_{0\leq l\leq k-2}$ is a $(k-1)$-band of $\nu$. So  we arrive at  $\tilde{r}_1(\nu)>\overline{(t-1)\eta}$, and this proves that  the condition (C) is valid.

{\noindent Condition   (D).}  There are two cases.

{\noindent Case  1.}  $\omega$ satisfies the   condition O  in Definition \ref{insertion}. Then $\nu$ is obtained from $\omega$ by inserting $\overline{t\eta}$   as an overlined part. From the condition (C), we know that $\tilde{r}_1(\nu)\leq \overline{t\eta}$ in this event.  Assume that $\{\nu_{i+l}\}_{0\leq l\leq k-2}$   is a  $(k-1)$-band   of $\nu$. We aim to show that $\{\nu_{i+l}\}_{0\leq l\leq k-2}$ is even in $\nu$.     Since $\tilde{r}_1(\nu)\leq \overline{t\eta}$ and there is a part in $\{\nu_{i+l}\}_{0\leq l\leq k-2}$ marked with $k-1$ in $RG(\nu)$,  we get $\nu_{i+k-2}\leq \overline{t\eta}$. There are two  subcases.

{\noindent Subcase   1.1.}   $\nu_{i+k-2}\leq \overline{(t-1)\eta}$. In this case, the assumption that $\{\nu_{i+l}\}_{0\leq l\leq k-2}$   is a  $(k-1)$-band  implies that $\nu_i\leq \nu_{i+k-2} $ with strict inequality if $\nu_i$ is overlined, and so $\nu_i<\overline{t\eta}$.
Recall that  $\nu$ is obtained from $\omega$ by inserting $\overline{t\eta}$   as an overlined part,
we find that $\{\nu_{i+l}\}_{0\leq l\leq k-2}$ is also a  $(k-1)$-band of $\omega$ and  $\overline{V}_\omega(\nu_{i})=\overline{V}_\nu(\nu_{i})$. Since $\omega\in\overline{\mathcal{B}}^{\,<}_0(\alpha_1,\ldots,\alpha_\lambda;\eta,k,r)$, we see that $\{\nu_{i+l}\}_{0\leq l\leq k-2}$  is even in $\omega$, and so   $\{\nu_{i+l}\}_{0\leq l\leq k-2}$  is also even in $\nu$.

{\noindent Subcase   1.2.}   $\overline{(t-1)\eta}<\nu_{i+k-2}\leq \overline{t\eta}$. In this case, using the same reasoning as in the Subcase 1.1, we obtain that $\nu_i\leq (t+1)\eta$, so that
 \[(t+1)\eta\geq\nu_i\geq \nu_{i+1}\geq \cdots \geq \nu_{i+k-2}>\overline{(t-1)\eta}.\]
Utilizing the condition (B), we deduce that $\nu_{i-1}>\overline{t\eta}$ and $\nu_{i+k-1}<t\eta$. Hence  the  inserted  part $\overline{t\eta}$ in $\nu$ belongs to  $\{\nu_{i+l}\}_{0\leq l\leq k-2}$, so we may write $\nu_m=\overline{t\eta}$, where $i\leq m\leq i+k-2$. By the construction of $\nu$, we see that
$\omega_l=\nu_l$ for $i\leq l\leq m-1$,  and $\omega_l=\nu_{l+1}$  for $m\leq l\leq i+k-3$. This implies that  \[(t+1)\eta\geq\omega_i\geq\cdots\geq  \omega_{i+k-3}>\overline{(t-1)\eta}\]
  are the $k-2$ consecutive parts of $\omega$ such that $\omega_i\leq \omega_{i+k-3}+\eta$ with strict inequality provided $\omega_i$ is overlined, and $\overline{V}_{\nu}(\nu_{i})=\overline{V}_{\omega}(\omega_{i})+1$.
   Under the assumption  that $\omega$ satisfies the condition O, we find that
  \[[\omega_{i}/\eta]+\cdots+[\omega_{i+k-3}/\eta]\equiv t+r+\overline{V}_{\omega}(\omega_{i})\pmod2.\]
Therefore,
\[
\begin{split}
&[\nu_i/\eta]+\cdots+[\nu_{m}/\eta]+\cdots+[\nu_{i+k-2}/\eta]\\[5pt]
&\qquad =[\omega_i/\eta]+\cdots+[\omega_{i+k-3}/\eta]+t\\[5pt]
&\qquad \equiv t+r+\overline{V}_{\omega}(\omega_{i})+t\\[5pt]
&\qquad \equiv r-1+\overline{V}_{\nu}(\nu_{i})\pmod{2},
\end{split}
\]
which means that  $\{\nu_{i+l}\}_{0\leq l\leq k-2}$  is  even in $\nu$.

{\noindent Case 2.}  $\omega$ satisfies the   condition U  in Definition \ref{insertion}.  In this case,
 $\nu$ is obtained from $\omega$ by inserting   ${t\eta}$   as a non-overlined part. For any $(k-1)$-band $\{\nu_{i+l}\}_{0\leq l\leq k-2}$ of $\nu$, we wish to show that  $\{\nu_{i+l}\}_{0\leq l\leq k-2}$ is even in $\nu$. There are two cases.

{\noindent Subcase 2.1.} $\nu_i<t\eta$. By construction of $\nu$, we see that  $\{\nu_{i+l}\}_{0\leq l\leq k-2}$ is also a $(k-1)$-band of $\omega$ and  $\overline{V}_\omega(\nu_{i})=\overline{V}_\nu(\nu_{i})$. Using the same argument as
in Subcase  1.1, it can be shown that $\{\nu_{i+l}\}_{0\leq l\leq k-2}$ is even in $\nu$.

{\noindent Subcase 2.2.} $\nu_i\geq t\eta$.  Since $\{\nu_{i+l}\}_{0\leq l\leq k-2}$ is a $(k-1)$-band of $\nu$, we deduce that   $\nu_{i+k-2}\geq \nu_i-\eta\geq {(t-1)\eta}$ and there is a part $\nu_{i+l_i}$ ($0\leq l_i\leq k-2$) marked with $k-1$ in $RG(\nu)$. Using the condition (C), we find that $\overline{(t-1)\eta}<\tilde{r}_1(\nu)\leq {t\eta}$. It follows that $\tilde{r}_1(\nu)=\nu_{i+l_i}$. Hence $\{\nu_{i+l}\}_{0\leq l\leq k-2}$ is a $(k-1)$-band of $\nu$ including $\tilde{r}_1(\nu)$.  As in the proof of the condition (C), we see that   $\tilde{r}_1(\nu)$ is also in  the $(k-1)$-band  $\{\nu_{s+l}\}_{0\leq l\leq k-2}$ in \eqref{even-moduli-over-011}. Therefore, it follows from Proposition  \ref{parity k-1 sequence-over}  that   $\{\nu_{i+l}\}_{0\leq l\leq k-2}$  and $\{\nu_{s+l}\}_{0\leq l\leq k-2}$ have same parity. Hence we only need  to show that $\{\nu_{s+l}\}_{0\leq l\leq k-2}$ in \eqref{even-moduli-over-011} is even in $\nu$.

 Assume that $c$ is the smallest integer such that $\nu_c=t\eta$. As in the proof of the condition (C), we find that $\nu_c$ belongs to    the $(k-1)$-band  $\{\nu_{s+l}\}_{0\leq l\leq k-2}$ in \eqref{even-moduli-over-011} and $\overline{V}_\nu(\nu_{s})=\overline{V}_\omega(\omega_s)$. Thus,
\begin{equation}\label{insert-lem}
[\nu_{s}/\eta]+\cdots+[\nu_{s+k-2}/\eta]=t+[\omega_s/\eta]+\cdots+[\omega_{s+k-3}/\eta].
\end{equation}
Substituting \eqref{even-moduli-over-0} into \eqref{insert-lem} and using  $\overline{V}_\nu(\nu_{s})=\overline{V}_\omega(\omega_s)$,  we are led to
\begin{equation*}
[\nu_{s}/\eta]+\cdots+[\nu_{s+k-2}/\eta]\equiv r-1+\overline{V}_\nu(\nu_{s})\pmod{2},
\end{equation*}
which means that the $(k-1)$-band $\{\nu_{s+l}\}_{0\leq l\leq k-2}$ in \eqref{even-moduli-over-011} is even, and so $\{\nu_{i+l}\}_{0\leq l\leq k-2}$ is even in $\nu$.

In either case, we have shown that  any $(k-1)$-band  of $\nu$  is even, and thus the condition (D) is verified.

By now, we have shown that   $\nu$ is an overpartition
 in $\overline{\mathcal{B}}^{\,=}_0(\alpha_1,\ldots,\alpha_\lambda;\eta,k,r|t)$.  Clearly,  $|\nu|=|\omega|+t\eta$ and $\ell(\nu)=\ell(\omega)+1$. This completes the proof.  \qed

We are now in a position to give a proof of Theorem \ref{theorem4.4} with the aid of  Lemma \ref{division} and Lemma \ref{combinatio}.

\noindent{\it Proof of Theorem \ref{theorem4.4}.} Let $\nu\in\overline{\mathcal{B}}^{\,=}_0(\alpha_1,\ldots,\alpha_\lambda;\eta,k,r|t)$. Invoking  Lemma \ref{division}, we know that $D_{t}(\nu)\in\overline{\mathcal{B}}^{\,<}_0(\alpha_1,\ldots,\alpha_\lambda;\eta,k,r|t)$.
Setting $\omega=D_{t}(\nu)$, by Lemma \ref{combinatio}, we see that $C_{t}(\omega) \in \overline{\mathcal{B}}^{\,=}_0(\alpha_1,\ldots,\alpha_\lambda;\eta,k,r|t)$. It remains to show that $\nu=C_{t}(\omega)$. We consider the following two cases.

{\noindent Case 1.}   $ol(\nu)<\overline{t\eta}$ and  $\overline{(t-1)\eta}<\tilde{r}_1(\nu)\leq{t\eta}$. In this case,  $\omega$ is obtained from $\nu$ by removing a non-overlined part ${t\eta}$.
To prove that $\nu=C_{t}(\omega)$, it suffices to show that  $\omega$ satisfies    the condition U  in Definition \ref{insertion}.

Assume that $\tilde{r}_1(\nu)$ is the $r_1$-th part of $\nu=(\nu_1,\nu_2,\ldots,\nu_\ell)$, that is, $\nu_{r_1}=\tilde{r}_1(\nu)$. Then the $(k-1)$-band induced by $\tilde{r}_1(\nu)$ consists of
\[\nu_{r_1-k+2}\geq \nu_{r_1-k+3}\geq\cdots\geq\nu_{r_1},\]
where $\nu_{r_1-k+2}\leq \nu_{r_1}+\eta$ with strict inequality if $\nu_{r_1}$ is overlined. Since  $\nu$ is an overpartition in  $\overline{\mathcal{B}}^{\,=}_0(\alpha_1,\ldots,\alpha_\lambda;\eta,k,r|t)$,   we deduce that $\{\nu_{r_1-l}\}_{0\leq l\leq k-2}$ is even, namely,
  \begin{equation}\label{k-2}
[\nu_{r_1-k+2}/\eta]+[\nu_{r_1-k+3}/\eta]+\cdots+[\nu_{r_1}/\eta]\equiv r-1+\overline{V}_\nu(\nu_{r_1-k+2})\pmod2.
\end{equation}

Under the assumption  $\overline{(t-1)\eta}<\tilde{r}_1(\nu)\leq{t\eta}$, we see that $\nu_{r_1}=\tilde{r}_1(\nu)>\overline{(t-1)\eta}$ and $\nu_{r_1-k+2}\leq \nu_{r_1}+\eta\leq (t+1)\eta$, and thus
\[(t+1)\eta\geq\nu_{r_1-k+2}\geq \nu_{r_1-k+3}\geq\cdots\geq\nu_{r_1}>\overline{(t-1)\eta}.\]
It is clear from  Proposition \ref{daxiao-0-5} that  $\nu$ contains a non-overlined part $t\eta$. Assume that  $m$ is the smallest integer such that   $\nu_{m}=t\eta$. The precondition $\nu\in\overline{\mathcal{B}}^{\,=}_0(\alpha_1,\ldots,\alpha_\lambda;\eta,k,r|t)$ implies that $\nu_{r_1-k+1}>\overline{t\eta}$ and $\nu_{r_1+1}<t\eta$. Hence   $r_1-k+2\leq m\leq r_1$,
$\omega_l=\nu_l$ for $r_1-k+2\leq l\leq m-1$, and $\omega_l=\nu_{l+1}$ for $m\leq l\leq r_1-1$.  Consequently,  \[(t+1)\eta\geq\omega_{r_1-k+2}\geq\cdots\geq  \omega_{r_1-1}>\overline{(t-1)\eta}\]
  are the $k-2$ consecutive parts of $\omega$ such that $\omega_{r_1-k+2}\leq \omega_{r_1-1}+\eta$ with strict inequality provided $\omega_{r_1-k+2}$ is overlined.  By the construction of $\omega$, we deduce that  $\overline{V}_\omega(\omega_{r_1-k+2})=\overline{V}_\nu(\nu_{r_1-k+2})$. Combining with \eqref{k-2}, we get
\[
 \begin{split}
 &[\omega_{r_1-k+2}/\eta]+[\omega_{r_1-k+1}/\eta]+\cdots+[\omega_{r_1-1}/\eta]\\
 &\qquad =[\nu_{r_1-k+2}/\eta]+[\nu_{r_1-k+1}/\eta]+\cdots+[\nu_{r_1}/\eta]-t\\
 &\qquad \equiv r-1+\overline{V}_\nu(\nu_{r_1-k+2})-t\\
 &\qquad \equiv t+r-1+\overline{V}_\omega(\omega_{r_1-k+2})\pmod2.
 \end{split}
\]
This implies that  $\omega$ satisfies    the condition U  in Definition \ref{insertion}, and so $\nu=C_t(\omega)$. Hence we conclude that  $C_t(D_t(\nu))=\nu$ for $\nu\in\overline{\mathcal{B}}^{\,=}_0(\alpha_1,\ldots,\alpha_\lambda;\eta,k,r|t)$.

{\noindent Case 2.}  $ol(\nu)=\overline{t\eta}$ and  $\tilde{r}_1(\nu)\leq\overline{t\eta}$. In this regard,  $\omega$ is obtained  from $\nu$ by removing $\overline{t\eta}$. To prove that $\nu=C_{t}(\omega)$, it is enough  to show that  $\omega$ satisfies    the condition O  in Definition \ref{insertion}.  Suppose to the contrary that  $\omega$ satisfies the   condition U  in Definition \ref{insertion}, that is, there exist $k-2$ consecutive parts of $\omega$, say
 \[(t+1)\eta\geq \omega_i\geq \omega_{i+1}\geq \cdots \geq \omega_{i+k-3}>\overline{(t-1)\eta},\]
such that $\omega_i\leq \omega_{i+k-3}+\eta$ with strict inequality if $\omega_i$ is overlined,  and
\begin{equation}\label{divcombin-temp}
 [\omega_{i}/\eta]+\cdots+[\omega_{i+k-3}/\eta]\equiv t+r-1+\overline{V}_{\omega}(\omega_{i})\pmod2.
\end{equation}
 Assume that $\overline{t\eta}$ is the $m$-th part $\nu_m$ of $\nu$. Since $\omega \in\overline{\mathcal{B}}^{\,<}_0(\alpha_1,\ldots,\alpha_\lambda;\eta,k,r|t)$, we see that $\tilde{r}_1(\omega)\leq \overline{(t-1)\eta}$, and so there are no $(k-1)$-bands in $(\overline{(t-1)\eta}, (t+1)\eta]$. It follows that $\omega_{i-1}>\overline{t\eta}$ and $\omega_{i+k-2}<t\eta$, which implies that $i\leq m\leq i+k-2$,   $\omega_{l}=\nu_{l}>\overline{t\eta}$ for $i\leq l< m$, and   $\omega_{l}=\nu_{l+1}\leq t\eta$ for $m\leq l\leq i+k-3$.
Thus,
\[(t+1)\eta\geq \nu_i\geq\cdots  \geq \nu_{i+k-2}>\overline{(t-1)\eta},\]
 where    $\nu_i\leq \nu_{i+k-2}+\eta$ with strict inequality if $\nu_i$ is overlined.
 In other words,  $\{\nu_{i+l}\}_{0\leq l\leq k-2}$ is a $(k-1)$-band of $\nu$. Moreover, we get $\overline{V}_\nu(\nu_i)=\overline{V}_\omega(\omega_i)+1$.
The precondition that $\nu$ is an overpartition in  $\overline{\mathcal{B}}^{\,=}_0(\alpha_1,\ldots,\alpha_\lambda;\eta,k,r|t)$ implies that $\{\nu_{i+l}\}_{0\leq l\leq k-2}$ is even, and so
\[
 \begin{split}
 &[\omega_i/\eta]+\cdots+[\omega_{i+k-3}/\eta]\\
 &\qquad =[\nu_i/\eta]+\cdots+[\nu_{i+k-2}/\eta]-t\\
 &\qquad \equiv r-1+\overline{V}_\nu(\nu_i)-t\\
 &\qquad \equiv t+r+\overline{V}_\omega(\omega_i)\pmod2,
 \end{split}
\]
which contradicts   \eqref{divcombin-temp}. Hence   $\omega$ satisfies the   condition O  in Definition \ref{insertion}, and so $\nu=C_t(\omega)$. This proves that $C_t(D_t(\nu))=\nu$ for $\nu\in\overline{\mathcal{B}}^{\,=}_0(\alpha_1,\ldots,\alpha_\lambda;\eta,k,r|t)$.

Conversely, let   $\omega\in\overline{\mathcal{B}}^{\,<}_0(\alpha_1,\ldots,\alpha_\lambda;\eta,k,r|t)$. By Lemma \ref{combinatio}, we find that $C_t(\omega)$ belongs to  $\overline{\mathcal{B}}^{\,=}_0(\alpha_1,\ldots,\alpha_\lambda;\eta,k,r|t)$.
By the definitions of $C_t$ and $D_t$, we deduce that $D_t(C_t(\omega))=\omega$. This completes the proof. \qed

The following proposition provides a criterion to determine whether an overpartition in $\overline{\mathcal{B}}^{\,=}_0(\alpha_1,\ldots,\alpha_\lambda;\eta,k,r|t)$ is also an overpartition in $\overline{\mathcal{B}}^{\,<} _0(\alpha_1,\ldots,\alpha_\lambda;\eta,k,r|t')$.

\begin{prop}\label{daxiao-0-4}
For $t\geq 1$, let $\nu$ be an overpartition in $\overline{\mathcal{B}}^{\,=}_0(\alpha_1,\ldots,\alpha_\lambda;\eta,k,r|t)$. Then $\nu$ is an overpartition in $\overline{\mathcal{B}}^{\,<} _0(\alpha_1,\ldots,\alpha_\lambda;\eta,k,r|t')$ if and only if $t<t'$.
\end{prop}

\pf  By definition, we see that $\nu$ is  an overpartition in $\overline{\mathcal{B}}^{\,=}_0(\alpha_1,\ldots,\alpha_\lambda;\eta,k,r|t)$ if and only if
$\nu$ is  an overpartition in $\overline{\mathcal{B}}_0(\alpha_1,\ldots,\alpha_\lambda;\eta,k,r)$ such that
\begin{equation}\label{daxiao-eqn-0-4}
 \max\{\lceil |ol(\nu)|/\eta\rceil,\lceil |\tilde{r}_1(\nu)|/\eta\rceil\}=t,
    \end{equation}
    where $|\cdot|$  signified the value of a part regardless of overline, and  $\lceil x\rceil$ denotes the smallest integer greater than or equal to $x$.

On the other hand,   $\nu$ is  an overpartition in $\overline{\mathcal{B}}^{\, <}_0(\alpha_1,\ldots,\alpha_\lambda;\eta,k,r|t')$ if and only if $\nu$ is  an overpartition in $\overline{\mathcal{B}}_0(\alpha_1,\ldots,\alpha_\lambda;\eta,k,r)$ such that
\begin{equation}\label{daxiao-eqn-0-4b}
 \max\{\lceil |ol(\nu)|/\eta\rceil,\lceil |\tilde{r}_1(\nu)|/\eta\rceil\}\leq t'-1.
    \end{equation}
 Combining  \eqref{daxiao-eqn-0-4} and \eqref{daxiao-eqn-0-4b} completes the proof.  \qed

\subsection{Proof of Theorem \ref{lem-b-1-over}}

In this subsection, we   demonstrate that   Theorem \ref{lem-b-1-over} can be
justified by  repeatedly  using   the $(k-1)$-reduction  and   the $(k-1)$-augmentation operations.

{\it \noindent Proof of Theorem \ref{lem-b-1-over}.} Let $\nu$ be an overpartition in $\mathcal{\overline{B}}_0(\alpha_1,\ldots,\alpha_\lambda;\eta,k,r)$. We  wish to construct a pair of overpartitions  $\Theta(\nu)=(\zeta,\omega)$ in  $\mathcal{D}_\eta \times {\mathcal{B}}_{1}(\alpha_1,\ldots,\alpha_\lambda;\eta,k-1,
 r-\delta_{r,k})$ such that $|\nu|=|\zeta|+|\omega|$ and $\ell(\nu)=\ell(\zeta)+\ell(\omega)$.  We consider the following two cases:

 Case 1: There are no     $(k-1)$-marked parts in $RG(\nu)$ and  there are no   overlined parts divisible by $\eta$ in $\nu$. Then set  $\zeta=\emptyset$ and   $\omega=\nu$.  By definition,  we see that  $\omega$ is an overpartition in ${\mathcal{B}}_{1}(\alpha_1,\ldots,\alpha_\lambda;\eta,k-1,
 r-\delta_{r,k})$. Moreover, $|\nu|=|\zeta|+|\omega|$ and $\ell(\nu)=\ell(\zeta)+\ell(\omega)$.

 Case 2: There exists a $(k-1)$-marked part  in $RG(\nu)$ or an overlined part  divisible by $\eta$ in $\nu$.   Set $b=0$, $\nu^{(0)} =\nu$, $\zeta^{(0)}=\emptyset$, and execute  the following procedure. Denote the intermediate pairs by $(\zeta^{(0)},\nu^{(0)}),(\zeta^{(1)},\nu^{(1)})$, and so on.
 \begin{itemize}

\item[(A)]  Set
     \[t_{b+1}=\max\{\lceil |ol(\nu^{(b)})|/\eta\rceil,\lceil |\tilde{r}_1(\nu^{(b)})|/\eta\rceil\}.\]
Since $\tilde{r}_1(\nu^{(b)})\geq \overline{\alpha_1}$ or $ol(\nu^{(b)})\geq \overline{\eta}$, we find that $t_{b+1}\geq 1$ and
\[\nu^{(b)}\in \overline{\mathcal{B}}^{\,=}_0(\alpha_1,\ldots,\alpha_\lambda;\eta,k,r|t_{b+1}).\]

Applying the $(k-1)$-reduction $D_{t_{b+1}}$ to $\nu^{(b)}$, we get
 \[\nu^{(b+1)}=D_{t_{b+1}}(\nu^{(b)}).\]
In view of Lemma \ref{division}, we deduce that
$\nu^{(b+1)}\in \overline{\mathcal{B}}^{\,<}_0(\alpha_1,\ldots,\alpha_\lambda;\eta,k,r|t_{b+1}),$
 \begin{equation}\label{theorem4.1weig}
  |\nu^{(b+1)}|=|\nu^{(b)}|-\eta t_{b+1},
 \end{equation}
 and
  \begin{equation}\label{theorem4.1len}
     \ell(\nu^{(b+1)})=\ell(\nu^{(b)})-1.
 \end{equation}
Then, insert $\eta t_{b+1}$  into $\zeta^{(b)}$ as a part to get $\zeta^{(b+1)}$.

\item[(B)] Replace $b$ by $b+1$. If there are no  $(k-1)$-marked parts in $RG(\nu^{(b)})$ and there are no  overlined parts divisible by $\eta$ in $\nu^{(b)}$, then we are done. Otherwise,
go back to (A).

 \end{itemize}

 Using Proposition \ref{daxiao-0-4}, we obtain that
\begin{equation}\label{theorem4.1rec}
 t_{b+1}>t_{b+2}\geq 1,
\end{equation}
for $b\geq 0$, which means  that the above procedure terminates. Assume that it terminates with  $b=c$, that is, there are no  $(k-1)$-marked parts in $RG(\nu^{(c)})$ and there are no  overlined parts divisible by $\eta$ in $\nu^{(c)}$. Set
    \[\omega=\nu^{(c)}\quad \text{and} \quad \zeta=\zeta^{(c)}=(\eta t_1,\ldots,\eta t_c).\]
Since there are no  $(k-1)$-marked parts in $RG(\nu^{(c)})$ and there are no  overlined parts divisible by $\eta$ in $\nu^{(c)}$, we conclude that  $\omega=\nu^{(c)} \in {\mathcal{B}}_{1}(\alpha_1,\ldots,\alpha_\lambda;\eta,k-1,
 r-\delta_{r,k})$.
In light of \eqref{theorem4.1rec},  we find that $ \zeta \in \mathcal{D}_\eta$.   Moreover, it is clear from \eqref{theorem4.1weig} and \eqref{theorem4.1len} that  $|\nu|=|\omega|+|\zeta|$ and $\ell(\nu)=\ell(\omega)+\ell(\zeta)$. Hence  $\Theta$ is the desired map from  $\mathcal{\overline{B}}_0(\alpha_1,\ldots,\alpha_\lambda;\eta,k,r)$ to $\mathcal{D}_\eta\times{\mathcal{B}}_{1}(\alpha_1,\ldots,\alpha_\lambda;\eta,k-1,r-\delta_{r,k})$.

To prove that $\Theta$ is a bijection, we  define a map $\Lambda$ from  $\mathcal{D}_\eta\times{\mathcal{B}}_{1}(\alpha_1,\ldots,\alpha_\lambda;\eta,k-1,r-\delta_{r,k})$ to $\mathcal{\overline{B}}_0(\alpha_1,\ldots,\alpha_\lambda;\eta,k,r)$ and intend to
show that it is the inverse map of $\Theta$. Given  an overpartition  $\omega$  in  ${\mathcal{B}}_{1}(\alpha_1,\ldots,\alpha_\lambda;\eta,k-1,
 r-\delta_{r,k})$ and  a partition $\zeta$ in   $\mathcal{D}_\eta$, we shall construct  an overpartition $\nu \in \mathcal{\overline{B}}_0(\alpha_1,\ldots,\alpha_\lambda;\eta,k,r)$  such that $|\nu|=|\zeta|+|\omega|$ and $\ell(\nu)=\ell(\zeta)+\ell(\omega)$.  There are  two cases.

 Case 1:  $\zeta=\emptyset$. Then set $\nu=\omega$. Clearly, $\nu\in\mathcal{\overline{B}}_0(\alpha_1,\ldots,\alpha_\lambda;\eta,k,r)$ since  there are no  $(k-1)$-bands  in $\omega$. Moreover, $|\nu|=|\zeta|+|\omega|$ and $\ell(\nu)=\ell(\zeta)+\ell(\omega)$.

 Case 2: $\zeta\neq\emptyset$. Assume that $\zeta=(\eta t_1,\eta t_2,\ldots,\eta t_c)$, where $ t_1> t_2>\cdots> t_c\geq 1$. Starting with $\omega$, apply the $(k-1)$-augmentation repeatedly    to get $\nu$. Denote the intermediate overpartitions by $\omega^{(0)},\ldots,\omega^{(c)}$ with $\omega^{(0)}=\omega$ and  $\omega^{(c)}=\nu$. Since   there are no  $(k-1)$-marked parts in $RG(\omega)$ and  there are no   overlined parts divisible by $\eta$ in $\omega$,  we have $\tilde{r}_1(\omega)=-\infty$ and $ol(\omega)=\overline{0}$, which yields
 $\omega^{(0)}=\omega\in\overline{\mathcal{B}}^{\,<}_0
    (\alpha_1,\ldots,\alpha_\lambda;\eta,k,r| t_{c}).$

 Set  $b=0$, and execute  the following procedure.

\begin{itemize}

\item[(A)] Set
    \[\omega^{(b+1)}=C_{ t_{c-b}}(\omega^{(b)}).\]
Since
\[\omega^{(b)} \in \overline{\mathcal{B}}^{\, <}_0    (\alpha_1,\ldots,\alpha_\lambda;\eta,k,r|
 t_{c-b}),\]
in light of Lemma \ref{combinatio}, we see that $\omega^{(b+1)}\in \overline{\mathcal{B}}^{\,=}_0(\alpha_1,\ldots,\alpha_\lambda;\eta,k,r| t_{c-b}),$

 \begin{equation}\label{insertaa-1}
 |\omega^{(b+1)}|=|\omega^{(b)}|+\eta  t_{c-b},
 \end{equation}
 and
 \begin{equation}\label{insertaa-1t}
\ell(\omega^{(b+1)})=\ell(\omega^{(b)})+1. \end{equation}

   \item[(B)] Replace $b$ by $b+1$. If $b=c$, then we are done. Otherwise,  since $ t_{c-b}> t_{c-b+1}$, it follows from Proposition \ref{daxiao-0-4}  that
  \[\omega^{(b)}\in \overline{\mathcal{B}}^{\,<}_0(\alpha_1,\ldots,\alpha_\lambda;\eta,k,r|
   t_{c-b}).\]
  Go back to (A).
     \end{itemize}
The above procedure generates an overpartition $\nu=\omega^{(c)}\in  \overline{\mathcal{B}}^{\,=}_0(\alpha_1,\ldots,\alpha_\lambda;\eta,k,r|
 t_{1})$, and so   $\nu$ is an overpartition in $\mathcal{\overline{B}}_0(\alpha_1,\ldots,\alpha_\lambda;\eta,k,r)$.  It is evident from \eqref{insertaa-1} and  \eqref{insertaa-1t} that
\begin{equation*}
     |\nu|=|\omega^{(c)}|=|\omega^{(0)}|+\eta t_{c}+\cdots+\eta t_{1}=|\omega|+|\zeta|,
    \end{equation*}
    and
 \begin{equation*}
     \ell(\nu)=\ell(\omega^{(c)})=\ell(\omega^{(0)})+c=\ell(\omega)+\ell(\zeta).
    \end{equation*}
Therefore,  $\Lambda$ is a map from  $\mathcal{D}_\eta\times{\mathcal{B}}_{1}(\alpha_1,\ldots,\alpha_\lambda;\eta,k-1,r-\delta_{r,k})$ to $\mathcal{\overline{B}}_0(\alpha_1,\ldots,\alpha_\lambda;\eta,k,r)$. By  Theorem \ref{theorem4.4},  we obtain that $\Lambda(\Theta(\nu))=\nu $ for $\nu \in \overline{\mathcal{B}}_0(\alpha_1,\ldots,\alpha_\lambda;\eta,k,r)$ and $\Theta(\Lambda(\zeta,\omega))=(\zeta,\omega) $ for $(\zeta,\omega)\in\mathcal{D}_\eta\times{\mathcal{B}}_{1}(\alpha_1,\ldots,\alpha_\lambda;\eta,k-1,r-\delta_{r,k})$. Hence $\Theta$ is a bijection   between $\mathcal{\overline{B}}_0(\alpha_1,\ldots,\alpha_\lambda;\eta,k,r)$ and $\mathcal{D}_\eta\times{\mathcal{B}}_{1}(\alpha_1,\ldots,\alpha_\lambda;\eta,k-1,r-\delta_{r,k})$. This completes the proof.   \qed

 \subsection{An example}

 We conclude this section with an example for  the bijection $\Theta$ in Theorem \ref{lem-b-1-over}. Let \[\nu=(\overline{50},\overline{30},\overline{23},20,20,\overline{10}
,\overline{7},\overline{3})\]
 be an overpartition in $\overline{\mathcal{B}}_0(3,7;10,4,3)$. We have
 \[RG(\nu)=(\overline{50}_1,\overline{30}_1,\overline{23}_2,{20}_1,{ {20}_3},\overline{10}_2
,\overline{7}_1,{ \overline{3}_3}).\]
The pair of overpartitions  $\Theta(\nu)=(\zeta,\omega)$ is obtained by successively applying the $(k-1)$-reduction to $\nu$. The detailed process is given below.

 \begin{itemize}
     \item   Set  $\nu^{(0)}=\nu$ and $\zeta^{(0)}=\emptyset$. Note that $ol(\nu^{(0)})=\overline{50}$ and $\tilde{r}_1(\nu^{(0)})=20$. Let
  \[t_1=\text{max}\{\lceil |ol(\nu^{(0)})|/{10}\rceil,\lceil |\tilde{r}_1(\nu^{(0)})|/{10}\rceil\}=5.\]
Now,  $\nu^{(0)}\in\overline{\mathcal{B}}^{\,=}_0(3,7;10,4,3|5)$.    Apply the $3$-reduction to $\nu^{(0)}$ to get $\nu^{(1)}$, namely, $\nu^{(1)}$ is obtained from $\nu^{(0)}$ by  removing $\overline{50}$. We get
\[RG(\nu^{(1)})=(\overline{30}_1,\overline{23}_2,{20}_1,{ {20}_3},\overline{10}_2
,\overline{7}_1,{ \overline{3}_3}).\]
 Setting $\zeta^{(1)}=(50)$ and using Lemma \ref{division}, we obtain that $\nu^{(1)}\in\overline{\mathcal{B}}^{\,<}_0(3,7;10,4,3|5)$.

   \item Since $ol(\nu^{(1)})=\overline{30}$ and $\tilde{r}_1(\nu^{(1)})=20$, we have
\[t_2=\text{max}\{ \lceil |ol(\nu^{(1)})|/{10}\rceil, \lceil |\tilde{r}_1(\nu^{(1)})|/{10}\rceil\}=3,\]
whence  $\nu^{(1)}\in\overline{\mathcal{B}}^{\,=}_0(3,7;10,4,3|3)$. Removing $\overline{30}$ from $\nu^{(1)}$, we  get $\nu^{(2)}$ and
\[RG(\nu^{(2)})=(\overline{23}_1,{20}_2,{ {20}_3},\overline{10}_1,\overline{7}_2,{ \overline{3}_3}).\]
 Setting $\zeta^{(2)}=(50,30)$ and using Lemma \ref{division}, we obtain that $\nu^{(2)}\in \overline{\mathcal{B}}^{\,<}_0(3,7;10,4,3|3)$.

\item  Since $ol(\nu^{(2)})=\overline{10}$ and $\tilde{r}_1(\nu^{(2)})=20$,  we have
\[t_3=\text{max}\{\lceil |ol(\nu^{(2)})|/{10}\rceil, \lceil |\tilde{r}_1(\nu^{(2)})|/{10}\rceil\}=2,\]
whence $\nu^{(2)}\in\overline{\mathcal{B}}^{\,=}_0(3,7;10,4,3|2)$. Removing a non-overlined part ${20}$ from $\nu^{(2)}$, we get $\nu^{(3)}$ and
\[RG(\nu^{(3)})=(\overline{23}_1,{20}_2,\overline{10}_1,\overline{7}_2,{ \overline{3}_3}).\]
 Setting $\zeta^{(3)}=(50,30,20)$ and using  Lemma \ref{division}, we obtain that $\nu^{(3)}\in \overline{\mathcal{B}}^{\,<}_0(3,7;10,4,3|2)$.

 \item  Since  $ol(\nu^{(3)})=\overline{10}$ and $\tilde{r}_1(\nu^{(3)})=\overline{3}$,  we have \[t_4=\text{max}\{\lceil |ol(\nu^{(3)})|/{10}\rceil, \lceil |\tilde{r}_1(\nu^{(3)})|/{10}\rceil\}=1,\]
whence $\nu^{(3)}\in\overline{\mathcal{B}}^{\,=}_0(3,7;10,4,3|1)$.  Removing  $\overline{10}$ from $\nu^{(3)}$, we get $\nu^{(4)}$ and
 \begin{equation*}
 RG(\nu^{(4)})=(\overline{23}_1,{20}_2,\overline{7}_1,\overline{3}_2).
 \end{equation*}
 Setting $\zeta^{(4)}=(50,30,20,10)$  and using Lemma \ref{division}, we obtain that $\nu^{(4)}\in \overline{\mathcal{B}}^{\,<}_0(3,7;10,4,\break$
 $3|1)$. Eventually,  there are no $3$-marked parts in $RG(\nu^{(4)})$ and there are no overlined parts divisible by $10$ in $\nu^{(4)}$.
 \end{itemize}
We now get  a  pair of partitions $(\zeta,\omega)$ with
\begin{equation}\label{omega-reverse}
 \zeta=\zeta^{(4)}=(50,30,20,10)\quad  \text{and} \quad  \omega=\nu^{(4)}=(\overline{23},{20},\overline{7},\overline{3})
\end{equation}
such that    $(\zeta, \omega) \in \mathcal{D}_{10} \times \overline{\mathcal{B}}_1(3,7;10,3,3)$, $|\nu|=|\zeta|+|\omega|$ and $\ell(\nu)=\ell(\omega)+\ell(\zeta)$.

Conversely, given   $(\zeta, \omega) \in \mathcal{D}_{10} \times \overline{\mathcal{B}}_1(3,7;10,3,3)$ as in \eqref{omega-reverse}, we may recover the overpartition $\nu$ by successively applying the $3$-augmentation
operation. More precisely, the reverse process goes as follows.

\begin{itemize}
    \item Insert $10$ into $\omega^{(0)}=\omega$ to get $\omega^{(1)}$.

     Since there are no $3$-marked parts in $RG(\omega^{(0)})$ and there are no overlined parts divisible by $10$ in $\omega^{(0)}$,  we have   $\tilde{r}_1(\omega^{(0)})=-\infty$ and $ol(\omega^{(0)})=\overline{0}$, which implies that $\omega^{(0)}\in\overline{\mathcal{B}}^{\,<}_0(3,7;10,4,
    3|1)$. Notice that  $\omega^{(0)}$  satisfies the   condition O in Definition \ref{insertion}. Then insert $\overline{10}$ into $\omega^{(0)}$  as an overlined part to  get
     \[\omega^{(1)}=C_{1}(\omega^{(0)})=(\overline{23},20,{\overline{10},\overline{7},\overline{3})}.\]
     Using Lemma \ref{combinatio}, we obtain that $\omega^{(1)}\in\overline{\mathcal{B}}^{\,=}_0(3,7;10,4,
    3|1)$.

     \item   Insert $20$ into $\omega^{(1)}$ to get $\omega^{(2)}$.     By Proposition \ref{daxiao-0-4}, we find that $\omega^{(1)}\in\overline{\mathcal{B}}^{\,<}_0(3,7;10,4,
    3|2)$.  Since  $\omega^{(1)}$  satisfies the   condition U in Definition \ref{insertion},   inserting ${20}$  into $\omega^{(1)}$ as a non-overlined part gives
     \[\omega^{(2)}=C_{2}(\omega^{(1)})=(\overline{23},20,20,{\overline{10},\overline{7},\overline{3})}.\]
In light of Lemma \ref{combinatio}, we deduce that $\omega^{(2)}\in\overline{\mathcal{B}}^{\,=}_0(3,7;10,4,3|2)$.

  \item   Insert $30$ into $\omega^{(2)}$ to get $\omega^{(3)}$.
By Proposition \ref{daxiao-0-4}, we find  that $\omega^{(2)}\in\overline{\mathcal{B}}^{\,<}_0(3,7;10,4,
    3|3)$.   Notice that  $\omega^{(2)}$  satisfies the   condition O in Definition \ref{insertion}. Then insert $\overline{30}$ into $\omega^{(2)}$ as an overlined part  to  get
     \[\omega^{(3)}=C_{3}(\omega^{(2)})=(\overline{30},\overline{23},20,20,{\overline{10},\overline{7},\overline{3})}.\]
     Using Lemma \ref{combinatio}, we obtain that $\omega^{(3)}\in\overline{\mathcal{B}}^{\,=}_0(3,7;10,4,
    3|3)$.

    \item Finally, insert $50$ into $\omega^{(3)}$ to get $\omega^{(4)}$.  By Proposition \ref{daxiao-0-4}, we find that $\omega^{(3)}\in\overline{\mathcal{B}}^{\,<}_0(3,7;10,4,
    3|5)$. Notice that  $\omega^{(3)}$  satisfies the   condition O in Definition \ref{insertion}. Then insert $\overline{50}$ into $\omega^{(3)}$  as an overlined part to  get
    \[\omega^{(4)}=C_{5}(\omega^{(3)})=(\overline{50},\overline{30},\overline{23},20,20,{\overline{10},\overline{7},\overline{3})}.\]
 Using Lemma \ref{combinatio}, we obtain that $\omega^{(4)}\in\overline{\mathcal{B}}^{\,=}_0(3,7;10,4,
    3|5)$.
\end{itemize}

Set $\nu=\omega^{(4)}$. Then   $\nu$ is an overpartition in $\overline{\mathcal{B}}_0(3,7;10,4,3)$ such that $|\nu|=|\omega|+|\zeta|$ and $\ell(\nu)=\ell(\omega)+\ell(\zeta)$.

\section{Proof of Theorem \ref{G-B-O-1}}

In this section, we will give a proof of Theorem \ref{G-B-O-1} by using Bailey pairs. It remains to be a question to find a combinatorial proof of this fact.     For historical perspectives and recent advances on Bailey pairs, we refer   to  Agarwal, Andrews and Bressoud \cite{Agarwal-Andrews-Bressoud}, Andrews \cite{Andrews-1986, Andrews-2000}, Bressoud, Ismail and Stanton \cite{Bressoud-Ismail-Stanton-2000}, Lovejoy \cite{Lovejoy-2004b}, Paule \cite{Paule-1987}, Warnaar \cite{Warnaar-2001},  to name of few.  A pair of sequences $(\alpha_n(a,q),\beta_n(a,q))$ is said to be a Bailey pair relative to $(a,q)$  if  for  $n\geq 0,$
\begin{equation*}\label{bailey pair}
\beta_n(a,q)=\sum_{r=0}^n\frac{\alpha_r(a,q)}{(q;q)_{n-r}(aq;q)_{n+r}}.
\end{equation*}

The following  formulation of  Bailey's lemma was given by Andrews  \cite{Andrews-1984, Andrews-1986}.

\begin{thm}[Bailey's lemma]\label{BL}
If $(\alpha_n(a,q),\beta_n(a,q))$ is  a Bailey pair  relative to  $(a,q)$,
then $(\alpha_n'(a,q),\beta_n'(a,q))$ is also a Bailey pair relative to  $(a,q)$, where
\begin{equation*}\label{BL-eq}
\begin{split}
\alpha_n'(a,q)&=\frac{(\rho_1;q)_n(\rho_2;q)_n}{(aq/\rho_1;q)_n(aq/\rho_2;q)_n}
\left(\frac{aq}
{\rho_1\rho_2}\right)^n\alpha_n(a,q),\\
\beta_n'(a,q)&=
\sum_{j=0}^{n}\frac{(\rho_1;q)_j(\rho_2;q)_j
(aq/\rho_1\rho_2;q)_{n-j}}
{(aq/\rho_1;q)_n(aq/\rho_2;q)_n(q;q)_{n-j}}
\left(\frac{aq}{\rho_1\rho_2}\right)^j\beta_j(a,q).
\end{split}
\end{equation*}
\end{thm}
When  $\rho_1,\rho_2\rightarrow \infty$, Bailey's lemma reduces to the  following form, which has been used by Andrews \cite{Andrews-1984}  to derive  the Andrews-Gordon identity \eqref{R-R-A1}    when $r=1$ or $r=k$.

 \begin{lem} \label{BL-1}
If $(\alpha_n(a,q),\beta_n(a,q))$ is  a Bailey pair  relative to  $(a,q)$,
then $(\alpha_n'(a,q),\beta_n'(a,q))$ is also a Bailey pair relative to  $(a,q)$, where
\begin{equation*}\label{BL-1-eq}
\begin{split}
\alpha_n'(a,q)&=a^nq^{n^2}\alpha_n(a,q),\\
\beta_n'(a,q)&=
\sum_{j=0}^{n}\frac{a^jq^{j^2}}
{(q;q)_{n-j}}
\beta_j(a,q).
\end{split}
\end{equation*}
\end{lem}

 Agarwal, Andrews and Bressoud \cite{Agarwal-Andrews-Bressoud} developed the technique of the Bailey lattice to establish the Andrews-Gordon identity \eqref{R-R-A1}  in  general for $1\leq r\leq k$.  Bressoud, Ismail and Stanton \cite{Bressoud-Ismail-Stanton-2000} found an alternative proof of  the Andrews-Gordon identity \eqref{R-R-A1}   in the general case by successively  using Bailey's lemma and the following proposition.

\begin{prop}\label{Bl-Bp}{\rm  \cite[Proposition 4.1]{Bressoud-Ismail-Stanton-2000}}
If  $(\alpha_n(1,q),\beta_n(1,q))$ is a Bailey pair relative to $(1,q)$, where
\[\alpha_n(1,q)=\left\{
             \begin{array}{ll}
               1, & \hbox{if $n=0$}, \\[5pt]
               (-1)^nq^{An^2}(q^{(A-1)n}+q^{-(A-1)n}), & \hbox{if $n\geq 1$,}
             \end{array}
           \right.
\]
then $(\alpha_n'(1,q),\beta_n'(1,q))$ is also a Bailey pair relative to $(1,q)$, where
\begin{align*}
 \alpha_n'(1,q)&=\left\{
             \begin{array}{ll}
               1, & \hbox{if $n=0$}, \\[3pt]
               (-1)^nq^{An^2}(q^{An}+q^{-An}), & \hbox{if $n\geq 1$,}
             \end{array}
           \right. \\[10pt]
       &\hskip -5cm \text{and for } n\geq 0,  \\[2pt]
 \beta_n'(1,q)&=q^n\beta_n(1,q).
\end{align*}
\end{prop}

To prove  Theorem \ref{G-B-O-1}, we also need the following proposition in \cite{He-Ji-Wang-Zhao-2019} and a limiting case of an identity  of Andrews \cite{Andrews-1984}.

\begin{prop}\label{cbl4}
If  $(\alpha_n(1,q),\beta_n(1,q))$ is a Bailey pair relative to $(1,q)$, where
\[\alpha_n(1,q)=\left\{
             \begin{array}{ll}
               1, & \hbox{if $n=0$}, \\[5pt]
               (-1)^nq^{An^2}(q^{(A-1)n}+q^{-(A-1)n}), & \hbox{if $n\geq 1$,}
             \end{array}
           \right.
\]
then $(\alpha_n'(1,q),\beta_n'(1,q))$ is also a Bailey pair relative to $(1,q)$, where
\begin{align*}
 \alpha_n'(1,q)&=\left\{
             \begin{array}{ll}
               1, & \hbox{if $n=0$}, \\[3pt]
               (-1)^nq^{An^2}(q^{(A-1)n}+q^{-An})(1+q^n)/2, & \hbox{if $n\geq 1$,}
             \end{array}
           \right. \\[10pt]
            &\hskip -4cm \text{and for } n\geq 0,  \\[2pt]
 \beta_n'(1,q)&=\beta_n(1,q)(1+q^n)/2.
\end{align*}
\end{prop}

\begin{thm}[Andrews]\label{Multi} If $(\alpha_n(1,q),\beta_n(1,q))$ is a Bailey pair relative to $(1,q)$, then for $N\geq 0$,
\begin{align}\nonumber
&\sum_{n \geq 0} \frac{\left(b_{1};q\right)_{n}\left(c_{1};q\right)_{n} \cdots\left(b_{k};q\right)_{n}\left(c_{k};q\right)_{n}\left(q^{-N};q\right)_{n}}{\left(a q /b_{1};q\right)_{n}\left(a q / c_{1};q\right)_{n}\cdots\left(a q / b_{k};q\right)_{n}\left(a q / c_{k};q\right)_{n}\left(a q^{N+1};q\right)_{n}}\\[5pt] \nonumber
& \qquad \times\left(\frac{a^{k} q^{k+N}}{b_{1} c_{1} \cdots b_{k} c_{k}}\right)^{n} q^{-{n\choose 2}}(-1)^{n} \alpha_{n}(1,q)\\[5pt] \nonumber
&\qquad \qquad =\frac{(a q;q)_{N}\left(a q / b_{k} c_{k};q\right)_{N}}{\left(a q / b_{k};q\right)_{N}\left(a q / c_{k};q\right)_{N}} \sum_{n_{k} \geq n_{k-1} \geq \cdots \geq n_{1} \geq 0} \frac{\left(b_{k};q\right)_{n_{k}}\left(c_{k};q\right)_{n_{k}} \cdots\left(b_{1};q\right)_{n_{1}}\left(c_{1};q\right)_{n_{1}}}
{(q;q)_{n_{k}-n_{k-1}}(q;q)_{n_{k-1}-n_{k-2}} \cdots(q;q)_{n_{2}-n_{1}}}\\[5pt] \nonumber
&\qquad \qquad \qquad \times \frac{\left(q^{-N};q\right)_{n_{k}}\left(a q / b_{k-1} c_{k-1};q\right)_{n_{k}-n_{k-1}} \cdots\left(a q / b_{1} c_{1};q\right)_{n_{2}-n_{1}}}{\left(b_{k} c_{k} q^{-N} / a;q\right)_{n_{k}}\left(a q / b_{k-1};q\right)_{n_{k}}\left(a q / c_{k-1};q\right)_{n_{k}} \cdots\left(a q / b_{1};q\right)_{n_{2}}\left(a q / c_{1};q\right)_{n_{2}}}\\[5pt] \label{Multi-eq}
&\qquad \qquad \qquad \qquad \times q^{n_{1}+\cdots+n_{k}} a^{n_{1}+\cdots+n_{k-1}}\left(b_{k-1} c_{k-1}\right)^{-n_{k-1}} \cdots\left(b_{1} c_{1}\right)^{-n_{1}}
\beta_{n_{1}}(1,q).
\end{align}
\end{thm}

Below is   a limiting case of  Theorem \ref{Multi}.
\begin{prop}\label{lc}If $(\alpha_n(1,q^\eta),\beta_n(1,q^\eta))$ is a Bailey pair relative to $(1,q^\eta)$, then for $r>\lambda\geq0$,
\begin{align}\nonumber
& \sum_{n=0}^\infty\frac{2q^{(r-\frac{\lambda +1}{2})\eta n^2+\frac{\lambda +1}{2}\eta n
-(\alpha_1+\cdots+\alpha_{\lambda})n}(-q^{\alpha_1};q^{\eta})_n\cdots(-q^{
\alpha_{\lambda}};q^{\eta})_n}{(1+q^{\eta n})(-q^{\eta-\alpha_1};q^{\eta})
_n\cdots(-q^{\eta-\alpha_{\lambda}};q^{\eta})_n}\alpha_n(1,q^\eta)
\\[5pt]\nonumber
&\qquad =\frac{(q^{\eta};q^{\eta})_{\infty}}{(-q^{\eta-\alpha_1};q^{\eta})_
{\infty}}\sum_{N_1\geq N_2\geq\cdots\geq N_r\geq 0}
\frac{q^{
\eta(N_{\lambda+2}^2+\cdots+N_{r}^2)+\eta\left({N_1+1 \choose2}+\cdots+{N_{\lambda+1}+1 \choose2}\right)
-(\alpha_1N_1+\cdots+\alpha_{\lambda}N_{\lambda})
}}{(q^{\eta};q^{\eta})_{N_1-N_{2}}\cdots(q^{\eta};q^{\eta})_{
N_{r-1}-N_{r}}}\\[5pt] \label{lem-cor}
&\qquad \qquad \times\frac{(-1;q^{\eta})_{N_{\lambda+1}
}(-q^{\alpha_1};q^{\eta})_{N_1}\cdots(-q^{\alpha_{\lambda
}};q^{\eta})_{N_{\lambda}}}
{(-q^{\eta};q^{\eta})_{N_{\lambda}}
(-q^{\eta-\alpha_2};q^
{\eta})_{N_1}\cdots(-q^{\eta-\alpha_{\lambda}};q^{\eta})_{
N_{\lambda-1}}}\beta_{N_r}(1,q^\eta),
\end{align}
where we assume that $N_{r+1}=0$.
\end{prop}
\pf Replacing $q$ by $q^\eta$ and setting $k=r,\ a=1,\  c_{r-\lambda}=-1$ and $c_{r-s+1}=-q^{\alpha_s}$ for $1\leq s\leq\lambda$, as $b_i\rightarrow \infty$ for $1\leq i\leq r$, $c_m\rightarrow \infty$ for $1\leq m\leq r-\lambda-1$ and $N\rightarrow \infty$, \eqref{Multi-eq} becomes
\begin{equation*}
  \begin{split}
   & \sum_{n=0}^\infty\frac{2q^{(r-\frac{\lambda +1}{2})\eta n^2+\frac{\lambda +1}{2}\eta n
    -(\alpha_1+\cdots+\alpha_{\lambda})n}(-q^{\alpha_1};q^{\eta})_n\cdots(-q^{
\alpha_{\lambda}};q^{\eta})_n}{(1+q^{\eta n})(-q^{\eta-\alpha_1};q^{\eta})
_n\cdots(-q^{\eta-\alpha_{\lambda}};q^{\eta})_n}\alpha_n(1,q^\eta)\\[5pt]
&\qquad =\frac{(q^{\eta};q^{\eta})_{\infty}}{(-q^{\eta-\alpha_1};q^{\eta})_
{\infty}}\sum_{n_r\geq n_{r-1}\geq\cdots\geq n_1\geq 0}
\frac{q^{
\eta(n_{1}^2+\cdots+n_{r-\lambda-1}^2)+\eta\left({n_{r-\lambda}+1 \choose2}+\cdots+{n_r+1 \choose2}\right)
-(\alpha_1n_r+\cdots+\alpha_{\lambda}n_{r-\lambda+1})
}}{(q^{\eta};q^{\eta})_{n_r-n_{r-1}}\cdots(q^{\eta};q^{\eta})_{
n_{2}-n_{1}}}\\[5pt]
&\qquad \qquad \times\frac{(-1;q^{\eta})_{n_{r-\lambda}
}(-q^{\alpha_1};q^{\eta})_{n_r}\cdots(-q^{\alpha_{\lambda
}};q^{\eta})_{n_{r-\lambda+1}}}
{(-q^{\eta};q^{\eta})_{n_{r-\lambda+1}}
(-q^{\eta-\alpha_2};q^
{\eta})_{n_r}\cdots(-q^{\eta-\alpha_{\lambda}};q^{\eta})_{
n_{r-\lambda+2}}}\beta_{n_1}(1,q^\eta).
  \end{split}
\end{equation*}
Writing $n_t=N_{r+1-t}$ for $1\leq t\leq r$, we are led to \eqref{lem-cor}.
This completes the proof. \qed

The following Bailey pair is also required in the proof of  Theorem \ref{G-B-O-1}.

\begin{prop}\label{BPG} For $k\geq r \geq 1$ and $n\geq 0$,
\begin{equation}{\label{BPG-eq}}
\begin{split}
\alpha_n(1,q)&=\left\{
             \begin{array}{ll}
               1, & \hbox{if $n=0$}, \\[5pt]
              (-1)^nq^{(k-r)  n^2}(q^{(k-r-1)  n}+q^{-(k-r)  n})(1+q^{  n})/2, & \hbox{if $n\geq1$,}
             \end{array}
               \right.\\[5pt]
\beta_n(1,q)&=\sum_{n\geq N_{r+1}\geq\cdots\geq N_{k-1}\geq0}\frac{(1+q^{  n})q^{ (N_{r+1}^2+\cdots+N_{k-1}^2+N_{r+1}+\cdots+N_{k-1})
}}{2(q;q)_{n-N_{r+1}}\cdots(q;q)_
{N_{k-2}-N_{k-1}}(q^{2 };q^{2 })_{N_{k-1}}}.
\end{split}
\end{equation}
is a Bailey pair relative to $(1,q)$.
\end{prop}
\pf We begin with the following Bailey pair \cite[E(5)]{Slater-1951},
\begin{equation}\label{Bp-1}
\begin{split}
&\alpha^{(0)}_n(1,q)=\left\{
             \begin{array}{ll}
               1, & \hbox{if $n=0$}, \\[5pt]
               (-1)^n(q^{-n}+q^{n}), & \hbox{if $n\geq1$,}
             \end{array}
               \right.\\[5pt]
&\beta^{(0)}_n(1,q)=\frac{(-1)^n}{q^n(q^2;q^2)_n}.
\end{split}
\end{equation}
Applying Proposition \ref{Bl-Bp} to \eqref{Bp-1}, we obtain  that
\begin{equation*}
\begin{split}
\alpha_n^{(1)}(1,q)&=\left\{
             \begin{array}{ll}
               1, & \hbox{if $n=0$}, \\[5pt]
               2(-1)^n, & \hbox{if $n\geq1$,}
             \end{array}
               \right.\\[5pt]
\beta_n^{(1)}(1,q)&=\frac{(-1)^n}{(q^{2 };q^{2 })_n}.
\end{split}
\end{equation*}
Using Lemma \ref{BL-1}, we get
\begin{equation}{\label{unitbp-2}}
\begin{split}
\alpha_n^{(2)}(1,q)&=\left\{
             \begin{array}{ll}
               1, & \hbox{if $n=0$}, \\[5pt]
               2(-1)^nq^{n^2}, & \hbox{if $n\geq1$,}
             \end{array}
               \right.\\[5pt]
\beta_n^{(2)}(1,q)&
=\sum_{j=0}^n\frac{(-1)^jq^{j^2}}{(q;q)_{n-j}(q^{2 };q^{2 })_j}.
\end{split}
\end{equation}
Employing the following $q$-Chu-Vandermonde formula with $c=-q$ and $a\rightarrow\infty$,
\[\sum_{j=0}^n\frac{(a;q)_j(q^{-n};q)_j}{(c;q)_j(q;q)_j}
\left(\frac{cq^n}{a}\right)^j=\frac{(c/a;q)_n}{(c;q)_n},\]
we find that
\[\beta_n^{(2)}(1,q)
=\frac{1}{(q^{2};q^{2 })_n}.\]
Applying  Proposition \ref{Bl-Bp} and  Lemma \ref{BL-1}  $k-r-1$  times to \eqref{unitbp-2} yields the following
Bailey pair
\begin{equation}{\label{2k-2r-1}}
\begin{split}
\alpha_n^{(2k-2r)}(1,q)&=\left\{
             \begin{array}{ll}
               1, & \hbox{if $n=0$}, \\[5pt]
               (-1)^nq^{(k-r)  n^2}
(q^{(k-r-1)  n}+q^{-(k-r-1)  n}), & \hbox{if $n\geq1$,}
             \end{array}
               \right.\\[5pt]
\beta_n^{(2k-2r)}(1,q)&=\hskip-0.5cm\sum_{n\geq N_{r+1}\geq \cdots\geq N_{k-1}\geq0}\frac{q^{ (N_{r+1}^2+\cdots+N_{k-1}^2+N_{r+1}
+\cdots+N_{k-1})}}{(q;q)_{n-N_{r+1}}\cdots
(q;q)_{N_{k-2}-N_{k-1}}(q^{2 };q^{2 })_
{N_{k-1}}}.
\end{split}
\end{equation}
By Proposition  \ref{cbl4} and \eqref{2k-2r-1}, we obtain  the following Bailey pair
\begin{equation*}{\label{2k-2r+1}}
\begin{split}
\alpha_n(1,q)&=\left\{
             \begin{array}{ll}
               1, & \hbox{if $n=0$}, \\[5pt]
              (-1)^nq^{(k-r)  n^2}(q^{(k-r-1)  n}+q^{-(k-r)  n})(1+q^{  n})/2, & \hbox{if $n\geq1$,}
             \end{array}
               \right.\\[5pt]
\beta_n(1,q)&=\sum_{n\geq N_{r+1}\geq\cdots\geq N_{k-1}\geq0}\frac{(1+q^{  n})q^{ (N_{r+1}^2+\cdots+N_{k-1}^2+N_{r+1}+\cdots+N_{k-1})
}}{2(q;q)_{n-N_{r+1}}\cdots(q;q)_
{N_{k-2}-N_{k-1}}(q^{2 };q^{2 })_{N_{k-1}}}.
\end{split}
\end{equation*}
  This completes the proof.\qed

We conclude this section with the proof of Theorem \ref{G-B-O-1} resorting to Proposition  \ref{lc} and Proposition \ref{BPG}.

\noindent{\it Proof of Theorem \ref{G-B-O-1}. } For $k\geq r>\lambda$, plugging $\alpha_n(1,q)$ in \eqref{BPG-eq}  with $q$ replaced by $q^\eta$ into the left-hand side of \eqref{lem-cor}, and  using the assumption  that  $\alpha_i+\alpha_{\lambda+1-i}=\eta$ for $1\leq i\leq \lambda$, the  left-hand side of \eqref{lem-cor} simplifies to
\begin{align}\nonumber
   & 1+\sum_{n=1}^\infty\frac{(-q^{\alpha_1};q^{\eta})_n\cdots(-q^{
\alpha_{\lambda}};q^{\eta})_n}{(-q^{\eta-\alpha_1};q^{\eta})
_n\cdots(-q^{\eta-\alpha_{\lambda}};q^{\eta})_n}\\[5pt] \nonumber
    &\qquad \times (-1)^nq^{(k-\frac{\lambda +1}{2})\eta n^2+\frac{\lambda}{2}\eta n
-(\alpha_1+\cdots+\alpha_{\lambda})n}(q^{(k-r-\frac{1}{2})\eta n}+q^{-(k-r-\frac{1}{2})\eta n})\\[5pt]\nonumber
&\qquad \qquad=1+\sum_{n=1}^{\infty}(-1)^nq^{(k-\frac{\lambda
+1}{2})\eta n^2}(q^{(k-r-\frac{1}{2})\eta n}+q^{-(k-r-\frac{1}{2})\eta n})\\[5pt] \label{eq-L}
&\qquad \qquad =(q^{(r-\frac{\lambda}{2})\eta},q^{(2k-r-1-
\frac{\lambda}{2})\eta},q^{(2k-\lambda-1)\eta};q^{(2k-\lambda
-1)\eta})_{\infty},
\end{align}
where the last equality follows from Jacobi's triple product identity \cite[Theorem 2.8]{Andrews-1976}.

On the other hand, substituting the expression for $\beta_n(1,q)$ in \eqref{BPG-eq}  with $q$ replaced by $q^\eta$ into the right-hand side of \eqref{lem-cor}, we get
\begin{align} \nonumber
&\frac{(q^{\eta};q^{\eta})_{\infty}}{(-q^{\eta-\alpha_1};q^{\eta})_
{\infty}}\sum_{N_1\geq \cdots\geq N_{k-1}\geq0}\frac{(1+q^{-\eta N_r})(-q^{\eta};q^{\eta})_{N_{\lambda+1}
-1}q^{\eta(N_{\lambda+2}^2+\cdots+N_{k-1}^2+N_r+\cdots+N_{k-1})}}
{(q^{\eta};q^{\eta})_{N_1-N_2}\cdots(q^{\eta};q^{\eta})_{N_{
k-2}-N_{k-1}}(q^{2\eta};q^{2\eta})_{N_{k-1}}}\\[5pt] \label{eq-R1}
&\qquad \times\frac{q^{\eta\left({N_1+1 \choose2}+\cdots+{N_{\lambda+1}+1 \choose2}\right)-(\alpha_1N_1+\cdots+\alpha_{\lambda}N_{\lambda})}
(-q^{\alpha_1};q^{\eta})_{N_1}\cdots(-q^{\alpha_{\lambda
}};q^{\eta})_{N_{\lambda}}}{(-q^{\eta};q^{\eta})_{
N_{\lambda}}(-q^{\eta-\alpha_2};q^{\eta})_{N_1}
\cdots(-q^{\eta-\alpha
_{\lambda}};q^{\eta})_{N_{\lambda-1}}}.
\end{align}
Observing that
\[
  (-q^{r};q^\eta)_{n}=q^{r n+\eta {n\choose 2}}(-q^{\eta-r-n\eta};q^\eta)_{n},\]
  and
  \[
  \frac{1}{(-q^{\eta-r};q^\eta)_n} =\frac{(-q^{\eta-r+n\eta};q^\eta)_\infty}
  {(-q^{\eta-r};q^\eta)_\infty},
\]
  the summation in \eqref{eq-R1} equals
\begin{align}\nonumber
&\sum_{N_1\geq\cdots\geq N_{k-1}\geq0}\frac{q^{\eta(N_1^2+\cdots+N_{k-1}^2+N_r+\cdots+N_{k-1})}
(1+q^{-\eta N_r})(-q^{\eta-\eta N_{\lambda+1}};q^{\eta})_{N_{\lambda+1}-1}}
{(q^{\eta};q^{\eta})_{N_1-N_2}\cdots(q^{\eta};q^{\eta})_
{N_{k-2}-N_{k-1}}(q^{2\eta};q^{2\eta})_{N_{k-1}}}\\[5pt] \label{eq-R2}
&\qquad \times\frac{(-q^{\eta+\eta N_{\lambda}};q^{\eta})_{\infty}\prod_{s=1}^{\lambda}(-q^{\eta-\alpha_s
-\eta N_{s}};q^{\eta})_{N_s}\prod_{s=2}^{\lambda}(-q^{\eta-\alpha_s+
\eta N_{s-1}};q^{\eta})_{\infty}}{(-q^\eta;q^\eta)_\infty\prod_{s=2}^{\lambda}(-q^{
\eta-\alpha_s};q^{\eta})_{\infty}}.
\end{align}
Combining \eqref{eq-L}  and \eqref{eq-R2}, we deduce that
\begin{align}\nonumber
&\frac{(q^{\eta};q^{\eta})_{\infty}}{(-q^{\eta-\alpha_1};q^{\eta})_
{\infty}}\sum_{N_1\geq\cdots\geq N_{k-1}\geq0}\frac{q^{\eta(N_1^2+\cdots+N_{k-1}^2+N_r+\cdots+N_{k-1})}
(1+q^{-\eta N_r})(-q^{\eta-\eta N_{\lambda+1}};q^{\eta})_{N_{\lambda+1}-1}}
{(q^{\eta};q^{\eta})_{N_1-N_2}\cdots(q^{\eta};q^{\eta})_
{N_{k-2}-N_{k-1}}(q^{2\eta};q^{2\eta})_{N_{k-1}}}\\[5pt] \nonumber
&\qquad\times\frac{(-q^{\eta+\eta N_{\lambda}};q^{\eta})_{\infty}\prod_{s=1}^{\lambda}(-q^{\eta-\alpha_s
-\eta N_{s}};q^{\eta})_{N_s}\prod_{s=2}^{\lambda}(-q^{\eta-\alpha_s+
\eta N_{s-1}};q^{\eta})_{\infty}}{(-q^\eta;q^\eta)_\infty\prod_{s=2}^{\lambda}(-q^{
\eta-\alpha_s};q^{\eta})_{\infty}}\\[5pt] \nonumber
&\qquad \qquad=(q^{(r-\frac{\lambda}{2})\eta},q^{(2k-r-1-
\frac{\lambda}{2})\eta},q^{(2k-\lambda-1)\eta};q^{(2k-\lambda
-1)\eta})_{\infty}.
\end{align}
Multiplying  both sides  by $$\frac{(-q^{\eta-\alpha_1},\ldots,-q^{\eta-\alpha_{\lambda}},
-q^\eta;q^\eta)_\infty}{(q^{\eta};q^{\eta})_{\infty}},$$
we obtain
\begin{align}\nonumber
& \sum_{N_1\geq\cdots\geq N_{k-1}\geq0}\frac{q^{\eta(N_1^2+\cdots+N_{k-1}^2+N_r+\cdots+N_{k-1})}
(1+q^{-\eta N_r})(-q^{\eta-\eta N_{\lambda+1}};q^{\eta})_{N_{\lambda+1}-1}(-q^{\eta+\eta N_{\lambda}};q^{\eta})_{\infty}}
{(q^{\eta};q^{\eta})_{N_1-N_2}\cdots(q^{\eta};q^{\eta})_
{N_{k-2}-N_{k-1}}(q^{2\eta};q^{2\eta})_{N_{k-1}}}\\[5pt] \nonumber
&\qquad\qquad\qquad\times\prod_{s=1}^{\lambda}(-q^{\eta-\alpha_s
-\eta N_{s}};q^{\eta})_{N_s}\prod_{s=2}^{\lambda}(-q^{\eta-\alpha_s+
\eta N_{s-1}};q^{\eta})_{\infty}\\[5pt] \nonumber
&\qquad \qquad=\frac{(-q^{\eta-\alpha_1},\ldots,-q^{\eta-\alpha_{\lambda}},
-q^\eta;q^\eta)_\infty(q^{(r-\frac{\lambda}{2})\eta},q^{(2k-r-1-
\frac{\lambda}{2})\eta},q^{(2k-\lambda-1)\eta};q^{(2k-\lambda
-1)\eta})_{\infty}}{(q^{\eta};q^{\eta})_{\infty}}.
\end{align}
But $\alpha_i+\alpha_{\lambda+1-i}=\eta$  for $1\leq i\leq\lambda$,
so we reach \eqref{G-B-O-1-eq} in Theorem \ref{G-B-O-1}.
This completes the proof. \qed

\section{Concluding remarks}

To conclude, we make a few remarks on the connection between the
main results of this paper and the original conjecture of Bressoud,
along with our subsequent work in this direction. Then we
mention some potential problems for future study.

It should be stressed that the overpartition analogues
considered in this paper are not merely a matter of
extension and specialization.
In fact, they play an essential role and serve as an indispensable
structure in tackling the conjecture of Bressoud formulated  in terms of ordinary partitions.

Based on the relationship between the overpartition analogue  $\overline{B}_1$ and Bressoud's function ${B}_{0}$ (Theorem \ref{rel-over1}),  we realize that the  case $j=0$ of Bressoud's conjecture    (that is, $A_0=B_0$) is a consequence of the relation
$\overline{A}_1=\overline{B}_1$ on overpartitions.
Nevertheless, the  case $j=1$ of Bressoud's conjecture  has been resolved by Kim \cite{Kim-2018} without resorting to overpartitions.
One is immediately led to show that $\overline{A}_1=\overline{B}_1$.
This is the objective of our subsequent paper \cite{He-Ji-Zhao}.
It is worth mentioning that the relation
$\overline{A}_1=\overline{B}_1$ can be regarded as
an overpartition analogue of Bressoud's conjecture for the case $j=1$.
In other words, we may say that Bressoud's conjecture
consists of two parts, one of which is the case $j=1$  settled by
Kim, and the other (the  case $j=0$) is an overpartition analogue.
Naturally, it would be interesting to give a direct combinatorial proof of the case  $j=0$ of Bressoud's conjecture  without relying
on the overpartition setting. Also, it would be desirable to give  direct combinatorial proofs of the generating functions of $\overline{B}_0$ and $\overline{B}_1$.

 \vskip 0.2cm
\noindent{\bf Acknowledgment.} This work
was supported by   the National Science Foundation of China. We are greatly indebted to the referees for their insightful suggestions leading to an improvement of an earlier version.

\end{document}